\documentclass[twoside,11pt]{article}

\usepackage{jmlr2e}
\usepackage{amssymb,latexsym,amsmath,bm,amsthm}   
\usepackage{booktabs}  
\usepackage{graphicx}
\usepackage{subcaption}
\usepackage{mathtools}
\usepackage{xcolor}
\usepackage[title]{appendix}
\usepackage{url}
\usepackage{hyperref}
	
\renewcommand{\limsup}{\varlimsup}

\newcommand{\mc}[1]{\mathcal{#1}}

\theoremstyle{lemma}
\newtheorem{definition}{Definition}
\theoremstyle{lemma}
\newtheorem{lemma}{Lemma}
\newtheorem{remark}{Remark}
\theoremstyle{plain}
\newtheorem{theorem}{Theorem}
\theoremstyle{plain}

\theoremstyle{plain}

\DeclareMathOperator*{\argmin}{arg\,min}

\usepackage{xargs}                      

\usepackage[colorinlistoftodos,prependcaption,textsize=tiny]{todonotes}
\newcommandx{\unsure}[2][1=]{\todo[linecolor=red,backgroundcolor=red!25,bordercolor=red,#1]{#2}}
\newcommandx{\change}[2][1=]{\todo[linecolor=blue,backgroundcolor=blue!25,bordercolor=blue,#1]{#2}}
\newcommandx{\info}[2][1=]{\todo[linecolor=OliveGreen,backgroundcolor=OliveGreen!25,bordercolor=OliveGreen,#1]{#2}}
\newcommandx{\improvement}[2][1=]{\todo[linecolor=Plum,backgroundcolor=Plum!25,bordercolor=Plum,#1]{#2}}
\newcommandx{\thiswillnotshow}[2][1=]{\todo[disable,#1]{#2}}

\ShortHeadings{EB Estimation in Discrete Linear Exponential Family}{Banerjee, Liu, Mukherjee and Sun}
\firstpageno{1}

\begin{document}

\title{A General Framework for Empirical Bayes Estimation in Discrete Linear Exponential Family}

\author{\name Trambak Banerjee  \email trambakb@usc.edu\\
	\addr Data Sciences and Operations\\
	University of Southern California\\
	Los Angeles, CA 90089, USA
	 \AND 
	 \name Qiang Liu \email lqiang@cs.utexas.edu\\
	 \addr Computer Science\\
	 University of Texas at Austin\\
	 Austin, Texas 78712, USA
	 \AND 
	 \name Gourab Mukherjee \email gourab@usc.edu\\
	 \addr Data Sciences and Operations\\
	 University of Southern California\\
	 Los Angeles, CA 90089, USA
	  \AND 
	  \name Wenguang Sun \email wenguans@marshall.usc.edu\\
       \addr Data Sciences and Operations\\
       University of Southern California\\
       Los Angeles, CA 90089, USA}

\editor{}
\maketitle
\begin{abstract}
We develop a Nonparametric Empirical Bayes (\texttt{NEB}) framework for compound estimation in the discrete linear exponential family, which includes a wide class of discrete distributions frequently arising from modern big data applications. We propose to directly estimate the Bayes shrinkage factor in the generalized Robbins' formula via solving a scalable convex program, which is carefully developed based on a RKHS representation of the Stein's discrepancy measure. The new \texttt{NEB} estimation framework is flexible for incorporating various structural constraints into the data driven rule, and provides a unified approach to compound estimation with both regular and scaled squared error losses. We develop theory to show that the class of \texttt{NEB} estimators enjoys strong asymptotic properties. Comprehensive simulation studies as well as analyses of real data examples are carried out to demonstrate the superiority of the \texttt{NEB} estimator over competing methods.

\end{abstract}

\medskip


%
%
%
%
%
%
%
%
\newpage
%
\section{Introduction}\label{sec:intro}
%

Shrinkage methods, exemplified by the seminal work of \citet{james1961estimation}, have received renewed attention  in modern large-scale inference problems (\citealp{efron2012large,fourdrinier2018shrinkage}). 
Under this setting, the classical Normal means problem has been extensively studied (\citealp{brown2008season,jiang2009general,brown2009nonparametric,efron2011tweedie,xie2012sure,weinstein2018group}). However, in a variety of applications, the observed data are often discrete. For instance, in the News Popularity study discussed in Section \ref{sec:realdata}, the goal is to estimate the popularity of a large number of news items based on their frequencies of being shared in social media platforms such as Facebook and LinkedIn. Another important application scenario arises from genomics research, where estimating the expected number of mutations across a large number of genomic locations can help identify key drivers or inhibitors of a given phenotype of interest. 

We mention two main limitations of existing shrinkage estimation methods. First, the methodology and theory developed for continuous variables, in particular for Normal means problem, may not be directly applicable to discrete models. Second, existing methods have focused on the squared error loss. However, the scaled loss \citep{clevenson1975simultaneous}, which effectively reflects the asymmetries in decision making [cf. Equation \eqref{scaled-loss}], is a more desirable choice for many discrete models such as Poisson, where the scaled loss corresponds to the local Kulback-Leibler distance. The scaled loss also provides a more desirable criterion in a range of sparse settings, for example, when the goal is to estimate the rates of rare outcomes in Binomial distributions \citep{fourdrinier1995intrinsic}. Much research is needed for discrete estimation problems under various loss functions. This article develops a general framework for empirical Bayes estimation for the discrete linear exponential (\texttt{DLE}) family, also known as the family of discrete power series distributions \citep{noack1950class}, under both regular and scaled error losses.

The \texttt{DLE} family includes a wide class of popular members such as the Poisson, Binomial, negative Binomial and Geometric distributions. Let $Y$ be a non-negative integer valued random variable. Then $Y$ is said to belong to a \texttt{DLE} family if its probability mass function (pmf) is of the form 
\begin{equation}
\label{eq:DLE}
p(y|\theta)=\dfrac{a_y\theta^y}{g(\theta)},~~y\in\{0, 1, 2, \cdots\},
\end{equation}
where $a_y$ and $g(\theta)$ are known functions such that $a_y\ge 0$ is independent of $\theta$ and $g(\theta)$ is a normalizing factor that is differentiable at every $\theta$. Special cases of \texttt{DLE} include the $\mbox{Poisson}(\lambda)$ distribution with $a_y=(y!)^{-1}$, $\theta=\lambda$ and $g(\theta)=\exp{(\theta)}$, and the $\mbox{Binomial}(m, q)$ distribution  with $a_y=\binom{m}{y}$, $\theta=q/(1-q)$ and $g(\theta)=(1+\theta)^m$. 

Suppose $Y_1, \ldots, Y_n$ obey the following hierarchical model 
\begin{equation}\label{eq:DLEmodel}
Y_i~|~\theta_i\stackrel{ind.}{\sim}\texttt{DLE}(\theta_i),~~\theta_i\stackrel{i.i.d}{\sim}G(\cdot),
\end{equation}
where $G(\cdot)$ is an unspecified prior distributusuaion on $\theta_i$. The problem of interest is to estimate $\bm{\theta}=(\theta_1,\ldots,\theta_n)$ based on $\bm Y=(Y_1,\ldots, Y_n)$. 
Empirical Bayes approaches to this compound decision problem date back to the famous Robbins' formula \citep{robbins1956empirical} under the Poisson model. Important recent progresses by \cite{brown2013poisson}, \cite{koenker2014convex} and \cite{koenker2017rebayes} show that Robbins' estimator can be vastly improved by incorporating smoothness and monotonicity adjustments. The main idea of existing works is to approximate the shrinkage factor in the Bayes estimator as smooth functionals of the unknown marginal pmf $p(y)$. The pmf can be estimated in various ways including the observed empirical frequencies \citep{robbins1956empirical}, the smoothness-adjusted estimator \citep{brown2013poisson} or the shape-constrained NPMLE approach \citep{koenker2014convex,koenker2017rebayes}. 

This article develops a general non-parametric empirical Bayes (\texttt{NEB}) framework for compound estimation in discrete models. We first derive generalized Robbins' formula (\texttt{GRF}) for the \texttt{DLE} Model \eqref{eq:DLEmodel}, and then implement \texttt{GRF} via solving a scalable convex program. The powerful convex program, which is carefully developed based on a reproducing kernel Hilbert space (RKHS) representation of Stein's discrepancy measure, leads to a class of efficient \texttt{NEB} shrinkage estimators. We develop theories to show that the \texttt{NEB} estimator is $\sqrt{n}$ consistent up to certain logarithmic factors and enjoys superior risk properties. Simulation studies are conducted to illustrate the superiority of the proposed \texttt{NEB} estimator when compared to existing state-of-the-art approaches such as \cite{brown2013poisson}, \cite{koenker2014convex} and \cite{koenker2017rebayes}. We show that the \texttt{NEB} estimator has smaller risk in all comparisons and the efficiency gain is substantial in many settings.

There are several advantages of the proposed \texttt{NEB} estimation framework. First, in contrast with existing methods such as the smoothness-adjusted Poisson estimator in \cite{brown2013poisson}, our methodology covers a much wider range of distributions and presents a unified approach to compound estimation in discrete models. Second, our proposed convex program is fast and scalable. It directly produces stable estimates of optimal Bayes shrinkage factors and can easily incorporate various structural constraints into the decision rule. By contrast, the three-step estimator in \cite{brown2013poisson}, which involves smoothing, Rao-Blackwellization and monotonicity adjustments, is complicated, computationally intensive and sometimes unstable (as the numerator and denominator of the ratio are computed separately). {Third, the RKHS representation of Stein's discrepancy measure provides a new analytical tool for developing theories such as asymptotic optimality and convergence rates. }
Finally, the \texttt{NEB} estimation framework is robust to departures from the true model due to its utilization of a generic quadratic program that does not rely on the specific form of a particular \texttt{DLE} family. Our numerical results in Section \ref{sec:numresults} demonstrate that the \texttt{NEB} estimator has significantly better risk performance than competitive approaches of \cite{efron2011tweedie}, \cite{brown2013poisson} and \cite{koenker2017rebayes} under a mis-specified Poisson model. 

An alternative approach to compound estimation in discrete models, as suggested and investigated by \cite{brown2013poisson}, is to employ variance stabilizing transformations, which converts the discrete problem to a classical normal means problem. This allows estimation via Tweedie's formula for normal variables \citep{efron2011tweedie}, where the marginal density can be estimated using NPMLE \citep{jiang2009general, koenker2014convex} or through kernel density methods \citep{brown2009nonparametric}. However, there are several drawbacks of this approach compared to our \texttt{NEB} framework. First, Tweedie's formula is not applicable to scaled error loss whereas our methodology is built upon the generalized Robbins' formula, which covers both regular and scaled squared error losses.  Second, there can be information loss in conventional data processing steps such as standardization, transformation and continuity approximation. While investigating the impact of information loss on compound estimation is of great interest, it is desirable to develop methodologies directly based on generalized Robbins' formula that is specifically derived and tailored for discrete variables. Finally, our \texttt{NEB} framework provides a convenient tool for developing asymptotic theories. By contrast, convergence rates are yet to be developed for normality inducing transformations, which can be highly non-trivial.  

The rest of the paper is organized as follows. In Section \ref{sec:2}, we introduce our estimation framework while Section \ref{sec:theory} presents a theoretical analysis of the \texttt{NEB} estimator. The numerical performance of our method is investigated using both simulated and real data in Sections \ref{sec:numresults} and \ref{sec:realdata} respectively. Additional technical details and proofs are relegated to the Appendices.

%
\section{A General Framework for Compound Estimation in \texttt{DLE} Family}\label{sec:2}
%

This section describes the proposed \texttt{NEB} framework for compound estimation in discrete models. We first introduce in Section \ref{sec:gen_robins_form} the generalized Robbins' formula for the \texttt{DLE} family \eqref{eq:DLEmodel}, then propose in Section \ref{sec:min_dkskd} a convex optimization approach for its practical implementation. Details for tuning parameter selection are discussed in Section \ref{sec:bandwidth}. 

\subsection{Generalized Robbins' formula for \texttt{DLE} models}\label{sec:gen_robins_form}

Denote $\delta_i$ an estimator of $\theta_i$. Consider a class of loss functions 
\begin{equation}\label{scaled-loss}
\ell^{(k)}(\theta_i,\delta_i)=\theta_i^{-k}(\theta_i-\delta_i)^2
\end{equation}
for $k\in \{0, 1\}$, where $\ell^{(0)}(\theta_i,\delta_i)$ is the usual squared error loss, and $\ell^{(1)}(\theta_i,\delta_i)=\theta_i^{-1}(\delta_i-\theta_i)^2$  corresponds to the scaled squared error loss \citep{clevenson1975simultaneous,fourdrinier1995intrinsic}. 
In compound estimation, one is concerned with the average loss 
$$
\mathcal{L}_n^{(k)}(\bm \theta,\bm \delta)=n^{-1}\sum_{i=1}^{n}\ell^{(k)}(\theta_i,\delta_i).
$$ 
The associated risk is denoted $R_n^{(k)}(\bm \theta,\bm \delta)=\mathbb{E_{\pmb Y|\pmb\theta}}\mathcal{L}_n^{(k)}(\bm \theta,\bm \delta)$. {Let $\pmb G(\pmb\theta)$ denote the joint distribution of $(\theta_1, \cdots, \theta_n)$.} The Bayes estimator $\bm \delta^{\pi}_{(k)}$ that minimizes the Bayes risk $B_n^{(k)}(\bm \theta) = \int R_n^{(k)}(\bm \theta,\bm \delta)d\pmb G(\bm \theta)$ is given by Lemma \ref{lem:bayes}.

\begin{lemma}[Generalized Robbins' formula]
	\label{lem:bayes}
	Consider the \texttt{DLE} Model \eqref{eq:DLEmodel}. Let $p(\cdot)=\int p(\cdot|\theta)dG(\theta)$ be the marginal pmf of $Y$. Define for $k\in\{0,1\}$,
$$
	w_p^{(k)}(y_i)=
	\dfrac{p(y_i-k)}{p(y_i+1-k)},\text{ for }y_i=k, k+1, \cdots.
$$
Then the Bayes estimator that minimizes the risk $B_n^{(k)}(\bm \theta)$ is given by $\bm \delta^{\pi}_{(k)}=\{\delta_{(k),i}^{\pi}(y_i): 1\leq i\leq n\}$, where 
\begin{equation}\label{wk}
	\delta_{(k),i}^{\pi}(y_i)=\begin{cases}
	\dfrac{a_{y_i-k}/a_{y_i+1-k}}{w_p^{(k)}(y_i)},\text{ for }y_i=k,k+1,\cdots\\
	0, \quad\quad\quad\quad\quad \quad \text{ for }y_i<k
	\end{cases}.
\end{equation}
\end{lemma}
\begin{remark}\rm{
Under the squared error loss $(k=0)$ with $Y_i~|~\theta_i\sim\text{Poi}(\theta_i)$ and $a_{y_i}=(y_i!)^{-1}$, Lemma \ref{lem:bayes} yields
\begin{equation}\label{Robbins-formula}
\delta_{(0), i}^{\pi}(y_i)=(y_i+1)\dfrac{p(y_i+1)}{p(y_i)},
\end{equation}
which recovers the classical Robbins' formula \citep{robbins1956empirical}. {In contrast, under the scaled loss, we have
\begin{equation}\label{Robbins-formula-scaled}
\delta_{(1), i}^{\pi}(y_i)=y_i\dfrac{p(y_i)}{p(y_i-1)} \mbox{ for } y_i>0 \mbox{ and } \delta_{(1), i}^{\pi}(y_i)=0 \mbox{ otherwise}.
\end{equation}
Under scaled error loss the estimator \eqref{Robbins-formula} can be much outperformed by \eqref{Robbins-formula-scaled} (and vice versa under the regular loss). We develop parallel results for the two types of loss functions. }
}
\end{remark}

Next we discuss related works for implementing Robbins' formula under the empirical Bayes (EB) estimation framework. Inspecting \eqref{wk} and \eqref{Robbins-formula}, we can view ${a_{y_i-k}/a_{y_i+1-k}}$ as a naive and known estimator of $\theta_i$. The ratio functional $w_p^{(k)}(y_i)$, which is unknown in practice, represents the optimal shrinkage factor that depends on $p(\cdot)$. Hence, a simple EB approach, as done in the classical Robbins' formula, is to estimate $w^{(k)}_p(y)$ by plugging-in empirical frequencies: $\hat w_n^{(0)}(y)=\hat p_n(y)/\hat p_n(y+1)$, where $\hat{p}_n(y)=n^{-1}\sum_{i=1}^{n}\mathbb I(y_i=y)$. It is noted by \cite{brown2013poisson} that this plug-in estimator can be highly inefficient especially when $\theta_i$ are small. Moreover, the numerator and denominator in $w_p^{(0)}(y)$ are estimated separately, which may lead to unstable ratios. \cite{brown2013poisson} showed that Robbins' formula can be dramatically improved by imposing additional smoothness and monotonicity adjustments. An alternative approach is to estimate $p(y)$ using NPMLE \citep{jiang2009general} under appropriate shape constraints \citep{koenker2014convex}. However, efficient estimation of $p(y)$ may not directly translate into an efficient estimation of the underlying ratio $p(y+1)/p(y)$. We recast the compound estimation problems as a convex program, which directly produces consistent estimates of the ratio functionals 
$$\bm w_p^{(k)}=\left\{w_p^{(k)}(y_1),\ldots,w_p^{(k)}(y_n)\right\}$$ from data. The estimators are shown to enjoy superior numerical and theoretical properties. Unlike existing works that are limited to regular loss and specific members in the \texttt{DLE} family, our method can handle a wide range of distributions and various types of loss functions in a unified framework. 

\subsection{Shrinkage estimation by convex optimization}\label{sec:min_dkskd}

This section focuses on the scaled squared error loss $(k=1)$. Methodologies and theories for the case with usual squared error loss $(k=0)$ can be derived similarly;  details are provided in Appendix \ref{sec:neb_sqaurederror}. We first introduce some notations and then present the \texttt{NEB} estimator in Definition \ref{def:neb_k1}.

Suppose $Y$ is a non-negative integer-valued random variable with pmf $p(\cdot)$. Define
\begin{equation}
\label{eq:h(y)}
h_0^{(1)}(y)=\begin{cases}
1~, \quad\quad\quad\quad\quad \text{if~}y=0 \\
1-w_p^{(1)}(y), \quad \text{if~}  y\in\{1, 2, \ldots\}. 
\end{cases}
\end{equation}
Let $\mc{K}_\lambda(y, y')=\exp\{-\frac 1 {2\lambda}(y-y')^2\}$ be the positive definite Radial Basis Function (RBF) kernel with bandwidth parameter $\lambda\in\Lambda$ where $\Lambda$ is a compact subset of $\mathbb{R}^+$ bounded away from $0$. 
Given observations $\bm y=(y_1,\ldots,y_n)$ from Model \eqref{eq:DLEmodel}, let $\bm h_0^{(1)} = \left\{h_{0}^{(1)}(y_1),\ldots,h_{0}^{(1)}(y_n)\right\}$. Define operators $\Delta_{y} \mc{K}_\lambda(y,y')=\mc{K}_\lambda(y+1,y')-\mc{K}_\lambda(y,y')$ and $$\Delta_{y,y'} \mc{K}_\lambda(y,y')=\Delta_{y'}\Delta_{y} \mc{K}_\lambda(y,y')=\Delta_{y}\Delta_{y'} \mc{K}_\lambda(y,y').$$ 
Consider the following $n\times n$ matrices, which are needed in the definition of the \texttt{NEB} estimator: 
$$
\bm K_{\lambda} = {n^{-2}}[\mc{K}_\lambda(y_i,y_j)]_{ij}, 
\quad 
{\Delta \bm K}_{\lambda}  = {n^{-2}}[\Delta_{y_i} \mc{K}_\lambda(y_i, y_j)]_{ij}, 
\quad 
{\Delta_2 \bm K}_{\lambda}  = {n^{-2}}[\Delta_{y_i,y_j}\mc{K}_\lambda(y_i, y_j)]_{ij}.
$$

\renewcommand*{\thedefinition}{1A} 
\begin{definition}[\texttt{NEB} estimator] 
	\label{def:neb_k1}	
Consider the \texttt{DLE} Model \eqref{eq:DLEmodel} with loss $\ell^{(1)}(\theta_i,\delta_i)$. For any fixed $\lambda\in\Lambda$, let $\hat{\bm h}_n^{(1)}(\lambda)=\left\{\hat{h}_{1}^{(1)}(\lambda),\ldots,\hat{h}_{n}^{(1)}(\lambda)\right\}$ be the solution to the following quadratic optimization problem:
	\begin{eqnarray}
	\label{eq:quad_opt}
	&&\min_{\bm h\in \bm H_n}~\hat{\mathbb{M}}_{\lambda,n}(\bm h)=\bm h^T\bm K_\lambda\bm h+2\bm h^T\Delta \bm K_\lambda\bm 1+\bm 1^T\Delta_2\bm K_\lambda\bm 1,
	\end{eqnarray}	
where $\bm H_n=\{\bm h=(h_1,\ldots,h_n):\mc{A}\bm h\preceq \bm b,~\mc{C}\bm h= \bm d\}$ is a convex set and $\mc{A},\mc{C},\bm b$ and $\bm d$ are known real matrices and vectors that enforce linear constraints on the components of $\bm h$. Define $\hat{w}_{i}^{(1)}(\lambda)=1-\hat{h}_{i}^{(1)}(\lambda)$. Then the \texttt{NEB} estimator is given by ${\bm\delta}^{\sf neb}_{(1)} (\lambda)=\left\{{\delta}_{(1),i}^{\sf neb}(\lambda): 1\leq i\leq n\right\}$, where
$$
{\delta}_{(1),i}^{\sf neb}(\lambda)=
	\dfrac{a_{y_i-1}/a_{y_i}}{\hat{w}_{i}^{(1)}(\lambda)}, \quad \mbox{ if } y_i\in\{1, 2, \ldots\}, 
$$
and ${\delta}_{(1),i}^{\sf neb}(\lambda)=0$ if $y_i=0$. 

\end{definition}

Next we provide some insights on why the optimization criterion \eqref{eq:quad_opt} works; theories are developed in Section \ref{sec:theory} to establish the properties of the \texttt{NEB} estimator rigorously. Denote $h_0^{(1)}$ and $\tilde{h}^{(1)}$ as the ratio functionals corresponding to pmfs $p$ and $\tilde{p}$, respectively. Suppose $Y_i$ are i.i.d. samples obeying $p(y)$. Theorem \ref{thm:spq} shows that 
$$
\hat{\mathbb{M}}_{\lambda,n}(\tilde{\bm h})={\mathbb{M}}_{\lambda}(\tilde{\bm h})+O_p\Big(\dfrac{\log^2n}{n^{1/2}}\Big),
$$
where $\hat{\mathbb{M}}_{\lambda,n}(\tilde{\bm h})$ is the objective function in \eqref{eq:quad_opt} and ${\mathbb{M}}_{\lambda}(\tilde{\bm h})$, also denoted $\mc{S}_\lambda[\tilde{p}](p)$, is the kernelized Stein's discrepancy (KSD). Roughly speaking, the KSD measures how different
one distribution $p$ is from another distribution $\tilde p$, with $\mc{S}_\lambda[\tilde{p}](p)=0$ if and only if $\tilde{p}=p$. A key feature of the KSD is that $\mc{S}_\lambda[\tilde{p}](p)$ can be equivalently represented by the discrepancy between the corresponding ratio functionals $h_0^{(1)}$ and $\tilde{h}^{(1)}$. Hence, optimizing \eqref{eq:quad_opt} is asymptotically equivalent to finding $\tilde{h}^{(1)}$ that is as close as possible to the true underlying $h_0^{(1)}$, which corresponds to the optimal shrinkage factor in the compound estimation problem. Theorems \ref{thm:w} and \ref{thm:bayesrisk_k1} demonstrate that \eqref{eq:quad_opt} is an effective convex program in the sense that the minimizer $\hat{\bm{h}}_n$ is $\sqrt{n}$ consistent with respect to $\bm{h}_0^{(1)}$, and the resultant \texttt{NEB} estimator converges to the Bayes estimator. 

\subsection{Structural constraints and bandwidth selection}\label{sec:bandwidth}

In problem \eqref{eq:quad_opt} the linear inequality $\mc{A}\bm h\preceq \bm b$ can be used to impose structural constraints on the \texttt{NEB} rule ${\bm\delta}^{\sf neb}_{(1)} (\lambda)$. The structural constraints, which may take the form of monotonicity constraints as pursued in, for example, \cite{brown2013poisson} and \cite{koenker2014convex}, have been shown to be effective for stabilizing the estimator and hence improving the accuracy. For example, when $Y_i~|~\theta_i~\sim~\text{Poi}(\theta_i)$ then $\delta_{(1),i}^{\pi}(y_i)=(y_i+1)/w_p^{(1)}(y_i)$ and $\mc{A}$, $\bm b$ can be chosen such that 
$$
h_{(i-1)}-\frac{y_{(i-1)}+1}{y_{(i)}+1} h_{(i)}\le 1-\frac{y_{(i-1)}+1}{y_{(i)}+1}, \text{ for }2\le i\le n
$$ 
and $y_{(1)}\le y_{(2)}\le \cdots \le y_{(n)}$. Moreover, when $y_i=0$ we set $\delta_{(1),i}^{\sf neb}(\lambda)=0$ by convention (see lemma \ref{lem:bayes}). The equality constraints $\mc{C}\bm h= \bm d$ accommodate such boundary conditions along with instances of ties for which we require $h_i=h_j$ whenever $y_i=y_j$.

The implementation of the quadratic program in \eqref{eq:quad_opt} requires the choice of a tuning parameter $\lambda$ in the RBF kernel. For practical applications, $\lambda$ must be determined in a data-driven fashion. For infinitely divisible random variables \citep{Klenke2014} such as Poisson variables, \cite{brown2013poisson} proposed a modified cross validation (MCV) method for choosing the tuning parameter. 
However, the MCV method cannot be applied to distributions with bounded support, e.g. variables that are not infinitely divisible \citep{sato1999levy} such as the Binomial distribution. {To provide a unified estimation framework for the \texttt{DLE} family, we develop an alternative method for choosing $\lambda$. The key idea is to derive an asymptotic risk estimate ${\sf ARE}^{(1)}_n(\lambda)$ that serves as an approximation of the true risk $\mc{R}_n^{(1)}(\bm \theta,\cdot)$. Then the tuning parameter is chosen to minimize ${\sf ARE}^{(1)}_n(\lambda)$.} 

The methodology based on ARE is illustrated below for Poisson and Binomial models under the scaled loss (see Definitions 2A and 3A, respectively). The ideas can be extended to other members in the \texttt{DLE} family. In Appendix \ref{sec:lam_k0}, we provide relevant details for choosing $\lambda$ under the regular loss $\mc{L}_n^{(0)}$.

\renewcommand*{\thedefinition}{2A}
\begin{definition}[ARE of ${\bm \delta}_{(1)}^{\sf neb}(\lambda)$ in the Poisson model]
	\label{def:are_pois_k1}	
	Suppose $Y_i~|~\theta_i\stackrel{ind.}{\sim} \texttt{Poi}(\theta_i)$. Under the loss $\ell^{(1)}(\theta_i,\cdot)$, an ARE of the true risk of ${\bm \delta}_{(1)}^{\sf neb}(\lambda)$ is
	$${\sf ARE}^{(1,\mc{P})}_n(\lambda)=\dfrac{1}{n}\Bigl\{ \sum_{i=1}^{n}y_i+\sum_{i=1}^{n}{\psi}_\lambda(y_i)-2\sum_{i=1}^{n}{\delta}_{(1),i}^{\sf neb}(\lambda)\Bigr\}, \mbox{ where}
	$$
 $$
\psi_\lambda(y_i)=\{{\delta}_{(1),j}^{\sf neb}(\lambda)\}^2/(y_i+1),~y_i= 0,1,\ldots.
$$ with $j\in\{1,\ldots,n\}$ such that $y_j=y_i+1$.
\end{definition}
For the Binomial model, we proceed along similar lines and consider the following asymptotic risk estimate of the true risk.

\renewcommand*{\thedefinition}{3A}
\begin{definition}[ARE of ${\bm \delta}_{(1)}^{\sf neb}(\lambda)$ in the Binomial model]
	\label{def:are_bin_k1}	
	Suppose $Y_i~|~q_i\sim \mc{B}in(m,q_i)$. Hence in Equations \eqref{eq:DLE} and \eqref{eq:DLEmodel} we have $a_{y_i}={m \choose y_i}$ and $\theta_i=q_i/(1-q_i)$. Under the loss $\ell^{(1)}(\theta_i,\cdot)$, an ARE of the true risk of ${\bm \delta}_{(1)}^{\sf neb}(\lambda)$ is
	$${\sf ARE}^{(1,\mc{B})}_n(\lambda)=\dfrac{1}{n}\Bigl\{ \sum_{i=1}^{n}\dfrac{y_i}{m-y_i+1}+\sum_{i=1}^{n}(m-y_i){\psi}_\lambda(y_i)-2\sum_{i=1}^{n}{\delta}_{(1),i}^{\sf neb}(\lambda)\Bigr\}, \mbox{ where}
	$$	 
	$$
	\psi_\lambda(y_i)=\{{\delta}_{(1),j}^{\sf neb}(\lambda)\}^2/(y_i+1),~y_i= 0,\ldots,m.\\
	$$ with $j\in\{1,\ldots,n\}$ such that $y_j=y_i+1$.
\end{definition}

\begin{remark}\rm{
Although the expression for $\psi_\lambda$ appears to be identical across Definitions \ref{def:are_pois_k1} and \ref{def:are_bin_k1}, it differs with respect to ${\bm \delta}_{(1)}^{\sf neb}(\lambda)$. Specifically, in definition \ref{def:are_pois_k1}, ${\bm \delta}_{(1)}^{\sf neb}(\lambda)$ is the \texttt{NEB} estimator of the Poisson means, whereas in Definition \ref{def:are_bin_k1}, ${\bm \delta}_{(1)}^{\sf neb}(\lambda)$ is the \texttt{NEB} estimator of the Binomial odds. }
\end{remark}

\begin{remark}\rm{
If for some $i$, $y_i+1$ is not available in the observed sample $\bm y$, ${\psi}_\lambda(y_i)$ can be calculated using cubic splines, and a linear interpolation can be used to tackle the boundary point of the observed sample maxima.
}
\end{remark}

We propose the following estimate of the tuning parameter $\lambda$ based on the ARE:
\begin{equation}
\label{eq:lam_k1}
\hat{\lambda}=\begin{cases}
\argmin_{\lambda\in{\Lambda}}{\sf ARE}^{(1,\mc{P})}_n(\lambda),~\text{if~}Y_i~|~\theta_i\stackrel{ind.}{\sim} \texttt{Poi}(\theta_i)\\
\argmin_{\lambda\in{\Lambda}}{\sf ARE}^{(1,\mc{B})}_n(\lambda),~\text{if~}Y_i~|~q_i\stackrel{ind.}{\sim} \texttt{Bin}(m,q_i)
\end{cases}.
\end{equation}
In practice we recommend using $\Lambda =[10,10^{2}]$, which works well in all our simulations and real data analyses. In Section \ref{sec:theory}, we present Lemmas \ref{lem:lam_bin_k1} and \ref{lem:lam_pois_k1} to provide asymptotic justifications for selecting $\lambda$ using equation \eqref{eq:lam_k1}.

\section{Theory}
\label{sec:theory}

This section studies the theoretical properties for the \texttt{NEB} estimator under the Poisson and Binomial models.  We first investigate the large-sample behavior of the KSD measure (Section \ref{sec:KSD}), then turn to the performance of the estimated risk ratios $\hat{\pmb w}_n$ (Section \ref{sec:w}), and finally establish the consistency and risk properties of the proposed estimator $\bm{\delta}^{\sf neb}$ (Section \ref{sec:neb-theory}). The accuracy of the ${\sf ARE}$ criteria, which are used in choosing tuning parameter $\lambda$, will also be investigated. 

\subsection{Theoretical properties of the KSD measure}\label{sec:KSD}

To provide motivation and theoretical support for Definition \ref{def:neb_k1}, we introduce the Kernelized Stein's Discrepancy (KSD) \citep{liu2016kernelized,chwialkowski2016kernel} and discuss its connection to the quadratic program \eqref{eq:quad_opt}.  While the KSD has been used in various contexts including goodness of fit tests \citep{liu2016kernelized}, variational inference \citep{liu2016stein} and Monte Carlo integration \citep{oates2017control}, our theory on its connection to the compound estimation problem and empirical Bayes methodology is novel.

Assume that $Y$ and $Y'$ are i.i.d. copies from the marginal pmf $p$. Consider $h_0$ defined in Equation \eqref{eq:h(y)}\footnote{In Section \ref{sec:KSD} we shall drop the superscript from $\bm h_0$, which is used to indicate whether the loss is scaled or regular. The simplification has no impact since the general idea holds for both types of losses and the discussion in this section focuses on the scaled  loss.}. Let $\tilde{p}$ denote a pmf on the support of $Y$, for which we similarly define $\tilde{h}$. The KSD, which is formally defined as 
\begin{equation}
\label{eq:ksd}
\mc{S}_\lambda[\tilde{p}](p)=
\mathbb{E}_p\Big[\left\{\tilde{h}(Y)-h_0(Y)\right\}\mc{K}_\lambda(Y,Y')\left\{\tilde{h}(Y')-h_0(Y')\right\}\Big],
\end{equation}
provides a discrepancy measure between $p$ and $\tilde{p}$ in the sense that (a)
$$
\mbox{$\mc{S}_\lambda[\tilde{p}](p)\ge0$ and $\mc{S}_\lambda[\tilde{p}](p)=0$ if and only if $p=\tilde{p}$}, 
$$
and (b) informally, $\mc{S}_\lambda[\tilde{p}](p)$ tends to increase when there is a bigger disparity between $h_0$ and $\tilde h$ (or equivalently, between $p$ and $\tilde p$).

The direct evaluation of $\mc{S}_\lambda[\tilde{p}](p)$ via Equation \eqref{eq:ksd} is difficult 
because $h_0$ is unknown. {Note that while the pmf $p$ can be learned well from a random sample $\{Y_1,\ldots,Y_n\} \sim p$, we introduce an alternative representation of KSD, developed by \cite{liu2016kernelized},  in a reproducing kernel Hilbert space (RKHS) that does not directly involve unknown $h_0$.} Concretely,
consider a smooth positive definite kernel function $\kappa_{\lambda}[\tilde{p}]$:
\begin{equation}\label{eq:kappa}
 \kappa_{\lambda}[\tilde{p}](u,v) = 
\tilde{h}(u)\tilde{h}(v) \mc{K}_\lambda(u,v)   
+    \tilde{h}(u) \Delta_{v} \mc{K}_\lambda(u,v) 
+  \tilde{h}(v) \Delta_{u} \mc{K}_\lambda(u,v)
+  \Delta_{u,v} \mc{K}_\lambda(u,v).
\end{equation}
For i.i.d. copies $(Y, Y^\prime)$ from distribution $p$, it can be shown that
\begin{eqnarray}\label{eq:sqp}
\mc{S}_\lambda[\tilde{p}](p) 
& = & \mathbb{E}_{(Y,Y^\prime)\stackrel{i.i.d.}{\sim}p}\Big[\kappa_{\lambda}[\tilde{p}](Y,Y^\prime)\Big]
 \\ 
& = & \dfrac{1}{n(n-1)}\mathbb{E}_p\Big[\sum_{1\le i\neq j\le n}\kappa_{\lambda}[\tilde{h}(Y_i),\tilde{h}(Y_j)](Y_i,Y_j)\Big] \nonumber \\
& \coloneqq & { {\mathbb{M}}_{\lambda}(\tilde{\bm h}) ,} \nonumber
\end{eqnarray}
where $\{Y_1,\ldots,Y_n\}$ is a random sample from $p$. It can be similarly shown that $\mathbb{M}_{\lambda}(\tilde{\bm h})=0$ if and only if $\tilde{\bm h}=\bm h_0$. 
Substituting the empirical distribution $\hat p_n$ in place of the pmf $p$ in \eqref{eq:sqp},  we obtain the following empirical evaluation scheme for $\mc{S}_\lambda[\tilde{p}](p)$ that is both intuitive and computationally efficient:
\begin{equation}\label{eq:est_sqp}
\mc{S}_\lambda[\tilde{p}](\hat{p}_n)=\dfrac{1}{n^2}\sum_{i=1}^{n}\sum_{j=1}^n\kappa_{\lambda}[{\tilde{h}(y_i),\tilde{h}(y_j)}](y_i,y_j)\coloneqq \hat{\mathbb{M}}_{\lambda,n}(\tilde{\bm h}).
\end{equation}
Note that \eqref{eq:est_sqp} is exactly the objective function of the quadratic program \eqref{eq:quad_opt}. 

The empirical representation of KSD \eqref{eq:est_sqp} provides an extremely useful tool for solving the discrete compound  decision problem under the EB estimation framework. A key observation is that the kernel function $\kappa_{\lambda}[\tilde{p}](u,v)$ depends on $\tilde{p}$ only through $\tilde{h}$. Meanwhile, the EB implementation of the generalized Robbins' formula [cf. Equations \eqref{wk} and \eqref{eq:h(y)}] essentially boils down to the estimation of $h_0$. 
Hence, if $\mc{S}_\lambda[\tilde{p}](\hat{p}_n)$ is asymptotically equal to $\mc{S}_\lambda[\tilde{p}](p)$, then minimizing $\mc{S}_\lambda[\tilde{p}](\hat{p}_n)$ with respect to the unknowns $\tilde{\bm h}=\left\{\tilde{h}(y_1),\ldots,\tilde{h}(y_n)\right\}$ is effectively the process of finding an $\tilde{h}$ that is as close as possible to $h_0$, which yields an asymptotically optimal solution to the EB estimation problem. 
Therefore our formulation of the \texttt{NEB} estimator $\bm{\delta}^{\sf neb}(\lambda)$ would be justified as long as we can establish the asymptotic consistency of the sample criterion $\mc{S}_\lambda[\tilde{p}](\hat{p}_n)$ around the population criterion 
$\mc{S}_\lambda[\tilde{p}](p)$ uniformly over $\lambda$ (Theorem \ref{thm:spq}). 



For a fixed mass function $\tilde{p}$ on the support of $Y$, we impose the following regularity conditions that are needed in our technical analysis.
\begin{itemize}
	\item[(A1)]$\mathbb{E}_p| \kappa_{\lambda}[\tilde{h}(U),\tilde{h}(V)](U,V)|^2<\infty$ for all $\lambda\in\Lambda$ where $\Lambda$ is a compact subset of $\mathbb R^+$ bounded away from $0$.
	\item[(A2)] For some $\epsilon\in(0,1)$, $\limsup_{n\to\infty}n^{-1}\sum_{i=1}^{n}\exp(\epsilon\theta_i)<\infty$.
	\item[(A3)] For any function $g$ that satisfies $0<\|g\|_2^2<\infty$, there exists a constant $c>0$ such that $\sum_{i,j=1}^{n}g(y_i)\mc{K}_\lambda(y_i,y_j)g(y_j)>c\|g\|_2^2$ for every $\lambda\in\Lambda$. 
	\item[(A4)]The feasible solutions $\bm h_n$  to equations \eqref{eq:quad_opt} and \eqref{eq:quad_opt_k0} satisfy $\sup_{\bm h_n\in\bm H_n}\|\bm h_n\|_1=O(n\log n)$.
\end{itemize}

\begin{remark}\rm{
Assumption (A1) is a standard moment condition on the kernel function related to V-statistics, see, for example, \cite{serfling2009approximation}. Assumption (A2) ensures that with high probability $\text{max}(Y_1,\ldots,Y_n)\le \log n$ as $n\to \infty$. This idea is formalized by Lemma \ref{lem:B} in Appendix \ref{sec:proofs}. Assumption (A3) is a standard condition which ensures that the KSD $\mc{S}_\lambda[\tilde{p}](p)$ is a valid discrepancy measure \citep{liu2016kernelized,chwialkowski2016kernel}. Assumption (A4) provides a control on the growth rate of the $\ell_1$ norm of the feasible solutions. In particular, both Assumptions (A3) and (A4) play a critical role in establishing point-wise Lipschitz stability of the optimal solution $\hat{\bm h}_n(\lambda)$ under perturbations on the tuning parameter $\lambda\in\Lambda$ (see Lemma \ref{lem:C} in Appendix \ref{sec:proofs}).}
\end{remark}

\renewcommand*{\thetheorem}{1} 
\begin{theorem}
	\label{thm:spq}
	If $p$ and $\tilde{p}$ are probability mass functions on the support of $Y$ then, under Assumptions (A1) and (A2), we have
	$$\sup_{\lambda\in\Lambda}\Big|\hat{\mathbb{M}}_{\lambda,n}(\tilde{\bm h})-{\mathbb{M}}_{\lambda}(\tilde{\bm h})\Big|=O_p\Big(\dfrac{\log^2 n}{\sqrt{n}}\Big).$$
\end{theorem}
In the context of our compound estimation framework, Theorem \ref{thm:spq} is significant because it guarantees that the empirical version of the KSD measure given by $\hat{\mathbb{M}}_{\lambda,n}(\tilde{\bm h})$ is asymptotically close to its population counterpart $\mathbb{M}_{\lambda}(\tilde{\bm h})$ uniformly in $\lambda\in\Lambda$. Moreover, along with the fact that $\mathbb{M}_{\lambda}(\bm h_0)=0$, Theorem \ref{thm:spq} establishes that $\hat{\mathbb{M}}_{\lambda,n}({\bm h})$ is the appropriate criteria to minimize with respect to $\bm h\equiv \tilde{\bm h}$. In Theorem \ref{thm:w}, we  further show that the resulting estimator of the ratio functionals $\bm w_p^{(1)}$ from equation \eqref{eq:quad_opt} are consistent.

\subsection{Theoretical properties of $\hat{\pmb w}_n$}\label{sec:w}

The optimization problem in \eqref{eq:quad_opt} is defined over a convex set $\bm H_n \subset \mathbb{R}^n$. However, the dimension of $\bm H_n$, denoted by $\texttt{dim}(\bm H_n)$, is usually much smaller than $n$. Consider the Binomial case where $Y_i|q_i\sim \texttt{Bin}(m_i, q_i)$ with $q_i\in(0,1)$, $m_i\leq m<\infty$ and $\theta_i=q_i/(1-q_i)$. Here $\texttt{dim}(\bm H_n)$ is at most $m$ since $\text{max}(Y_1,\ldots,Y_n)\le m$. 
While the boundedness of the support is not always available outside the Binomial case, in most practical applications it is reasonable to assume that the distribution of $\theta_i$ has some finite moments, which ensures that $\texttt{dim}(\bm H_n)$ grows slower than $\log n$; see Assumption (A2). In  Lemma \ref{lem:B} we make this precise. The next theorem establishes the asymptotic consistency of  $\hat{\bm w}_n^{(1)}(\lambda)$. 
\renewcommand*{\thetheorem}{2A} 
\begin{theorem}
	\label{thm:w}
	Let $\mc{K}_\lambda(\cdot,\cdot)$ be the positive definite RBF kernel with bandwidth parameter $\lambda\in\Lambda$. If $\lim_{n\to\infty}c_nn^{-1/2}\log^2n=0$ then, under Assumptions (A1) - (A3), we have for any $\lambda\in\Lambda$,
$$\lim_{n\to\infty}\mathbb{P}\left\{\Big\|\hat{\bm w}_n^{(1)}(\lambda)-{\bm w}_p^{(1)}\Big\|_2\ge c_n^{-1}\epsilon\right\}=0,~\text{for any }\epsilon>0, $$
where $\hat{\bm w}_n^{(1)}(\lambda)=1-\hat{\bm h}_n^{(1)}(\lambda)$. 
\end{theorem}

Theorem \ref{thm:w} shows that under the scaled squared error loss, $\hat{\bm w}_n^{(1)}(\lambda)$, the optimizer of quadratic form \eqref{eq:quad_opt}, provides a consistent estimator of ${\bm w}_p^{(1)}$, the optimal shrinkage factor in the Bayes rule (Lemma \ref{lem:bayes}). Theorem \ref{thm:w} is proved in Appendix \ref{sec:proof_thm2}, where we also include relevant details for proving a companion result under the regular squared error loss.

\begin{remark}\rm{
The estimation framework in Definition \ref{def:neb_k1} may be used for producing consistent estimators for any member in the \texttt{DLE} family. This allows the corresponding \texttt{NEB} estimator to cover a much wider class of discrete distributions than previously proposed. Compared to existing methods \citep{efron2011tweedie,brown2013poisson,koenker2014convex,koenker2017rebayes},
our proposed \texttt{NEB} estimation framework is robust against departures from the true data generating process. This is due to the fact that the quadratic optimization problem in \eqref{eq:quad_opt} does not rely on the specific form of the distribution of $Y|\theta$, and that the shrinkage factors are estimated in a non-parametric fashion. 
The robustness of the estimator is corroborated by our numerical results in Section \ref{sec:numresults}.}
\end{remark}

 
\subsection{Properties of the \texttt{NEB} estimator}\label{sec:neb-theory}

In this section we discuss the risk properties of the \texttt{NEB} estimator. We begin with two lemmas showing that uniformly in $\lambda\in\Lambda$, the gap between the estimated risk ${\sf ARE}^{(1)}_n(\lambda)$ 
and true risk is asymptotically negligible. This justifies our proposed methodology for choosing the tuning parameter $\lambda$ in Section \ref{sec:bandwidth}. In the following two lemmas, we let $c_n$ be a sequence satisfying $\lim_{n\to\infty}c_nn^{-1/2}\log^{5/2}n=0$. 
\renewcommand*{\thelemma}{2}
\begin{lemma}
	\label{lem:lam_bin_k1}
Under Assumptions (A3) and (A4) and the Binomial model, we have
	\begin{eqnarray}
	(a). & c_n\sup_{\lambda\in\Lambda}\Big|{\sf ARE}^{(1,\mc{B})}_n(\lambda,\bm Y)-\mc{R}_n^{(1)}(\bm \theta,{\bm \delta}_{(1)}^{\sf neb}(\lambda))\Big|=o_p(1);\nonumber\\
	(b). & c_n\sup_{\lambda\in\Lambda}\Big|{\sf ARE}^{(1,\mc{B})}_n(\lambda,\bm Y)-\mc{L}_n^{(1)}(\bm \theta,{\bm \delta}_{(1)}^{\sf neb}(\lambda))\Big|=o_p(1).\nonumber
	\end{eqnarray}
\end{lemma}
\setcounter{lemma}{3}
\renewcommand*{\thelemma}{3}
\begin{lemma}
	\label{lem:lam_pois_k1}
Under Assumptions (A2), (A3) and (A4) and the Poisson model, we have
	\begin{eqnarray}
	(a). & c_n\sup_{\lambda\in\Lambda}\Big|{\sf ARE}^{(1,\mc{P})}_n(\lambda,\bm Y)-\mc{R}_n^{(1)}(\bm \theta,{\bm \delta}_{(1)}^{\sf neb}(\lambda))\Big|=o_p(1);\nonumber\\
	(b). & c_n\sup_{\lambda\in\Lambda}\Big|{\sf ARE}^{(1,\mc{P})}_n(\lambda,\bm Y)-\mc{L}_n^{(1)}(\bm \theta,{\bm \delta}_{(1)}^{\sf neb}(\lambda))\Big|=o_p(1). \nonumber
	\end{eqnarray}
\end{lemma}

To analyze the quality of the data-driven bandwidth $\hat{\lambda}$ [cf. Equation \eqref{eq:lam_k1}], we consider an oracle loss estimator $\bm \delta^{\sf or}_{(1)}\coloneqq {\bm \delta_{(1)}^{\sf neb}}(\lambda^{\sf orc}_1)$, where 
$$
\lambda^{\sf orc}_1\coloneqq\argmin_{\lambda\in\Lambda}\mc{L}_n^{(1)}\left\{\bm \theta,{\bm \delta}_{(1)}^{\sf neb}(\lambda)\right\}.
$$ 
The oracle bandwidth $\lambda^{\sf orc}_1$ is not available in practice since it requires the knowledge of unknown $\pmb\theta$. However, it provides a benchmark for assessing the effectiveness of the data-driven bandwidth selection procedure in Section \ref{sec:bandwidth}. The following lemma shows that the loss of ${\bm \delta}^{\sf neb}_{(1)}(\hat{\lambda})$ converges in probability to the loss of $\bm \delta^{\sf or}_{(1)}$.
\renewcommand*{\thelemma}{4} 
\begin{lemma}
	\label{lem:horc_k1}
	Under Assumptions (A2) - (A4), if $\lim_{n\to\infty}c_nn^{-1/2}\log^{5/2} n=0$,  then for both the Poisson and Binomial models, we have 
	$$
	\lim_{n\to\infty}\mathbb{P}\Big[\mc{L}_n^{(1)}\left\{\bm \theta,{\bm \delta}_{(1)}^{\sf neb}(\hat{\lambda})\right\}\ge \mc{L}_n^{(1)}(\bm \theta,{\bm \delta}_{(1)}^{\sf or})+c_n^{-1}\epsilon \Big]=0\text{ for any }\epsilon>0.
	$$
\end{lemma}

Obviously, the estimator $\bm\delta_{(1)}^{\sf neb}(\lambda^{\sf orc}_1)$
 is lower bounded by the risk of the optimal solution $\bm \delta_{(1)}^{\pi}$ (Lemma \ref{lem:bayes}). Next we study the \emph{asymptotic optimality} of $\bm \delta^{\sf neb}_{(1)}$, which aims to provide decision theoretic guarantees on $\bm \delta^{\sf neb}_{(1)}$ in relation to $\bm \delta_{(1)}^{\pi}$. Theorem \ref{thm:bayesrisk_k1} establishes the optimality theory by showing that (a) the largest coordinate-wise gap between  ${\bm \delta}^{\sf neb}_{(1)}({\hat{\lambda}})$ and $\bm \delta_{(1)}^{\pi}$ is asymptotically small, and (b) the estimation loss of the \texttt{NEB} estimator converges in probability to the loss of the corresponding Bayes estimator as $n\to \infty$. 

\renewcommand*{\thetheorem}{3A} 
\begin{theorem}
	\label{thm:bayesrisk_k1}
	Under the conditions of Theorem \ref{thm:w}, if {$\lim_{n\to\infty}c_nn^{-1/2}\log^{4} n=0$}, then for both the Poisson and Binomial models, we have
		$$ c_n\Big\|{\bm\delta}_{(1)}^{\sf neb}({\hat\lambda})-\bm\delta_{(1)}^{\pi}\Big\|_\infty=o_p(1).$$
	Furthermore, under the same conditions, we have for both the Poisson and Binomial models,
	$$\lim_{n\to\infty}\mathbb{P}\Big[\mc{L}_n^{(1)}(\bm \theta,{\bm \delta}_{(1)}^{\sf neb}(\hat{\lambda}))\ge \mc{L}_n^{(1)}(\bm \theta,\bm \delta_{(1)}^{\pi})+c_n^{-1}\epsilon \Big]=0\text{ for any }\epsilon>0.$$
\end{theorem}

\begin{remark}\rm{
The second part of the statement of Theorem \ref{thm:bayesrisk_k1} follows from the first part and Lemma \ref{lem:horc_k1}.  In Appendix \ref{sec:lam_k0}, we discuss the counterpart to Theorem \ref{thm:bayesrisk_k1} under the squared error loss $\mc{L}_n^{(0)}$.}
\end{remark}

\section{Numerical Results}
\label{sec:numresults}
In this section we first discuss, in Section \ref{sec:rpackage}, the implementation details of the convex program \eqref{eq:quad_opt} and bandwidth selection process \eqref{eq:lam_k1} (see also \eqref{eq:lam_k0} in Appendix \ref{sec:lam_k0}). Then we investigate the numerical performance of the \texttt{NEB} estimator for Poisson and Binomial compound decision problems, respectively in Sections \ref{sec:sims} and \ref{simu-binom}. Both regular and scales losses will be considered. Our numerical results demonstrate that the proposed the \texttt{NEB} estimator enjoys superior numerical performance and the efficiency gain over competitive methods is substantial in many settings. 

We have developed an R package, \href{https://github.com/trambakbanerjee/npeb#npeb}{\texttt{npeb}}, to implement our proposed \texttt{NEB} estimator in Definitions \ref{def:neb_k1} (and Definition \ref{def:neb_k0} in Appendix \ref{sec:neb_sqaurederror}), for the Poisson and Binomial models under both regular and scaled losses. 
Moreover, the R code that reproduces the numerical results in simulations can be downloaded from the following link:

\url{https://github.com/trambakbanerjee/DLE_paper}.

\subsection{Implementation Details}
\label{sec:rpackage}

For a fixed $\lambda$ we use the R-package \texttt{CVXR} \citep{fu2017cvxr} to solve the optimization problem in Equations \eqref{eq:quad_opt} [and \eqref{eq:quad_opt_k0} in Appendix \ref{sec:neb_sqaurederror}]. As discussed in Section \ref{sec:min_dkskd}, in the implementation under the scaled squared error loss ($k=1$) the linear inequality constraints, given by $\mc{\bm A}\bm h\preceq\bm b$, ensure that the resulting decision rule ${\bm \delta}^{\sf neb}_{(1)}(\lambda)$ is monotonic, while the equality constraints $\mc{\bm C}\bm h=\bm d$ handle boundary cases that involve $y_i=0$ and ties. Moreover, since $w_p^{(1)}(y)>0$, the inequality constraints also ensure that $h_i<1$ whenever $y_i>0$. Implementation under the squared error loss ($k=0$) follows along similar lines and the inequality constraints in this case ensure that $h_i+y_i>0$ whenever $y_i\ge 0$. 
 
A data-driven choice of the tuning parameter $\lambda$ is obtained by first solving problems \eqref{eq:quad_opt} and \eqref{eq:quad_opt_k0} over a grid of $\lambda$ values, i.e. $\{\lambda_1,\ldots,\lambda_s\}$, and then computing the corresponding asymptotic risk estimate $\texttt{ARE}_n^{(k)}(\lambda_j)$ for $j=1,\ldots,s$. Then $\lambda$ is chosen according to 
$$
\hat{\lambda}_k\coloneqq\argmin_{\lambda\in{\lambda_1,\ldots,\lambda_s}}\texttt{ARE}_n^{(k)}(\lambda),
$$ 
where $k\in \{0,1\}$. For all simulations and real data analyses considered in this paper, we have fixed $s=10$ and employed an equi-spaced grid over $[10,10^{2}]$.

\subsection{Simulations: Poisson Distribution}
\label{sec:sims}
In this section, we generate observations $Y_i~|~\theta_i\stackrel{ind.}{\sim}\texttt{Poi}(\theta_i)$ for $i=1,\ldots,n$ and vary $n$ from $500$ to $5000$ in increments of $500$. We consider four different scenarios to simulate $\theta_i$. For each scenario, the following competing estimators are considered: 
\begin{itemize}
\item the proposed estimator, denoted \texttt{NEB};
\item {the oracle \texttt{NEB} estimator $\bm \delta^{\sf or}\coloneqq {\bm \delta^{\sf neb}}(\lambda^{\sf orc})$, denoted \texttt{NEB OR}};
\item the estimator of Poisson means based on \cite{brown2013poisson}, denoted \texttt{BGR};
\item the estimator of Poisson means based on \cite{koenker2017rebayes}, denoted \texttt{KM};  
\item Tweedie's formula based on \cite{efron2011tweedie} for the Poisson model, denoted \texttt{TF OR}; 
\item Tweedie's formula for the Normal means problem based on transformed data and the convex optimization approach in \cite{koenker2014convex}, denoted \texttt{TF Gauss}. The approach using transformation was suggested by \cite{brown2013poisson}. 
\end{itemize}

The risk performance of the \texttt{TF OR} method relies heavily on the choice of a suitable bandwidth parameter $h>0$. We use the oracle loss estimate $h^{\sf orc}$, which is obtained by minimizing the true loss $\mc{L}_n^{(0)}$. The \texttt{TF Gauss} methodology is only applicable for the Normal means problem, and uses a variance stabilization transformation on $Y_i$ to get $Z_i=2\sqrt{Y_i+0.25}$. The $Z_i$ are then treated as approximate Normal random variables with mean $\mu_i$ and variances $1$. To estimate normal means $\mu_i$, we use the NPMLE approach proposed by \cite{koenker2014convex}. Finally, $\theta_i$ is estimated as $0.25\hat{\mu_i}^2$.

It is important to note that the competitors to our \texttt{NEB} estimator only focus on the regular loss $\mc{L}_n^{(0)}$. Nevertheless, in our simulation we assess the performance of these estimators for estimating $\bm \theta$ under both $\mc{L}_n^{(0)}$ and $\mc{L}_n^{(1)}$. Consider the following settings: 
 
\begin{description}

\item \textit{Scenario 1: }We generate $\theta_i\stackrel{iid.}{\sim}\texttt{Unif}(0.5,15)$ for $i=1,\ldots,n$.

\item \textit{Scenario 2: }We generate $\theta_i\stackrel{iid.}{\sim}0.75~\texttt{Gamma}(5,1)+0.25~\texttt{Gamma}(10,1)$ for $i=1,\ldots,n$. 

\medskip

In the next two scenarios we assess the robustness of the five competing estimators to departures from the Poisson model. Specifically consider the Conway-Maxwell-Poisson distribution \citep{shmueli2005useful} $\texttt{CMP}(\theta_i,\nu)$. The \texttt{CMP} distribution is a generalization of some well-known discrete distributions. With $\nu<1$, CMP represents a discrete distribution that has longer tails than the Poisson distribution with parameter $\theta_i$.

\medskip

\item \textit{Scenario 3: }We generate 
$
\theta_i\stackrel{iid.}{\sim}0.5~\delta_{\{10\}}+0.5~\texttt{Gamma}(5,2)
$
for each $i$ and let 
$$Y_i~|~\theta_i\stackrel{ind.}{\sim}0.8~\texttt{Poi}(\theta_i)+0.2~\texttt{CMP}(\theta_i,\nu),
$$ 
where we let $\nu=0.8$ for the CMP distribution. 

\item \textit{Scenario 4: } We let $\bm \theta$ to be an equi-spaced vector of length $n$ in $[1,5]$ and let $Y_i~|~\theta_i$ be distributed as the \texttt{CMP} distribution with parameters $\theta_i$ and $\nu=0.8$. 

\end{description}

 \begin{figure}[!h]
	\centering
	\begin{subfigure}{.46\textwidth}
		\centering
		\includegraphics[width=1\linewidth]{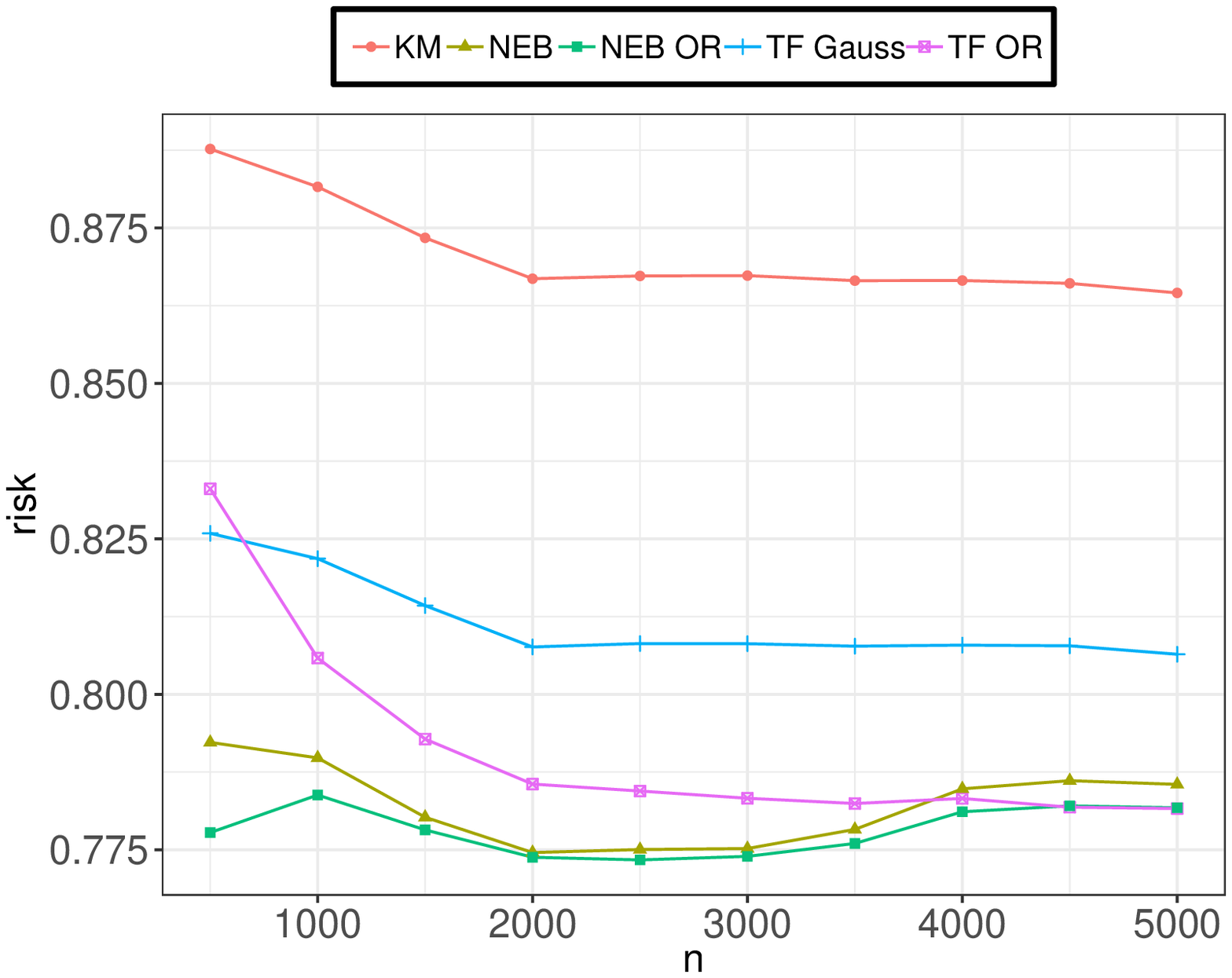}
		\caption{Scenario 1: Estimation of $\bm \theta$ under loss $\mc{L}_n^{(1)}$ where $\theta_i\stackrel{iid.}{\sim}\text{Unif}(0.5,15)$.}
		\label{fig:pois_k1_exp1}
	\end{subfigure}%
	\hfill
	\begin{subfigure}{0.46\textwidth}
		\centering
		\includegraphics[width=1\linewidth]{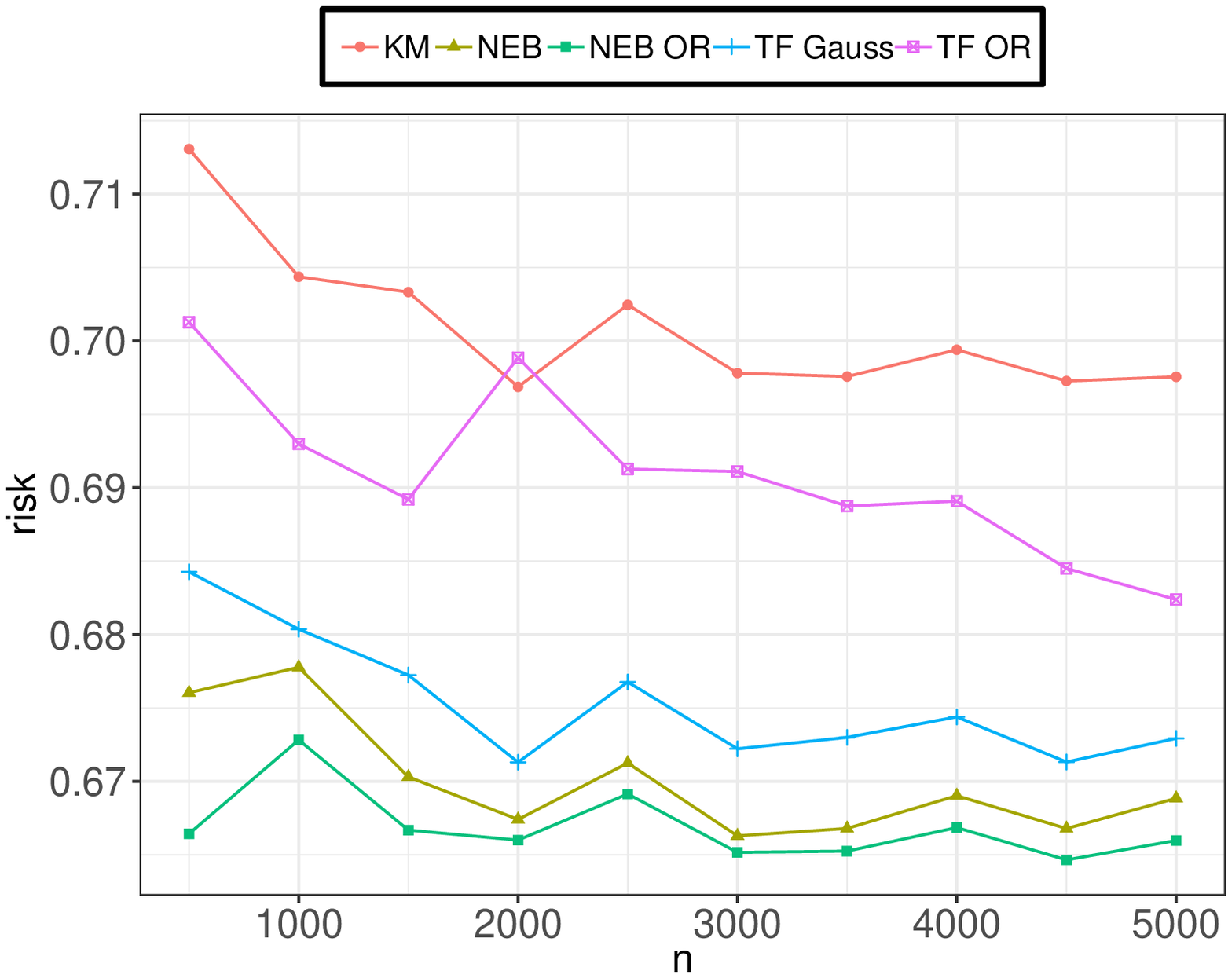}
		\caption{Scenario 2: Estimation of $\bm \theta$ under loss $\mc{L}_n^{(1)}$ where $\theta_i\stackrel{iid.}{\sim}0.75~\texttt{Gamma}(5,1)+0.25~\texttt{Gamma}(10,1)$.}
		\label{fig:pois_k1_exp3}
	\end{subfigure}
	\begin{subfigure}{.46\textwidth}
	\centering
	\includegraphics[width=1\linewidth]{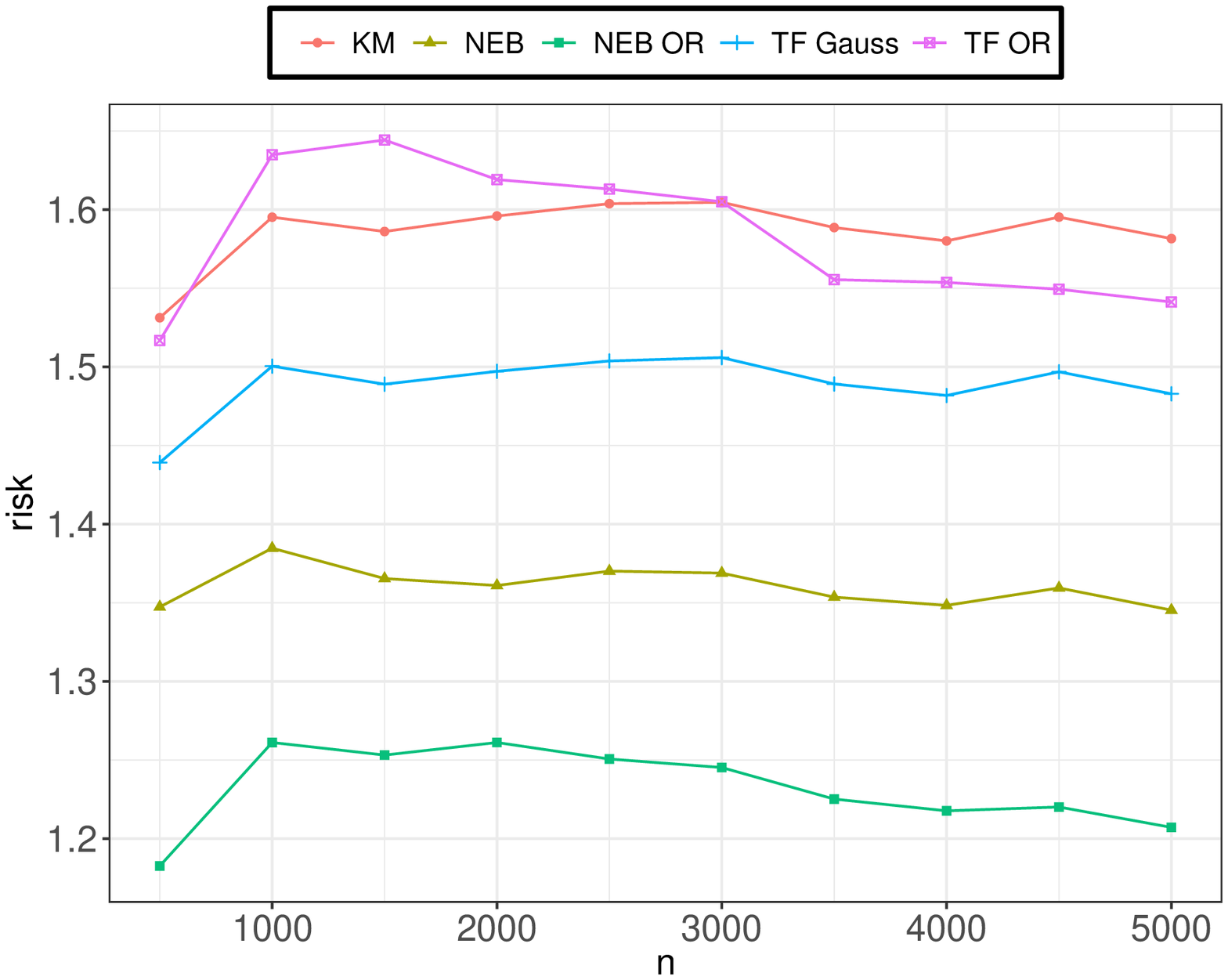}
	\caption{Scenario 3: Estimation of $\bm \theta$ under loss $\mc{L}_n^{(1)}$ where $\theta_i\stackrel{iid.}{\sim}0.5~\delta_{\{10\}}+0.5~\texttt{Gamma}(5,2)$.}
	\label{fig:pois_k1_exp2}
\end{subfigure}%
\hfill
\begin{subfigure}{0.46\textwidth}
	\centering
	\includegraphics[width=1\linewidth]{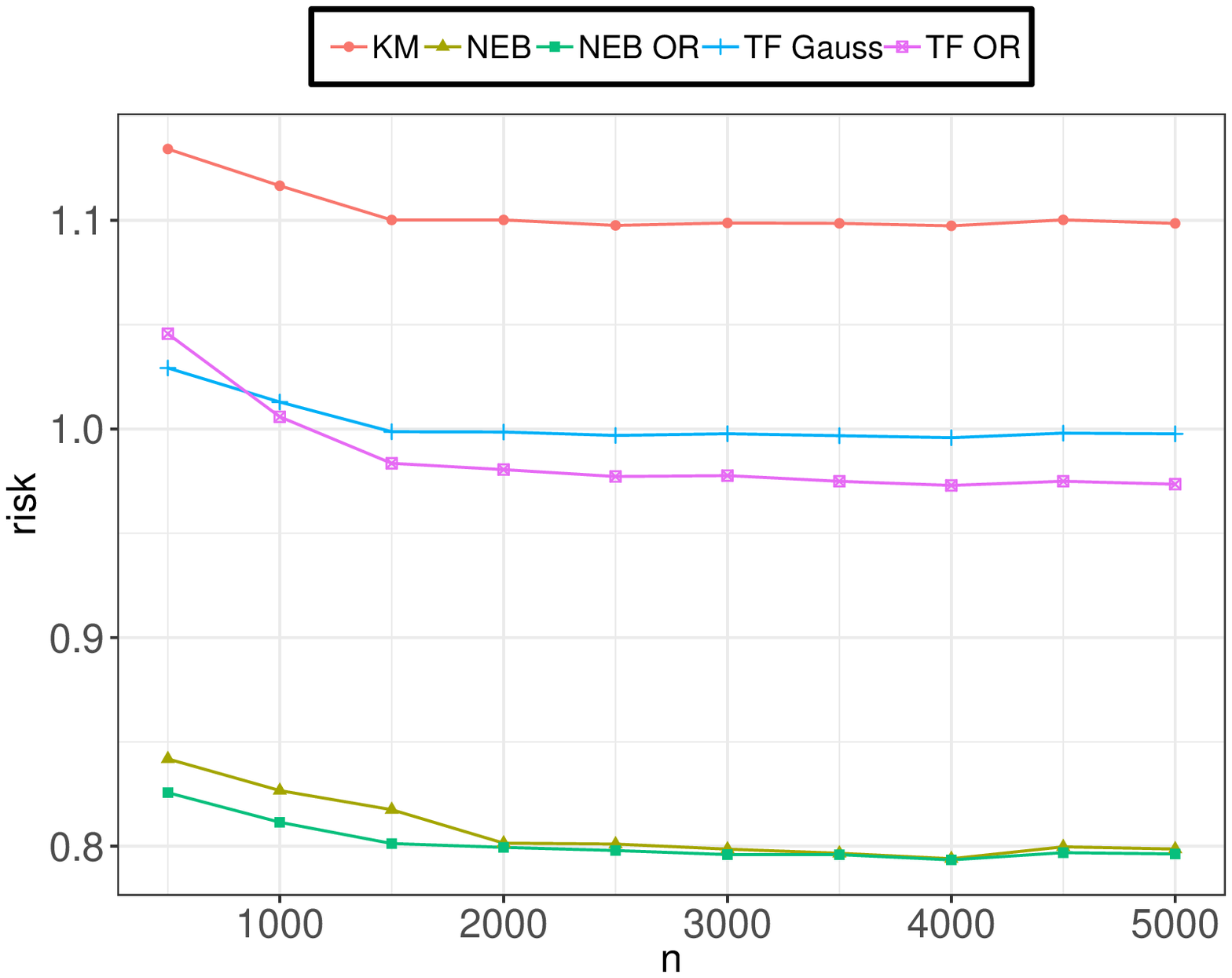}
	\caption{Scenario 4: Estimation of $\bm \theta$ under loss $\mc{L}_n^{(1)}$ where $\bm \theta$ is an equi-spaced vector of length $n$ in $[1,5]$ and $Y_i|\theta_i\stackrel{ind.}{\sim} \texttt{CMP}(\theta_i,0.8)$ .}
	\label{fig:pois_k1_exp5}
\end{subfigure}
\caption{Poisson compound decision problem under scaled squared error loss: Risk estimates of the various estimators for scenarios 1 to 4.}
\label{fig:pois1}
\end{figure}
 \begin{figure}[!h]
	\centering
	\begin{subfigure}{.46\textwidth}
		\centering
		\includegraphics[width=1\linewidth]{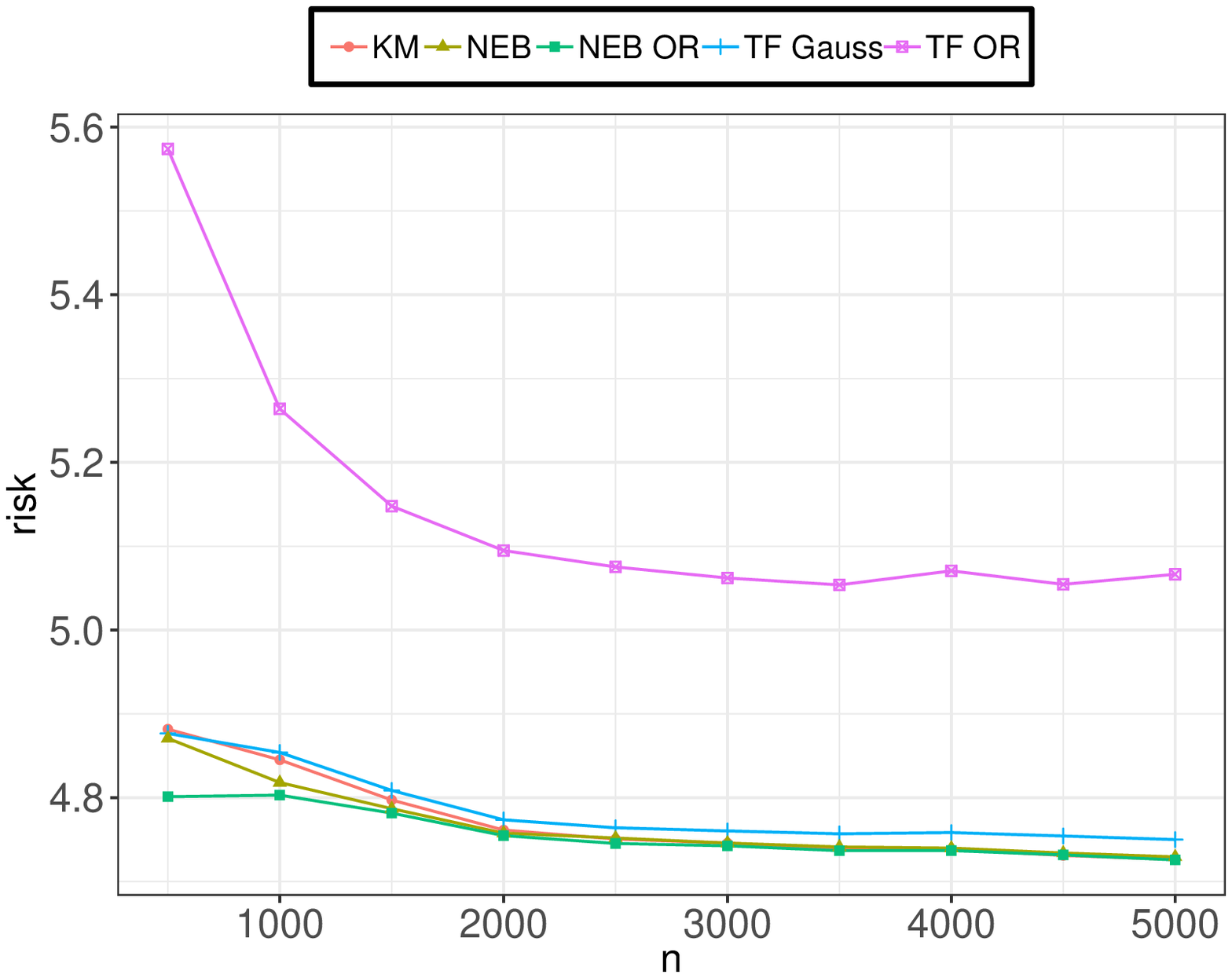}
		\caption{Scenario 1: Estimation of $\bm \theta$ under loss $\mc{L}_n^{(0)}$ where $\theta_i\stackrel{iid.}{\sim}\texttt{Unif}(0.5,15)$.}
		\label{fig:pois_k0_exp1}
	\end{subfigure}%
	\hfill
	\begin{subfigure}{0.46\textwidth}
		\centering
		\includegraphics[width=1\linewidth]{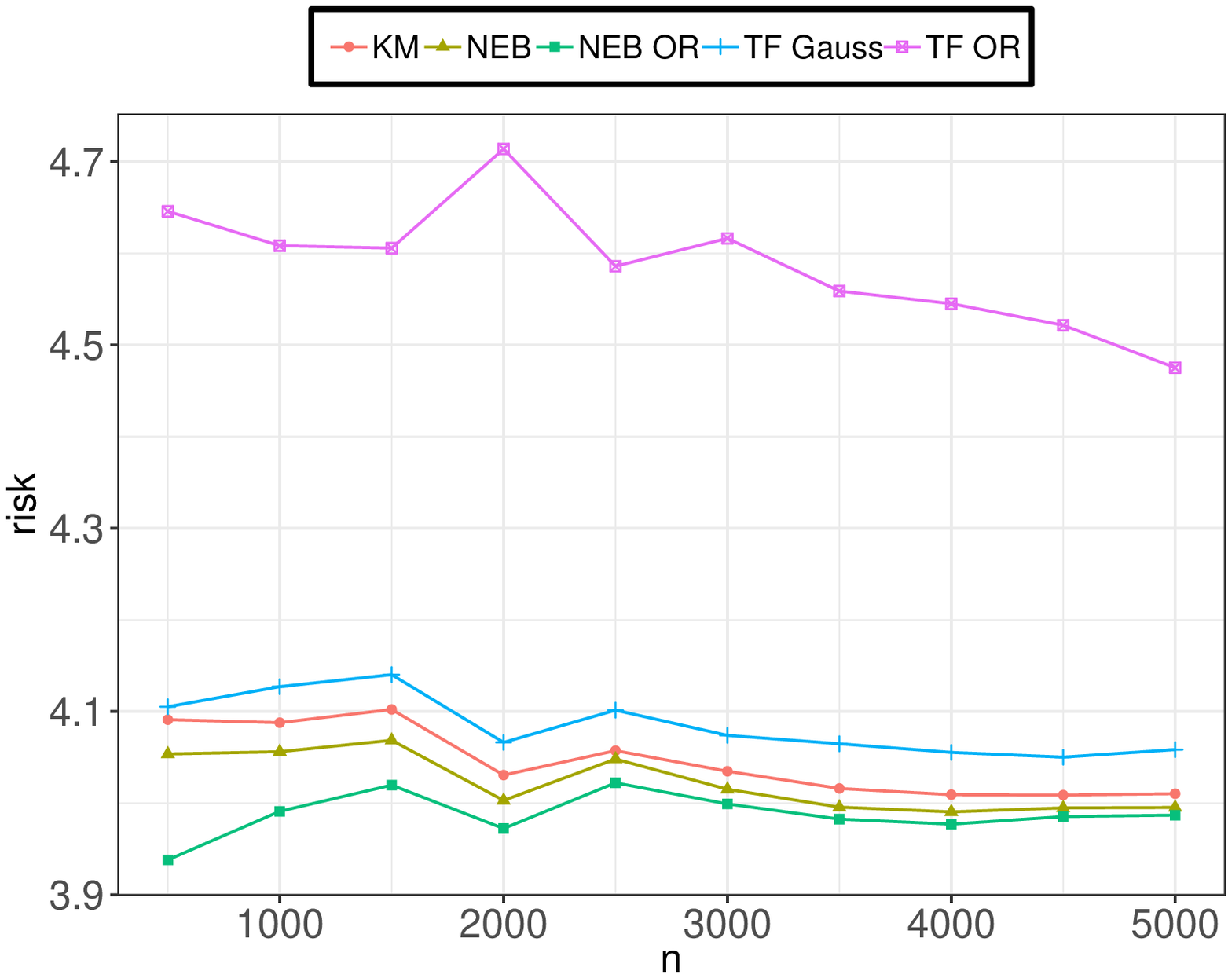}
		\caption{Scenario 2: Estimation of $\bm \theta$ under loss $\mc{L}_n^{(0)}$ where $\theta_i\stackrel{iid.}{\sim}0.75~\texttt{Gamma}(5,1)+0.25~\texttt{Gamma}(10,1)$.}
		\label{fig:pois_k0_exp3}
	\end{subfigure}
	\begin{subfigure}{.46\textwidth}
	\centering
	\includegraphics[width=1\linewidth]{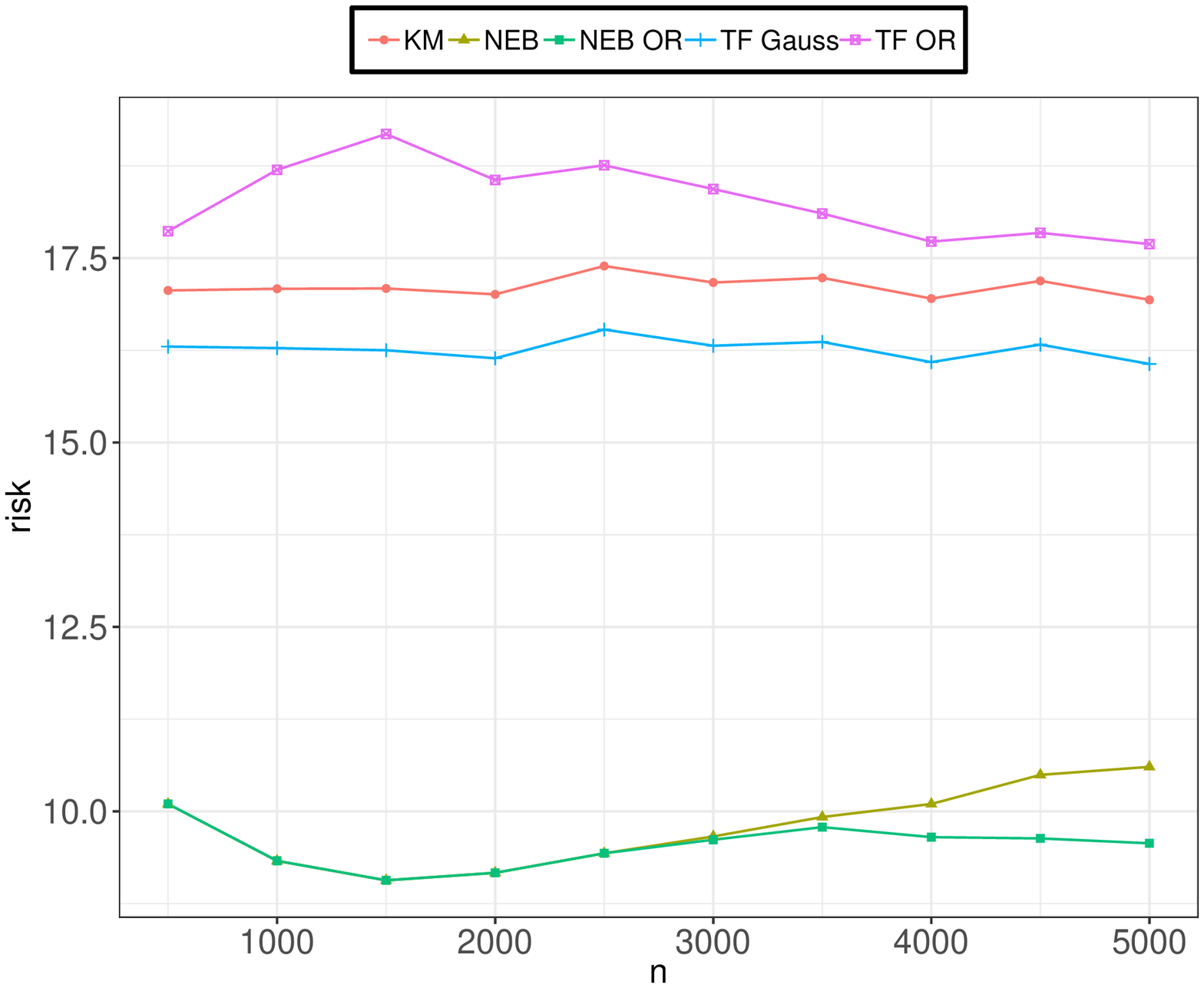}
	\caption{Scenario 3: Estimation of $\bm \theta$ under loss $\mc{L}_n^{(0)}$ where $\theta_i\stackrel{iid.}{\sim}0.5~\delta_{\{10\}}+0.5~\texttt{Gamma}(5,2)$.}
	\label{fig:pois_k0_exp2}
\end{subfigure}%
\hfill
\begin{subfigure}{0.46\textwidth}
	\centering
	\includegraphics[width=1\linewidth]{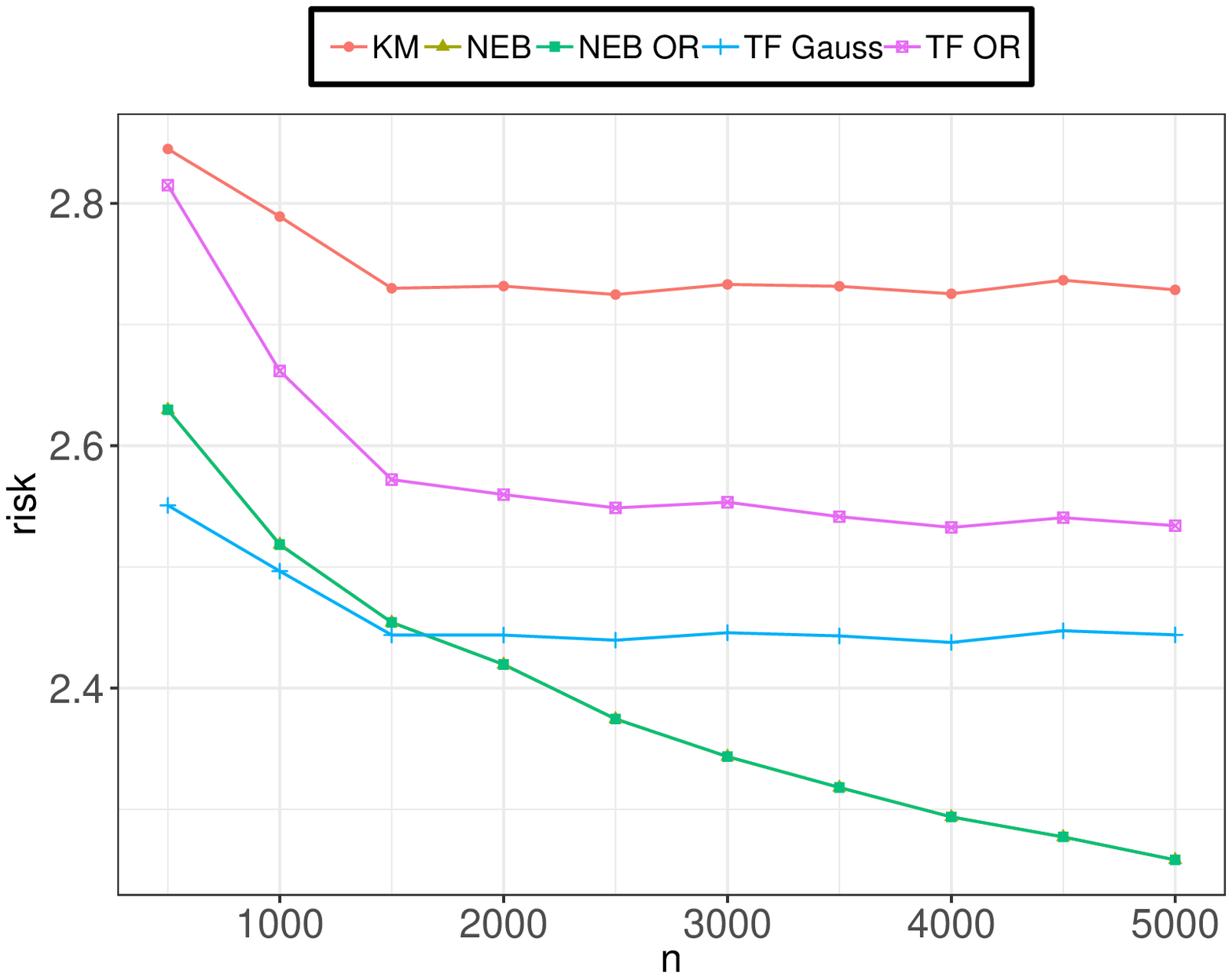}
	\caption{Scenario 4: Estimation of $\bm \theta$ under loss $\mc{L}_n^{(0)}$ where $\bm \theta$ is an equi-spaced vector of length $n$ in $[1,5]$ and $Y_i|\theta_i\stackrel{ind.}{\sim} \texttt{CMP}(\theta_i,0.8)$ .}
	\label{fig:pois_k0_exp5}
\end{subfigure}
\caption{Poisson compound decision problem under squared error loss: Risk estimates of the various estimators for scenarios 1 to 4.}
\label{fig:pois2}
\end{figure}

\begin{table}[!h]
	\begin{minipage}{.47\linewidth}
		\centering
		\caption{Poisson compound decision problem under scaled squared error loss: Risk ratios $\mc{R}_n^{(1)}(\bm \theta,\cdot)/\mc{R}_n^{(1)}(\bm \theta,\bm \delta^{\sf neb}_k)$ at $n=5000$ for estimating $\bm \theta$.}
		\scalebox{0.95}{
			\begin{tabular}{lcccc}
				& \multicolumn{4}{c}{Scenario} \\
				\toprule
				Method & 1     & 2     & 3     & 4 \\
				\midrule
				\texttt{KM}   & 1.10  & 1.04  &  1.17 & 1.37 \\
				\texttt{TF Gauss} & 1.03  & 1.01  & 1.10  & 1.25 \\
				\texttt{TF OR} & 1.00  & 1.02  &1.14   & 1.22 \\
				\texttt{BGR}   & 1.22  & 1.07  & 1.16  & 1.37 \\
				\midrule
				\texttt{NEB}   & 1.00  & 1.00  & 1.00  & 1.00 \\
				\texttt{NEB OR} & 1.00  & 1.00  & 0.90  & 1.00 \\
				\bottomrule
		\end{tabular}}%
		\label{tab:pois1}%
	\end{minipage}%
	\hfill
	\begin{minipage}{.47\linewidth}
		\centering
		\caption{Poisson compound decision problem under squared error loss: Risk ratios $\mc{R}_n^{(0)}(\bm \theta,\cdot)/\mc{R}_n^{(0)}(\bm \theta,\bm \delta^{\sf neb}_k)$ at $n=5000$ for estimating $\bm \theta$.}
		\scalebox{0.95}{
			\begin{tabular}{lcccc}
				& \multicolumn{4}{c}{Scenario} \\
				\toprule
				Method & 1     & 2     & 3     & 4 \\
				\midrule
				\texttt{KM}   & 1.00  & 1.00  & 1.59  & 1.21 \\
				\texttt{TF Gauss} & 1.00  & 1.01  & 1.51  & 1.08 \\
				\texttt{TF OR} & 1.07  & 1.12  & 1.66  & 1.12 \\
				\texttt{BGR}   & 1.01  & 1.01  & 1.55  & 1.15 \\
				\midrule
				\texttt{NEB}   & 1.00  & 1.00  & 1.00  & 1.00 \\
				\texttt{NEB OR} & 1.00  & 1.00  & 0.90  & 1.00 \\
				\bottomrule
		\end{tabular}}%
		\label{tab:pois2}%
	\end{minipage} 
\end{table}
The performances of these four estimators are presented in figures \ref{fig:pois1} and \ref{fig:pois2} wherein the risk $\mc{R}_n^{(k)}(\bm \theta,\cdot)$ of the various estimators is estimated using $50$ Monte Carlo repetitions for varying $n$. Tables \ref{tab:pois1} and \ref{tab:pois2} report the risk ratios $\mc{R}_n^{(k)}(\bm \theta,\cdot)/\mc{R}_n^{(k)}(\bm \theta,\bm \delta^{\sf neb}_k)$ at $n=5000$ and for $k=1,0$ respectively, where a risk ratio bigger than $1$ demonstrates a smaller estimation risk for the \texttt{NEB} estimator. For \texttt{BGR} the modified cross validation approach of choosing the bandwidth parameter was extremely slow in our simulations and we therefore report its risk performance only at $n=5000$. From figures \ref{fig:pois1}, \ref{fig:pois2} and tables \ref{tab:pois1}, \ref{tab:pois2}, we note that the \texttt{NEB} estimator demonstrates an overall competitive risk performance. In particular, we see that when estimation is conducted under loss $\mc{L}_n^{(1)}$ the risk ratios of the competing estimators in table \ref{tab:pois1} reflect a relatively better performance of the \texttt{NEB} estimator which is not surprising considering the fact that \texttt{KM, TF Gauss} and \texttt{TF OR} are designed to estimate $\bm \theta$ under loss $\mc{L}_n^{(0)}$. We note that \texttt{TF Gauss} is highly competitive against \texttt{KM} \citep{koenker2014convex} and this observation was also reported in \cite{brown2013poisson}. Of particular interest are Scenarios 3 and 4, which reflect the relative performance of these estimators under departures from the Poisson model. The \texttt{NEB} estimator has a significantly better risk performance in these settings across both types of losses.

\subsection{Simulations: Binomial Distribution}\label{simu-binom}

In this section, we generate $Y_i~|~q_i\stackrel{ind.}{\sim}\mathcal{B}in(m_i,~q_i)$ for $i=1,\ldots,n$ and vary $n$ from $500$ to $5000$ in increments of $500$. We consider four different scenarios to simulate $\theta_i=q_i/(1-q_i)$ and for each scenario we consider the following competing estimators:

\begin{itemize}
\item the proposed estimator, denoted \texttt{NEB};
\item {the oracle \texttt{NEB} estimator $\bm \delta^{\sf or}\coloneqq {\bm \delta^{\sf neb}}(\lambda^{\sf orc})$, denoted \texttt{NEB OR}};
 \item the estimator of Binomial means based on \cite{koenker2017rebayes}, denoted \texttt{KM};  
\item Tweedie's formula for Binomial log odds based on \cite{efron2011tweedie} and \cite{funonparametric2018}, denoted \texttt{TF OR}; 
\item Tweedie's formula for the Normal means problem based on transformed data and the convex optimization approach in \cite{koenker2014convex}, denoted \texttt{TF Gauss}. 
\end{itemize}

For \texttt{TF OR}, analogous to the Poisson case, we continue to use the oracle loss estimate $h^{\sf orc}$ as a choice for the bandwidth parameter. Since the \texttt{TF Gauss} methodology is only applicable for the Normal means problem, it uses a variance stabilization transformation on $Y_i$ to get $Z_i={\sf arcsin}\sqrt{(Y_i+0.25)/(m_i+0.5)}$. The $Z_i$ are then treated as approximate Normal random variables with mean $\mu_i$, variances $(4m_i)^{-1}$, and estimate of the means $\mu_i$'s are obtained using the NPMLE approach of \cite{koenker2014convex}. Finally, $q_i$ is estimated as $\{\sin(\hat\mu_i)\}^2$. We note that unlike the Poisson case discussed earlier, the competitors to our \texttt{NEB} estimator do not directly estimate the odds $\bm \theta$. For instance, under a squared error loss both \texttt{KM} and \texttt{TF Gauss} estimate the success probabilities $\bm q$ while \texttt{TF OR} estimates $\log \bm\theta$. Nevertheless, in this simulation experiment we assess the performance of these estimators for estimating the odds under both squared error loss and its scaled version.

The following settings are considered in our simulation: 

\begin{description}
 
 \item \textit{Scenario 1: }We generate $q_i\stackrel{iid.}{\sim}0.4~\delta_{\{0.5\}}+0.6~\texttt{Beta}(2,5)$ and fix $m_i=5$ for $i=1,\ldots,n$.  \\[1ex]

 \item \textit{Scenario 2: }We let $\theta_i\stackrel{iid.}{\sim}0.8~\delta_{\{0.5\}}+0.2~\texttt{Gamma}(1,2)$ and fix $m_i=10$ for $i=1,\ldots,n$. In this scenario we let the odds $\theta_i$ arise from a mixture model that has $80\%$ point mass at $0.5$.
 
 \item \textit{Scenario 3: }
 $\theta_i\stackrel{iid.}{\sim}\chi^2_2$
 and fix $m_i=5$ for $i=1,\ldots,n$. This scenario is similar to scenario 2 where we let the odds $\theta_i$ arise from a Chi-square distribution with 2 degrees of freedom.
 
  \item \textit{Scenario 4: }We generate $q_i\stackrel{iid.}{\sim}0.5~\texttt{Beta}(1,1)+0.5~\texttt{Beta}(1,3)$ and fix $m_i=10$ for $i=1,\ldots,n$. 

\end{description}

\begin{figure}[!h]
 	\centering
 	\begin{subfigure}{.46\textwidth}
 		\centering
 		\includegraphics[width=1\linewidth]{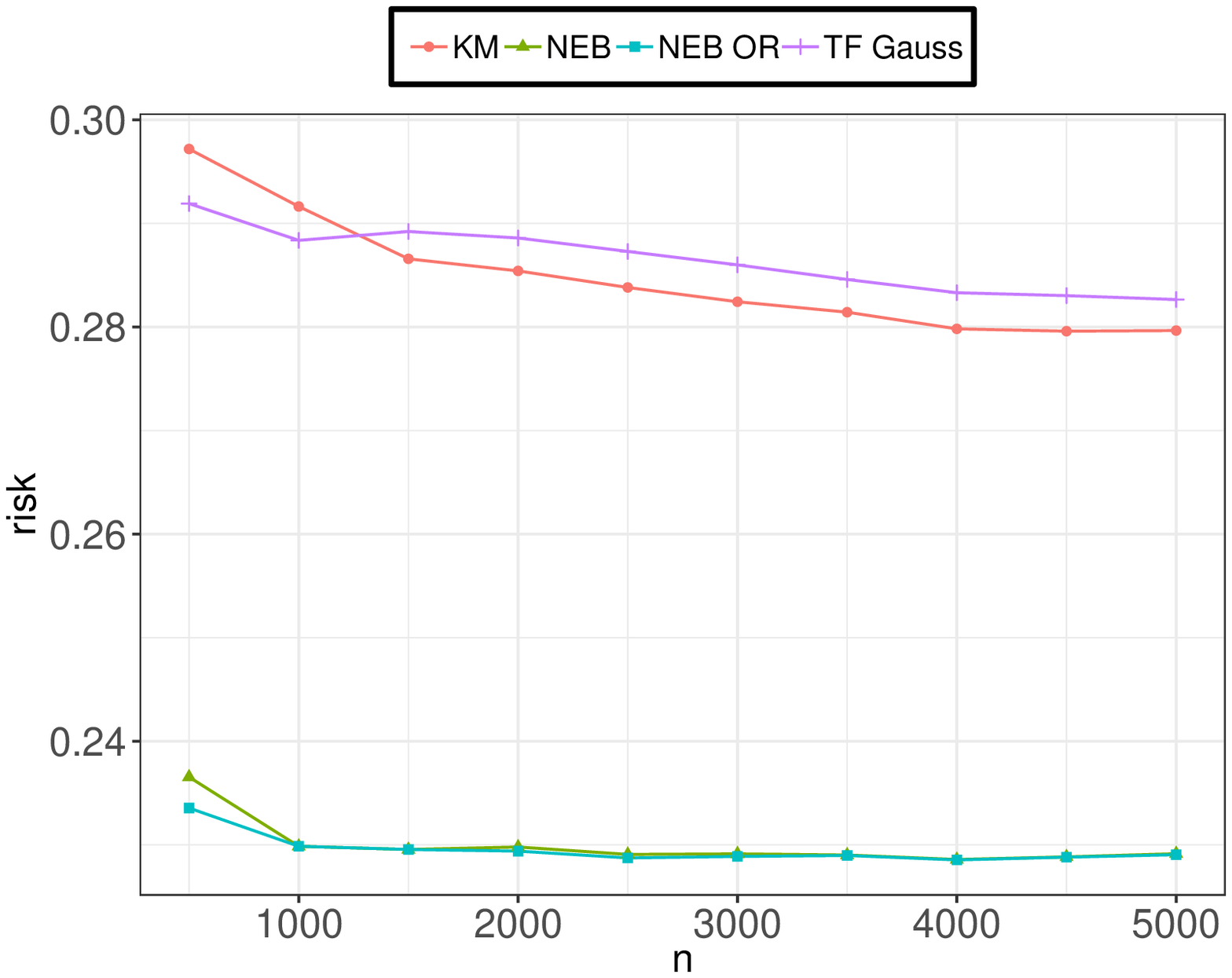}
 		\caption{Scenario 1: Estimation of odds $\bm \theta$ under loss $\mc{L}_n^{(1)}$ where $q_i\stackrel{iid.}{\sim}0.4~\delta_{\{0.5\}}+0.6~\texttt{Beta}(2,5)$ and $m_i=5$.}
 		\label{fig:bin_k1_exp1}
 	\end{subfigure}%
 \hfill
 	\begin{subfigure}{0.46\textwidth}
 		\centering
 		\includegraphics[width=1\linewidth]{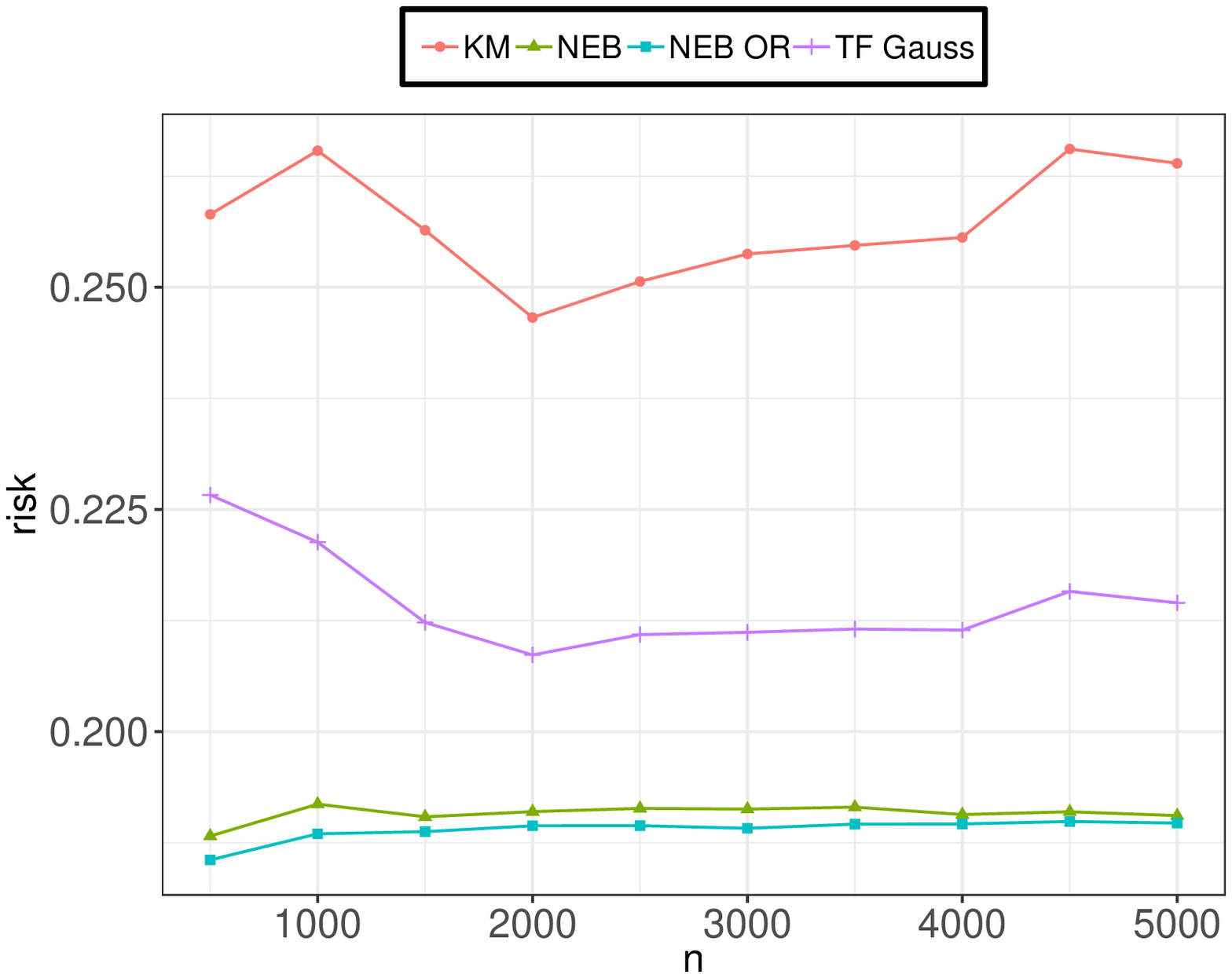}
 		\caption{Scenario 2: Estimation of odds $\bm \theta$ under loss $\mc{L}_n^{(1)}$ where $\theta_i\stackrel{iid.}{\sim}0.8~\delta_{\{0.5\}}+0.2~\texttt{Gamma}(1,2)$ and fix $m_i=10$.}
 		\label{fig:bin_k1_exp2}
 	\end{subfigure}
	\begin{subfigure}{.46\textwidth}
	\centering
	\includegraphics[width=1\linewidth]{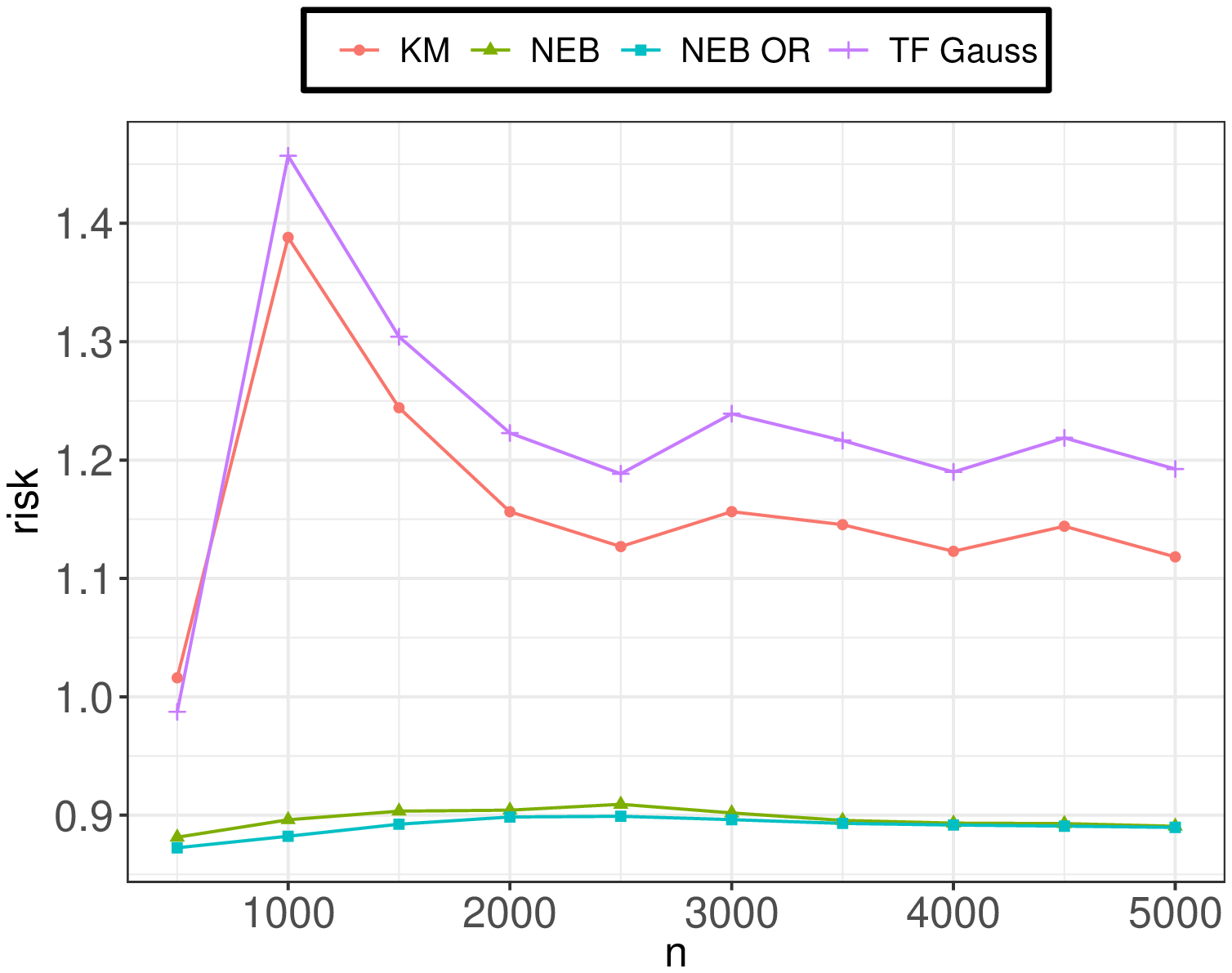}
	\caption{Scenario 3: Estimation of odds $\bm \theta$ under loss $\mc{L}_n^{(1)}$ where $\theta_i\stackrel{iid.}{\sim}\chi^2_2$ and $m_i=5$.}
	\label{fig:bin_k1_exp3}
\end{subfigure}%
\hfill
 \begin{subfigure}{0.46\textwidth}
 	\centering
 	\includegraphics[width=1\linewidth]{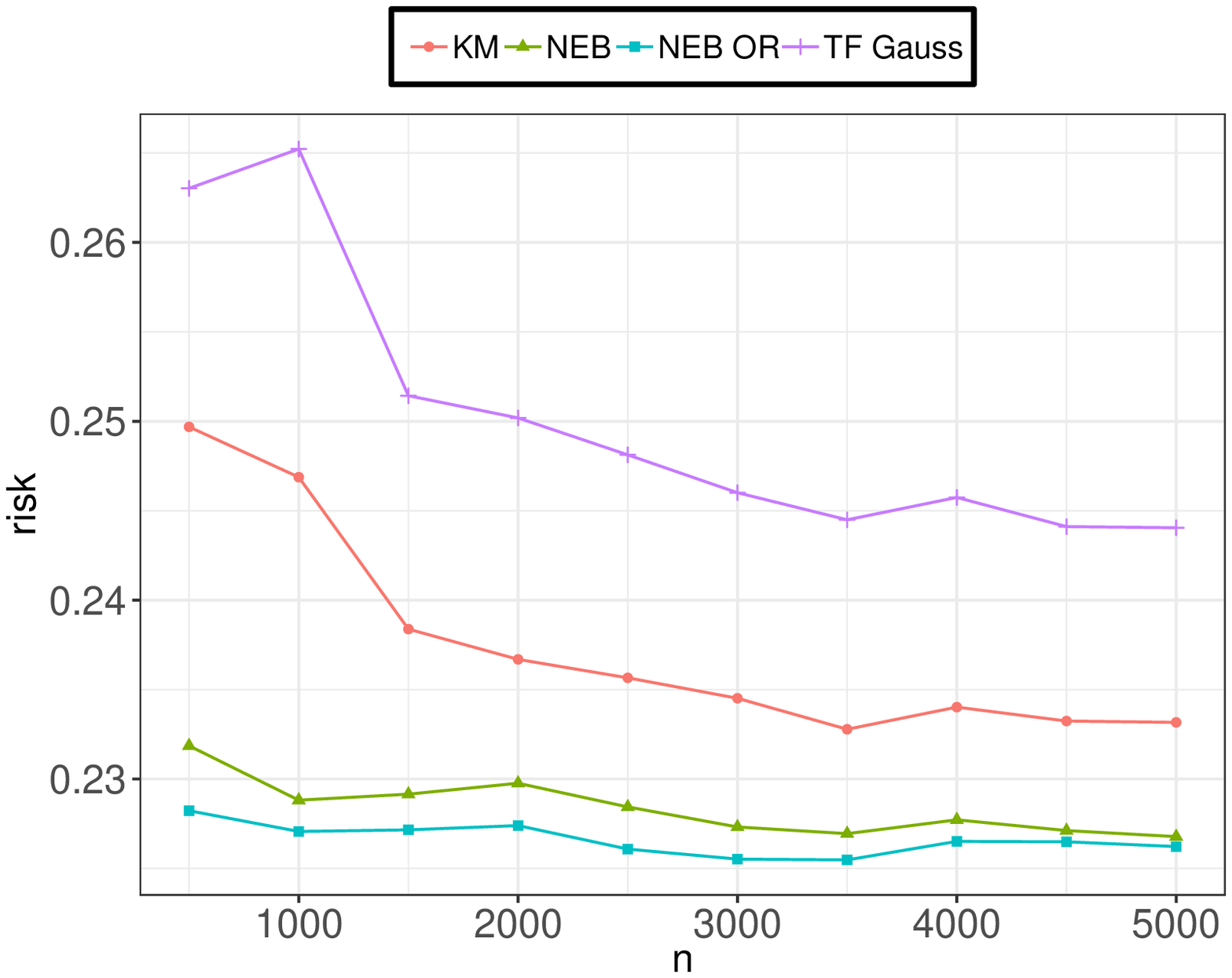}
 	\caption{Scenario 4: Estimation of odds $\bm \theta$ under loss $\mc{L}_n^{(1)}$ where $q_i\stackrel{iid.}{\sim}0.5~\mc{B}eta(1,1)+0.5~\texttt{Beta}(1,3)$ and $m_i=10$.}
 	\label{fig:bin_k1_exp4}
 \end{subfigure}
	\caption{Binomial compound decision problem under scaled squared error loss: Risk estimates of the various estimators for Scenarios 1 to 4.}
\label{fig:bin1}
\end{figure}

\begin{figure}[!t]
	\centering
	\begin{subfigure}{.46\textwidth}
		\centering
		\includegraphics[width=1\linewidth]{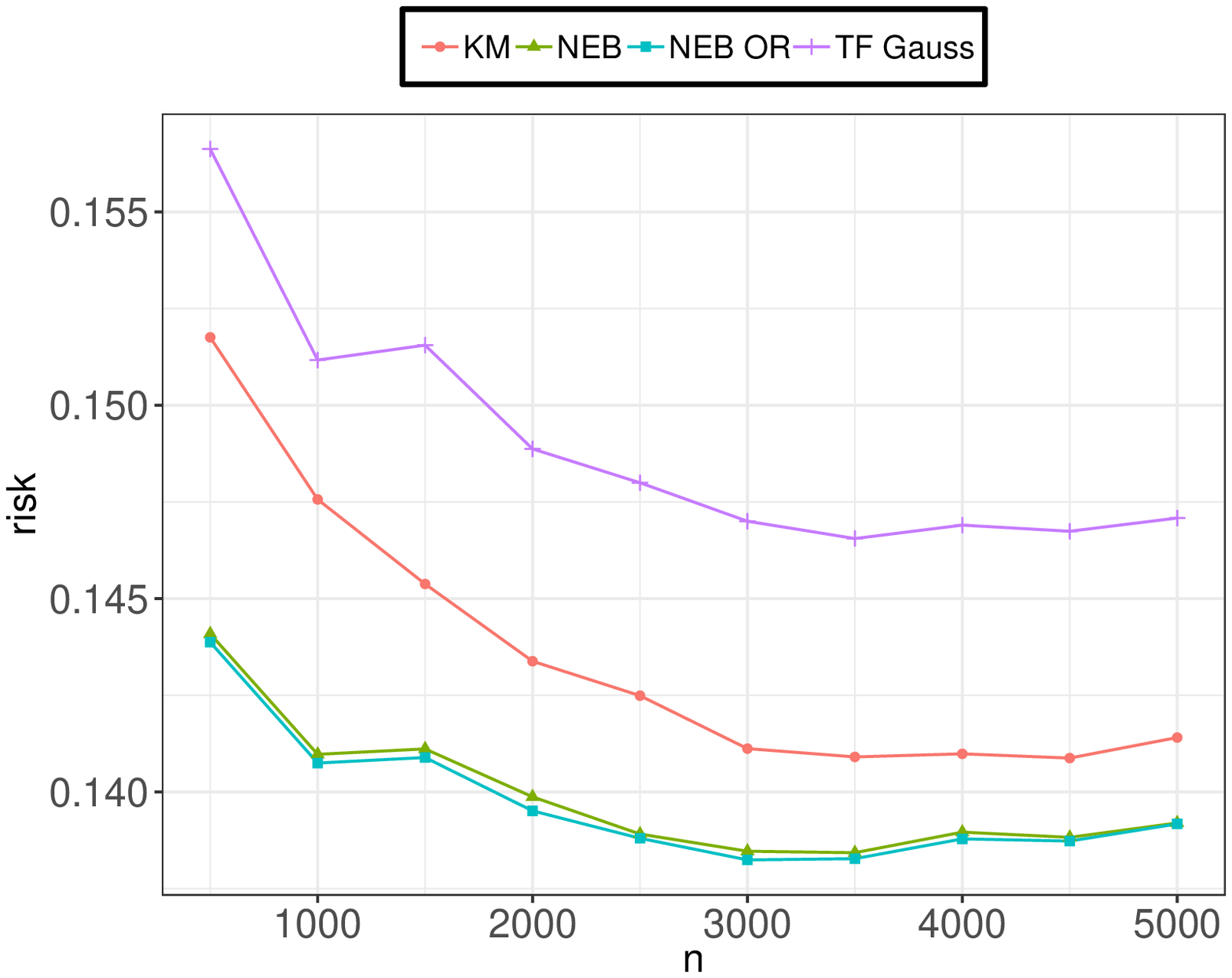}
		\caption{Scenario 1: Estimation of odds $\bm \theta$ under loss $\mc{L}_n^{(0)}$ where $q_i\stackrel{iid.}{\sim}0.4~\delta_{\{0.5\}}+0.6~\texttt{Beta}(2,5)$ and $m_i=5$.}
		\label{fig:bin_k0_exp1}
	\end{subfigure}%
\hfill
	\begin{subfigure}{0.46\textwidth}
		\centering
		\includegraphics[width=1\linewidth]{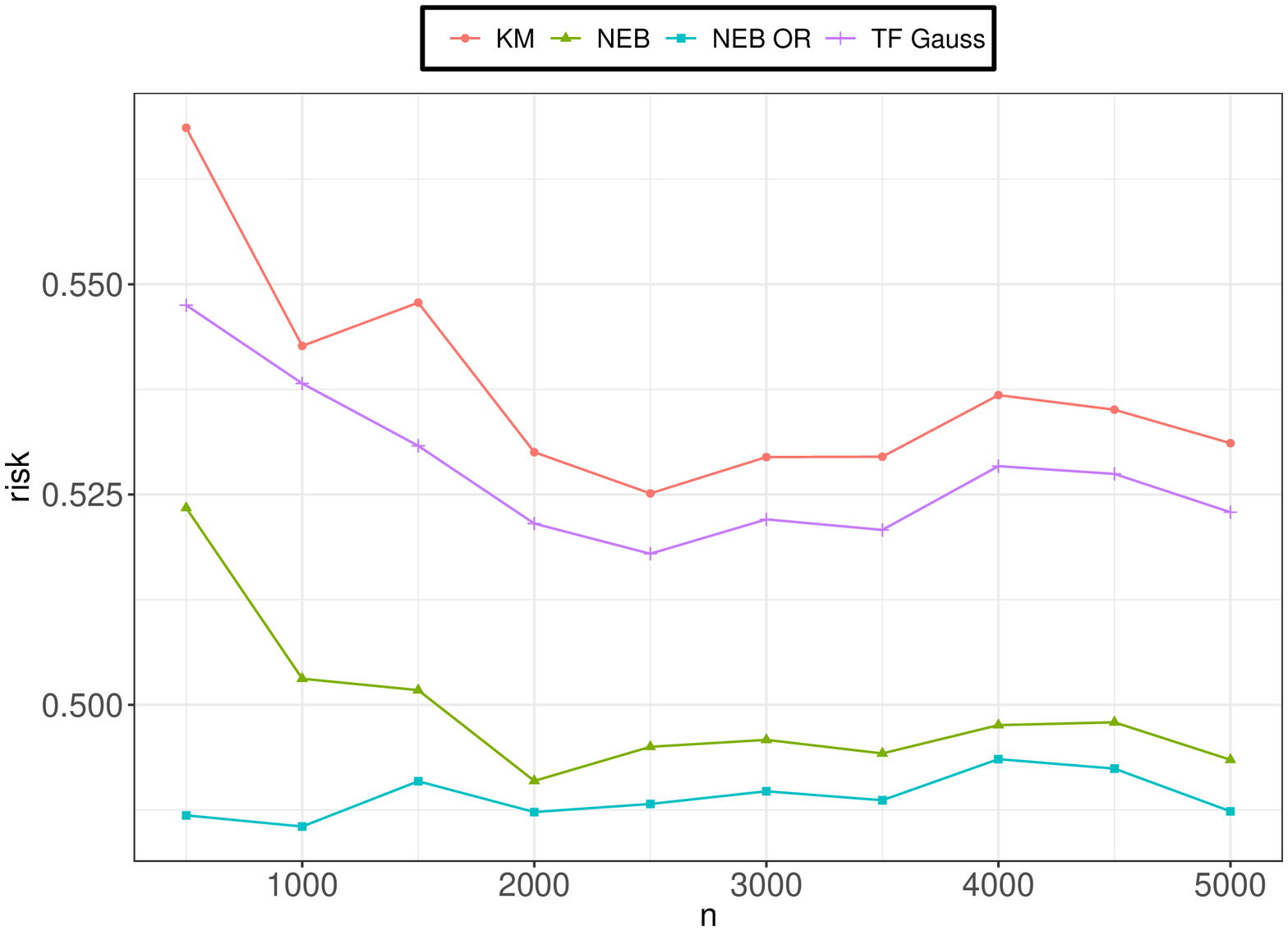}
		\caption{Scenario 2: Estimation of odds $\bm \theta$ under loss $\mc{L}_n^{(0)}$ where $\theta_i\stackrel{iid.}{\sim}0.8~\delta_{\{0.5\}}+0.2~\texttt{Gamma}(1,2)$ and $m_i=10$.}
		\label{fig:bin_k0_exp2}
	\end{subfigure}
	 	\begin{subfigure}{.46\textwidth}
	\centering
	\includegraphics[width=1\linewidth]{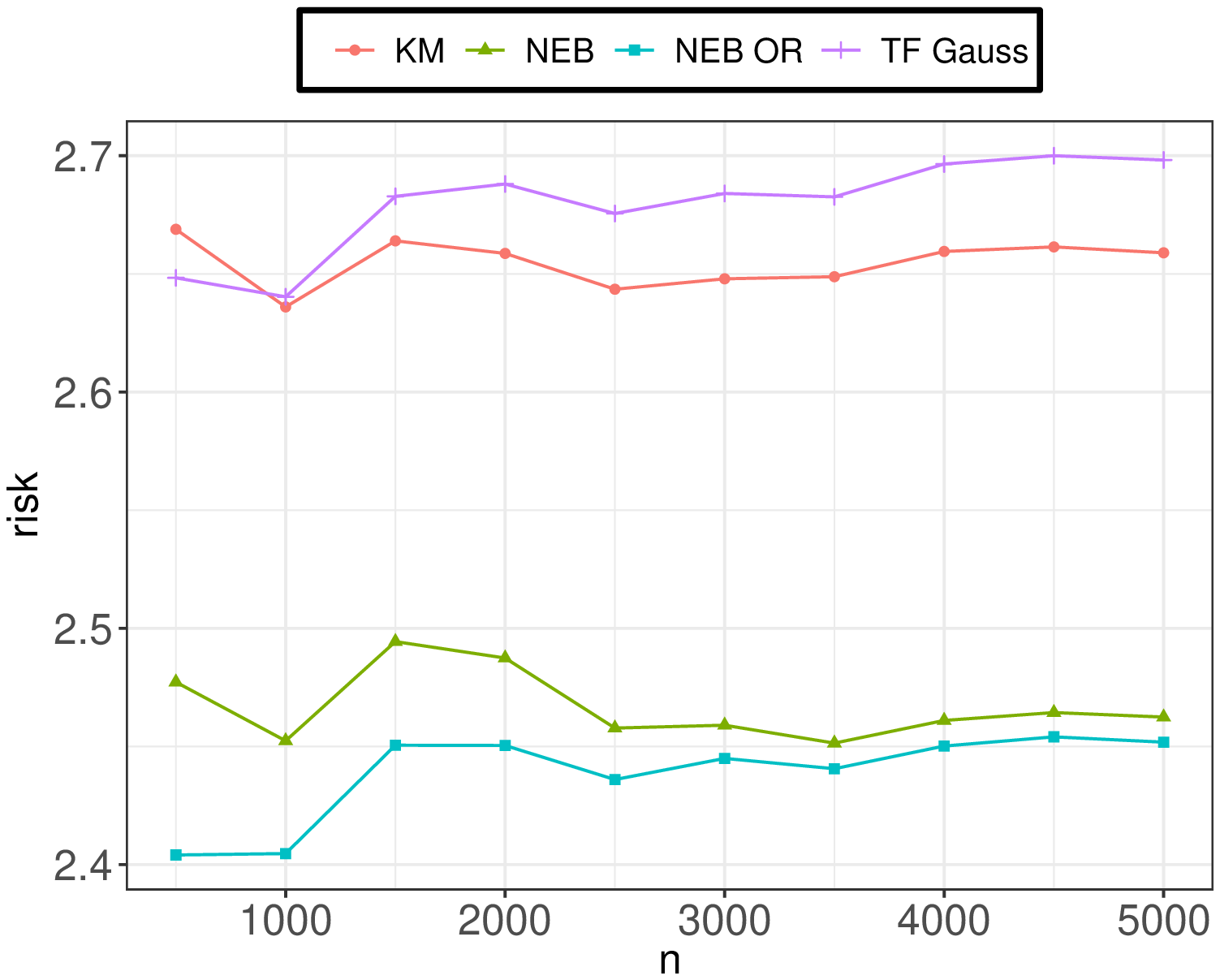}
	\caption{Scenario 3: Estimation of odds $\bm \theta$ under loss $\mc{L}_n^{(0)}$ where $\theta_i\stackrel{iid.}{\sim}\chi^2_2$ and $m_i=5$.}
	\label{fig:bin_k0_exp3}
\end{subfigure}%
	\hfill
	\begin{subfigure}{0.46\textwidth}
		\centering
		\includegraphics[width=1\linewidth]{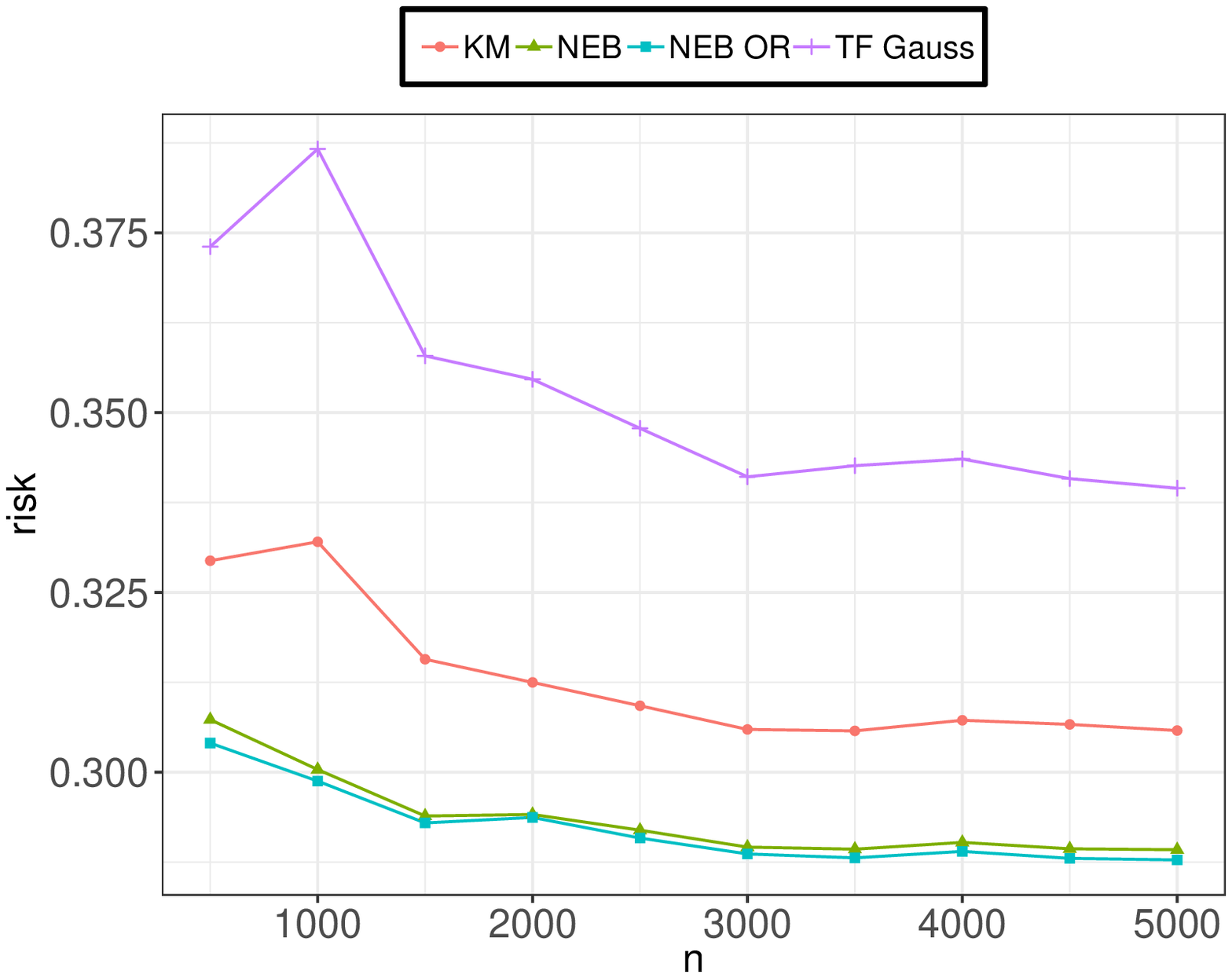}
		\caption{Scenario 4: Estimation of odds $\bm \theta$ under loss $\mc{L}_n^{(0)}$ where $q_i\stackrel{iid.}{\sim}0.5~\mc{B}eta(1,1)+0.5~\texttt{Beta}(1,3)$ and $m_i=10$.}
		\label{fig:bin_k0_exp4}
	\end{subfigure}
	\caption{Binomial compound decision problem under squared error loss: Risk estimates of the various estimators for Scenarios 1 to 4.}
	\label{fig:bin2}
\end{figure}

The simulation results are presented in Figures \ref{fig:bin1} and \ref{fig:bin2} wherein the risks  of various estimators are calculated by averaging over $50$ Monte Carlo repetitions for varying $n$. Tables \ref{tab:binom1} and \ref{tab:binom2} report the risk ratios $\mc{R}_n^{(k)}(\bm \theta,\cdot)/\mc{R}_n^{(k)}(\bm \theta,\bm \delta^{\sf neb}_k)$ at $n=5000$ and for $k=1,0$ respectively, where a risk ratio bigger than $1$ demonstrates a smaller estimation risk for the \texttt{NEB} estimator.

\begin{table}[!h]
	\begin{minipage}{.47\linewidth}
		\centering
		\caption{The Binomial compound decision problem {under scaled squared error loss}: Risk ratios $\mc{R}_n^{(1)}(\bm \theta,\cdot)/\mc{R}_n^{(1)}(\bm \theta,\bm \delta^{\sf neb}_k)$ at $n=5000$ for estimating $\bm \theta$.} 
		\scalebox{0.95}{
			\begin{tabular}{lcccc}
				& \multicolumn{4}{c}{Scenario} \\
				\toprule
				Method & 1     & 2     & 3     & 4 \\
				\midrule
				\texttt{KM}   & 1.22  & 1.38  & 1.26  & 1.03 \\
				\texttt{TF Gauss} & 1.23  & 1.51  & 1.34  & 1.08 \\
				\texttt{TF OR} & $>10$  & $>10$  & $>10$  & $>10$ \\
				\midrule
				\texttt{NEB}   & 1.00  & 1.00  & 1.00  & 1.00 \\
				\texttt{NEB OR} & 1.00  & 1.00  & 1.00  & 1.00 \\
				\bottomrule
		\end{tabular}}%
		\label{tab:binom1}%
	\end{minipage}%
	\hfill
	\begin{minipage}{.47\linewidth}
		\centering
		\caption{The Binomial compound decision problem {under the usual squared error loss}: Risk ratios $\mc{R}_n^{(0)}(\bm \theta,\cdot)/\mc{R}_n^{(0)}(\bm \theta,\bm \delta^{\sf neb}_k)$ at $n=5000$ for estimating $\bm \theta$}
		\scalebox{0.95}{
			\begin{tabular}{lcccc}
				& \multicolumn{4}{c}{Scenario} \\
				\toprule
				Method & 1     & 2     & 3     & 4 \\
				\midrule
				\texttt{KM}   & 1.01  & 1.08  & 1.08  & 1.06 \\
				\texttt{TF Gauss} & 1.06  & 1.06  & 1.09  & 1.17 \\
				\texttt{TF OR} & $>10$  & $>10$  & $>10$  & $>10$ \\
				\midrule
				\texttt{NEB}   & 1.00  & 1.00  & 1.00  & 1.00 \\
				\texttt{NEB OR} & 1.00  & 0.99  & 0.99  & 1.00 \\
				\bottomrule
		\end{tabular}}%
		\label{tab:binom2}%
	\end{minipage} 
\end{table}

We can see from the simulation results that the \texttt{NEB} estimator demonstrates an overall superior risk performance than its competitors. In particular, we see that when estimation is conducted under loss $\mc{L}_n^{(1)}$ the risk ratios of the competing estimators in Table \ref{tab:binom1} reflect a significantly better performance of the \texttt{NEB} estimator. This is not surprising because \texttt{KM, TF Gauss} and \texttt{TF OR} are designed to estimate $\bm q$ under loss $\mc{L}_n^{(0)}$. This also explains the relatively improved performance of these estimators as seen through their risk ratios in table \ref{tab:binom2} wherein the estimation is conducted under the usual squared error loss $\mc{L}_n^{(0)}$. Across the four scenarios, \texttt{TF OR} exhibits the poorest performance and appears to suffer from the fragmented approach of estimating the gradient of the log density $\log p(\bm y)$ wherein $p(\bm y)$ and its first derivative with respect to $\bm y$ are estimated separately using a Gaussian kernel with common bandwidth $h^{\sf orc}$. The approach of using a variance stabilizing transformation to convert the data to approximate normality renders \texttt{TF Gauss} highly competitive to \texttt{KM} \citep{koenker2014convex}. A similar phenomenon was also reported in \cite{brown2013poisson} in the context of the Poisson model. However, under the Binomial model, when the primary goal is to estimate the odds $\bm \theta$, the risk ratios reported in Tables \ref{tab:binom1} and \ref{tab:binom2} suggest that the proposed \texttt{NEB} estimator is by far the best amongst these competitors under both types of losses.

\section{Real Data Analyses}
\label{sec:realdata}

This section illustrates the proposed method for estimating the Juvenile Delinquency rates from Poisson models and news popularity in Binomial models.

\subsection{Estimation of Juvenile Delinquency rates}
\label{sec:realdata_juvenile}
In this section we consider an application for analysis of the Uniform Crime Reporting Program (UCRP) Database \citep{us2012uniform} that holds county-level counts of arrests and offenses ranging from robbery to weapons violations in 2012. The database is maintained by the National Archive of Criminal Justice Data (NACJD) and is one of the most widely used database for research related to factors that affect juvenile delinquency (JD) rates across the United States; see for example \cite{aizer2015juvenile} and \cite{damm2014does,koski2018state}. A preliminary and important goal in these analyses is to estimate the JD rates based on the observed arrest data and determine the counties that are amongst the worst or least affected. However with almost 3,000 counties being evaluated the JD rates are susceptible to selection bias, wherein some of the data points are in the extremes merely by chance and traditional estimators may underestimate or overestimate the corresponding means, especially in counties with fewer total number of arrests across all age groups. 
 \begin{figure}[!t]
	\centering
	\includegraphics[width=0.7\linewidth]{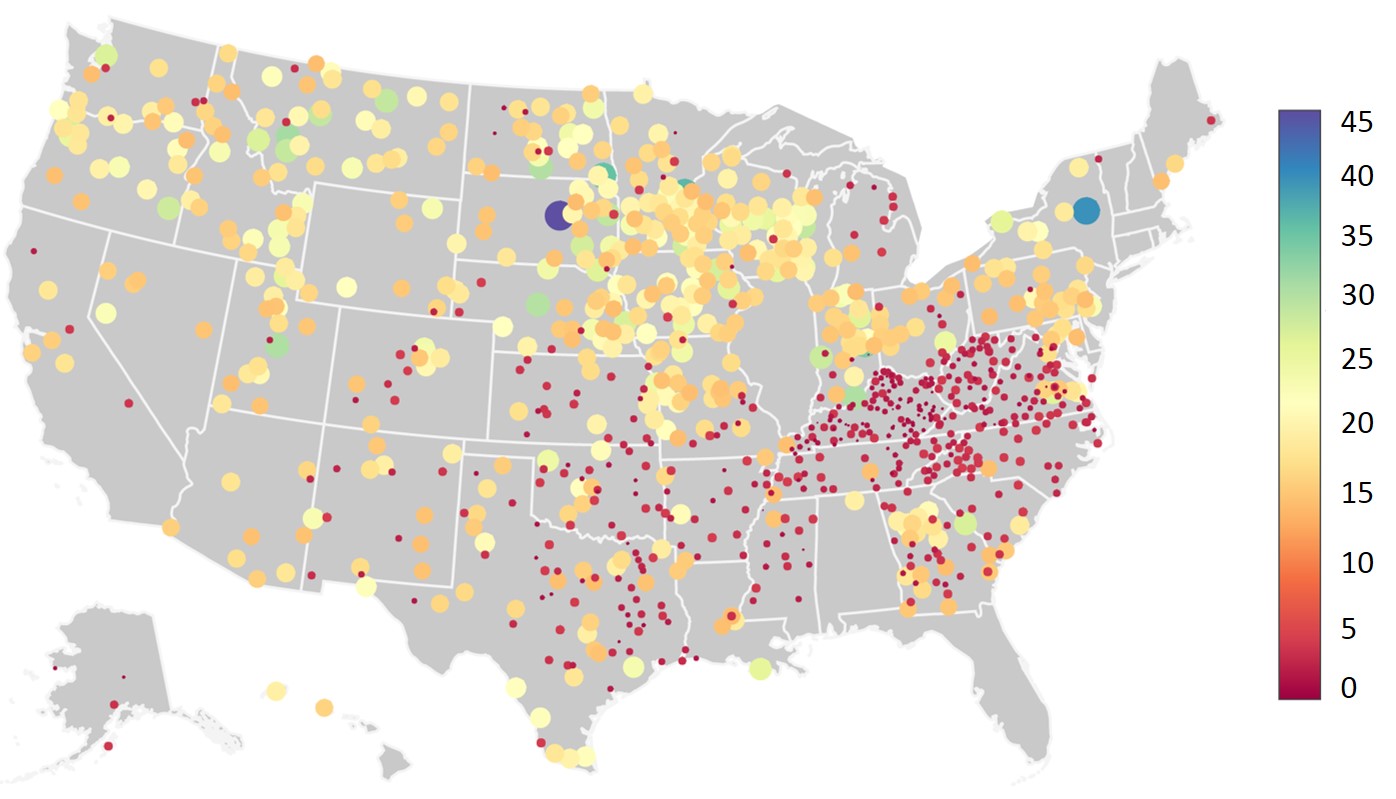}
	\caption{Observed Juvenile Delinquency rates in 2012. The top $500$ and bottom $500$ counties are plotted. The data on Florida arrests is not available in the \cite{us2012uniform} database.}
	\label{fig:crime1}

	\includegraphics[width=1.01\linewidth]{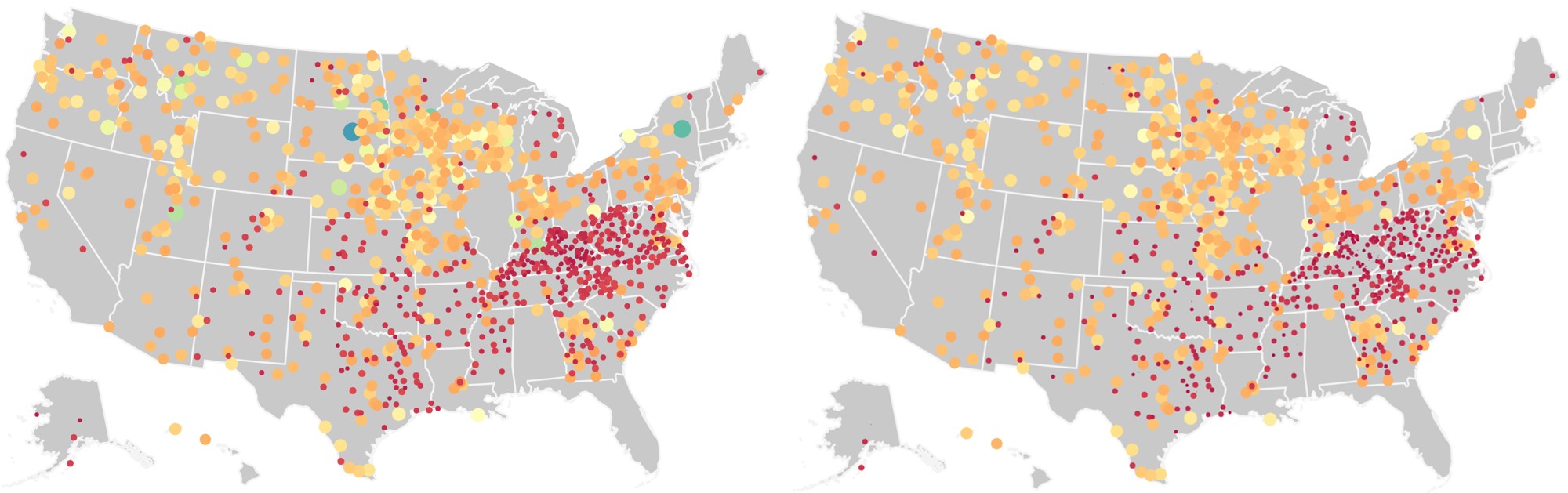}
	\caption{Estimated Juvenile Delinquency rates of the $1000$ counties exhibited in figure \ref{fig:crime1}. Left: Estimation under squared error loss $(k=0)$. Right: Estimation under scaled squared error loss $(k=1)$.}
	\label{fig:crime2}
\end{figure}
\begin{table}[!h]
	\centering
	\caption{Loss ratios of the competing methods for estimating $\bm \theta$.}
	\begin{tabular}{lcc}
		\toprule
		$(n=2,803)$& \multicolumn{2}{c}{Loss  ratios} \\
		\midrule
		Method & $k=1$   & $k=0$ \\
		\midrule
		\texttt{NEB}  & 1.00  & 1.00 \\
		\texttt{BGR}   & 1.18  & 1.03 \\
		\texttt{KM}  & 1.19  & 1.03 \\
		\texttt{TF Gauss} & 1.12  & 1.01 \\
		\texttt{TF OR} & 1.11  & 1.05 \\
		\bottomrule
	\end{tabular}%
	\label{tab:crime}%
\end{table}%

For the purpose of our analyses, we use the 2012 UCRP data that spans $n=3,178$ counties in the U.S. and consider estimating the mean JD rate $\bm \theta$ as a vector of Poisson means. The observed data for county $i$ is denoted $y_i$, which represents the number of juvenile arrests expressed as a percentage of total arrests in that county in Year 2012. We assume that $Y_i~|~\theta_i\stackrel{ind.}{\sim}\texttt{Poi}(\theta_i)$ for $i=1,\ldots,n$. Figure \ref{fig:crime1} plots the observed data for the top $500$ and the bottom $500$ counties that have at least $1$ juvenile arrest. Campbell county in South Dakota, followed by Fulton county in New York, exhibits the highest observed JD rates in Year 2012. 

As discussed in Section \ref{sec:sims}, we consider the following 5 estimators of $\bm \theta$: \texttt{NEB}, \texttt{BGR} \citep{brown2013poisson}, \texttt{KM} \citep{koenker2017rebayes}, \texttt{TF OR} \citep{efron2011tweedie} and \texttt{TF Gauss} \citep{koenker2014convex,brown2013poisson}. We use the 2014 UCRP data \citep{us2014uniform} to compare their estimation accuracies under both $\mc{L}_n^{(0)}$ and $\mc{L}_n^{(1)}$ losses. The data were cleaned prior to any analyses which ensured that all counties in the year 2012 had at least one arrest (juvenile or not). This resulted in $n=2803$ counties where all methods are applied to. 

In Figure \ref{fig:crime2} we visualize the shrinkage estimates of JD rates for those $1000$ counties considered in Figure \ref{fig:crime1}. The left plot presents the estimates under the squared error loss while right plot presents the results under the scaled squared error loss. Notably, the scaled error loss exhibits a larger magnitude of shrinkage for the bigger observations than the squared error loss. 
Table \ref{tab:crime} reports the loss ratios $\mc{L}_n^{(k)}(\bm \theta,\bm \delta)/\mc{L}_n^{(k)}(\bm \theta,{\bm \delta}^{\sf neb}_{(k)})$ where for any estimator $\bm \delta$ of $\bm\theta$, a ratio bigger than $1$ indicates a smaller estimation loss for ${\bm \delta}^{\sf neb}_{(k)}$. We can see that for estimating $\bm \theta$, all four competitors exhibit loss ratios bigger than $1$ under the scaled squared error loss ($k=1$). This is not surprising since these competitors are designed to estimate $\bm \theta$ under the regular squared error loss ($k=0$). Interestingly, even under the regular loss, the \texttt{NEB} estimator continues to provide a better estimation accuracy than \texttt{TF OR}, \texttt{BGR} and \texttt{KM}, and demonstrates a competitive performance against \texttt{TF Gauss}.

\subsection{News popularity in social media platforms}
\label{sec:realdata_news}

Journalists and editors often face the critical task of assessing the popularity of various news items and determining which articles are likely to become popular; hence existing content generation resources can be efficiently managed and optimally allocated to avenues with maximum potential. Due to the dynamic nature of the news articles, popularity is usually measured by how quickly the article propagates (frequency) and the number of readers that the article can reach (severity) through social media platforms like Twitter, Youtube, Facebook and LinkedIn. As such predicting these two aspects of popularity based on early trends is extremely valuable to journalists and content generators \citep{bandari2012pulse}. 

%
 \begin{table}[!h]
	\begin{minipage}{.47\linewidth}
		\centering
		\caption{Loss ratios of the competing methods for estimating $\bm \theta$. News article genre: \texttt{Economy} and social media: \texttt{Facebook}}
		\scalebox{0.95}{
	\begin{tabular}{lcc}
		\toprule
		$(n=3,972)$ & \multicolumn{2}{c}{Loss Ratios} \\
		\midrule
		Method & $k=1$   & $k=0$ \\
		\midrule
		\texttt{NEB}   & 1.00  & 1.00 \\
		\texttt{KM}  & 12.25  & 13.25 \\
		\texttt{TF Gauss} & 4.06  & 3.25 \\
		\texttt{TF OR} & 81.57 & 36.24 \\
		\bottomrule
	\end{tabular}}%
	\label{tab:news_economy}%
	\end{minipage}%
	\hfill
	\begin{minipage}{.47\linewidth}
		\centering
		\caption{Loss ratios of the competing methods for estimating $\bm \theta$. News article genre: \texttt{Microsoft} and social media: \texttt{LinkedIn}}
		\scalebox{0.95}{
		\begin{tabular}{lcc}
			\toprule
		$(n=3,850)$	& \multicolumn{2}{c}{Loss Ratios} \\
			\midrule
			Method & $k=1$   & $k=0$ \\
			\midrule
			\texttt{NEB}   & 1.00  & 1.00 \\
			\texttt{KM}  & 41.06  & 59.64 \\
			\texttt{TF Gauss} & 9.31  & 7.44 \\
			\texttt{TF OR} & 36.14  & 12.76 \\
			\bottomrule
		\end{tabular}}%
		\label{tab:news_micro}%
	\end{minipage} 
\end{table}

In this section, we assess the popularity of several news items based on their frequency of propagation and analyze a dataset from \cite{moniz2018multi} that holds $48$ hours worth of social media feedback data on a large collection of news articles since the time of first publication. For the purposes of our analysis, we consider two popular genres of news from this data set: \texttt{Economy} and \texttt{Microsoft}, and examine how frequently these articles were shared in Facebook and LinkedIn, respectively, over a period of $48$ hours from the time of their first publication. Each news article in the data has a unique identifier and $16$ consecutive time intervals, each of length $180$ minutes, to detect whether the article was shared at least once in that time interval. Let $Z_{ij}=1$ if article $i$ was shared in time interval $j$ and $0$ otherwise, where $i=1,\ldots,n$ and $j=1,\ldots, 16$. Suppose $q_{ij}\in[0,1]$ denote the probability that news article $i$ is shared in interval $j$. We let $q_{ij}=q_i$ for all $j=1,\ldots, 16$ and assume that for each $i$, $Z_{ij}$ are independent realizations from $\texttt{Ber}(q_i)$. It follows that $Y_i=\sum_{j=1}^{8}Z_{ij}\stackrel{ind.}{\sim}\texttt{Bin}(8,q_i)$. 

To assess the popularity of article $i$, we estimate its odds of sharing given by $\theta_i=q_i/(1-q_i)$ and consider the following 4 estimators of $\bm \theta$: \texttt{NEB}, \texttt{KM} \citep{koenker2017rebayes}, \texttt{TF Gauss} \citep{koenker2014convex} and \texttt{TF OR} \citep{efron2011tweedie,funonparametric2018}. We use the data on time points $j=9,\ldots,16$ to compare the estimation accuracy of these estimators under both $\mc{L}_n^{(0)}$ and $\mc{L}_n^{(1)}$ losses. Tables \ref{tab:news_economy} and \ref{tab:news_micro} report the loss ratios $\mc{L}_n^{(k)}(\bm \theta,\bm \delta)/\mc{L}_n^{(k)}(\bm \theta,{\bm \delta}^{\sf neb}_{(k)})$ where for any estimator $\bm \delta$ of $\bm\theta$, a ratio bigger than $1$ indicates a smaller estimation loss for ${\bm \delta}^{\sf neb}_{(k)}$. We observe that the three competitors to the \texttt{NEB} estimator exhibit loss ratios substantially bigger than $1$ under both the losses. This is not surprising since these competitors are designed to estimate $\bm q$ and $\log \bm \theta$ under a squared error loss ($k=0$). 
However, when the primary goal is to estimate the odds, the proposed \texttt{NEB} estimator is by far the best amongst these competitors under both losses.

\section*{Acknowledgments}

G. Mukherjee's work was partly supported by the NSF grant DMS-1811866. W. Sun's work was partly supported by the NSF grant DMS-1712983.
Q. Liu's work is supported in part by NSF CRII 1830161 and NSF CAREER 1846421.

\newpage

\begin{appendices}

\begin{center}\LARGE
		Supplementary Material for ``EB Estimation in Discrete Linear Exponential Family''
	\end{center}

In this supplement, we first present in Appendix \ref{l0:sec} the results for the \texttt{NEB} estimator under the regular squared loss, then provide in Appendix \ref{sec:proofs} the proofs and technical details of all theories in the main text and Appendix \ref{l0:sec}.


\section{Results Under the Squared Error Loss}\label{l0:sec}


\subsection{The \texttt{NEB} estimator}
\label{sec:neb_sqaurederror}

In this section we discuss the estimation of $\bm w_p^{(0)}$ that appear in lemma \ref{lem:bayes} under the usual squared error loss $(k=0)$. 
Let $Y$ be a non-negative integer-valued random variable with probability mass function (pmf) $p$ and define
\begin{equation}
\label{eq:h(y)_k0}
h_0^{(0)}(y)=\dfrac{y+1}{w_p^{(0)}(y)}-y~,y\in\{0\}\cup\mathbb{N}
\end{equation}
Suppose $\mc{K}_\lambda(y, y')=\exp\{-0.5\lambda^{-2}(y-y')^2\}$ be the positive definite RBF kernel with bandwidth parameter $\lambda\in\Lambda$ where $\Lambda$ is a compact subset of $\mathbb{R}^+$ bounded away from $0$. 
Given observations $\bm y=(y_1,\ldots,y_n)$ from model \eqref{eq:DLEmodel}, let $\bm h_0^{(0)} = (h_{0}^{(0)}(y_1),\ldots,h_{0}^{(0)}(y_n))$ and define the following $n\times n$ matrices: $n^2\bm K_{\lambda} = [\mc{K}_\lambda(y_i,y_j)]_{ij}$, $n^2{\Delta \bm K}_{\lambda}  = [\Delta_{y_i} \mc{K}_\lambda(y_i, y_j+1)]_{ij}$ and $n^2{\Delta_2 \bm K}_{\lambda}  = [\Delta_{y_i,y_j}\mc{K}_\lambda(y_i, y_j)]_{ij}$ where $\Delta_{y} \mc{K}_\lambda(y,y')=\mc{K}_\lambda(y+1,y')-\mc{K}_\lambda(y,y')$ and $\Delta_{y,y'} \mc{K}_\lambda(y,y')=\Delta_{y'}\Delta_{y} \mc{K}_\lambda(y,y')=\Delta_{y}\Delta_{y'} \mc{K}_\lambda(y,y')$. 
\renewcommand*{\thedefinition}{1B} 
\begin{definition}[\texttt{NEB} estimator of $\theta_i$] 
	\label{def:neb_k0}	
	Consider the \texttt{DLE} Model \eqref{eq:DLEmodel} with loss $\ell^{(0)}(\theta_i,\delta_i)$. For a fixed $\lambda\in\Lambda$, let $\hat{w}_{i}^{(0)}(\lambda)=(y_i+1)/(y_i+\hat{h}_{i}^{(0)}(\lambda))$ and $\hat{\bm h}_n^{(0)}(\lambda)=\left\{\hat{h}_{1}^{(0)}(\lambda),\ldots,\hat{h}_{n}^{(0)}(\lambda)\right\}$ be the solution to the following quadratic optimization problem:
	\begin{eqnarray}
	\label{eq:quad_opt_k0}
	&&\min_{\bm h\in \bm H_n}~\hat{\mathbb{M}}_{\lambda,n}(\bm h)=\bm h^T\bm K_\lambda\bm h+2\bm h^T\Delta \bm K_\lambda\bm y+\bm y^T\Delta_2\bm K_\lambda\bm y,
	\end{eqnarray}	
	where $\bm H_n=\{\bm h=(h_1,\ldots,h_n):\mc{A}\bm h\preceq \bm b,~\mc{C}\bm h= \bm d\}$ is a convex set and $\mc{A},\mc{C},\bm b,\bm d$ are known real matrices and vectors that enforce linear constraints on the components of $\bm h$. Then the \texttt{NEB} estimator for a fixed $\lambda$ is given by ${\bm\delta}^{\sf neb}_{(0)} (\lambda)=\left\{{\delta}_{(0),i}^{\sf neb}(\lambda): 1\leq i\leq n\right\}$, where
	$$
	{\delta}_{(0),i}^{\sf neb}(\lambda)=
	\dfrac{a_{y_i}/a_{y_i+1}}{\hat{w}_{i}^{(0)}(\lambda)}, \quad \mbox{ if } y_i\in\{0,1, 2, \ldots\}
	$$
\end{definition}

\renewcommand*{\thetheorem}{2B} 
\begin{theorem}
	\label{thm:w_k0}
		Let $\mc{K}_\lambda(\cdot,\cdot)$ be the positive definite RBF kernel with bandwidth parameter $\lambda\in\Lambda$. If $\lim_{n\to\infty}c_nn^{-1/2}\log n=0$ then, under assumptions {$A1 - A3$}, we have for any $\lambda\in\Lambda$,
	$$\lim_{n\to\infty}\mathbb{P}\Big[\Big\|\hat{\bm w}_n^{(0)}(\lambda)-{\bm w}_p^{(0)}\Big\|_2\ge c_n^{-1}\epsilon\Big]=0,~\text{for any }\epsilon>0 $$
	where $\hat{\bm w}_n^{(0)}(\lambda)=[(Y_i+1)/(\hat{h}_i^{(0)}(\lambda)+Y_i)]_i$. 
\end{theorem}
We now provide some motivation behind the minimization
problem in definition \ref{def:neb_k0} for estimating the ratio functionals ${\bm w}_p^{(0)}$. Suppose $\tilde{p}$ be a probability mass function on the support of $Y$ and define
\begin{equation}
\label{eq:ksd_k0}
\mc{S}_\lambda[\tilde{p}](p)=
\mathbb{E}_p\Big[(\tilde{h}^{(0)}(Y)-h_0^{(0)}(Y))\mc{K}_\lambda(Y+1,Y'+1)(\tilde{h}^{(0)}(Y')-h_0^{(0)}(Y'))\Big]
\end{equation}
where $h_0^{(0)},\tilde{h}^{(0)}$ are as defined in equation \eqref{eq:h(y)_k0} and $Y,Y'$ are i.i.d copies from the marginal distribution that has mass function $p$. $\mc{S}_\lambda[\tilde{p}](p)$ in equation \eqref{eq:ksd_k0} is the Kernelized Stein's Discrepancy (KSD) measure that can be used to distinguish between two distributions with mass functions $p,\tilde{p}$ such that $\mc{S}_\lambda[\tilde{p}](p)\ge0$ and $\mc{S}_\lambda[\tilde{p}](p)=0$ if and only if $p=\tilde{p}$ \citep{liu2016kernelized,chwialkowski2016kernel}. Moreover $\mc{S}_\lambda[\tilde{p}](p)=\mathbb{E}_{(U,V)\stackrel{i.i.d}{\sim}}p\Big[\kappa_\lambda[\tilde{p}](U,V)\Big]$ where $\kappa_\lambda[\tilde{p}](u,v)$ is
\begin{equation}
\label{eq:kappa_k0}
\tilde{h}^{(0)}(u)\tilde{h}^{(0)}(v) \mc{K}_\lambda(u, v)   
+    \tilde{h}^{(0)}(u)v \Delta_{v} \mc{K}_\lambda(u+1,v) 
+  \tilde{h}^{(0)}(v)u \Delta_{u} \mc{K}_\lambda(u,v+1)
+  uv\Delta_{u,v} \mc{K}_\lambda(u,v)
\end{equation}
An empirical evaluation scheme for $\mc{S}_\lambda[\tilde{p}](p)$ is given by $\mc{S}_\lambda[\tilde{p}](\hat{p}_n)$ where
\begin{equation}
\label{eq:est_sqp_k0}
\mc{S}_\lambda[\tilde{p}](\hat{p}_n)=\dfrac{1}{n^2}\sum_{i,j=1}^{n}\kappa_\lambda[\tilde{p}](y_i,y_j)
\end{equation}
and $\bm Y=(Y_1,\ldots,Y_n)$ is a random sample from the marginal distribution with mass function $p$ with empirical CDF $\hat{p}_n$.
Note that $\kappa_{\lambda}[\tilde{p}](u,v)$ in equation \eqref{eq:kappa_k0} involves $\tilde{p}$ only through $\tilde{h}^{(0)}$ and may analogously be denoted by $\kappa_{\lambda}[\tilde{h}(u),\tilde{h}(v)](u,v)\coloneqq \kappa_{\lambda}[\tilde{p}](u,v)$ where we have dropped the superscript from $\tilde{h}$ that indicates that the loss in question is the regular squared error loss. This slight abuse of notation is harmless as the discussion in this section is geared towards the squared error loss only.

Under the empirical Bayes compound estimation framework of model \eqref{eq:DLEmodel}, our goal is to estimate $h_0^{(0)}$. To do that we minimize $\mc{S}_\lambda[\tilde{p}](\hat{p}_n)$ in equation \eqref{eq:est_sqp_k0} with respect to the unknowns $\bm \tilde{h}=(\tilde{h}(y_1),\ldots,\tilde{h}(y_n))$ and the sample criteria
$$\hat{\mathbb{M}}_{\lambda,n}(\tilde{\bm h})=\dfrac{1}{n^2}\sum_{i,j=1}^{n}\kappa_{\lambda}[\bm \tilde{h}(y_i),\tilde{h}(y_j)](y_i,y_j)$$ is thus the objective function of the optimisation problem in equation \eqref{eq:quad_opt_k0} with optimisation variables $h_i\equiv \tilde{h}(y_i)$ for $i=1,\ldots,n$.

Note that $\hat{\mathbb{M}}_{\lambda,n}(\tilde{\bm h})$ above is a V-statistic and a biased estimator of the population criteria $\mathbb{M}_\lambda(\tilde{\bm h})$ that is defined in equation \eqref{eq:sqp} with $\kappa_{\lambda}[\tilde{h}(u),\tilde{h}(v)](u,v)$ given by equation \eqref{eq:kappa_k0}.

\subsection{Bandwidth choice and asymptotic properties}
\label{sec:lam_k0}
We propose the following asymptotic risk estimate of the true risk of ${\bm \delta}_{(0)}^{\sf neb}(\lambda)$ in the Poisson and Binomial model. 
\renewcommand*{\thedefinition}{2B}
\begin{definition}[ARE of ${\bm \delta}_{(0)}^{\sf neb}(\lambda)$ in the Poisson model]
	\label{def:are_pois_k0}	
	Suppose $Y_i~|~\theta_i\stackrel{ind.}{\sim} \texttt{Pois}(\theta_i)$. Under the loss $\ell^{(0)}(\theta_i,\cdot)$ an ARE of the true risk of ${\bm \delta}_{(0)}^{\sf neb}(\lambda)$ is
	$${\sf ARE}^{(0,\mc{P})}_n(\lambda)=\dfrac{1}{n}\Bigl\{ \sum_{i=1}^{n}y_i(y_i-1)-2\sum_{i=1}^{n}y_i{\psi}_\lambda(y_i)+\sum_{i=1}^{n}\{{\delta}_{(0),i}^{\sf neb}(\lambda)\}^2\Bigr\}
	$$	where $$
	\psi_\lambda(y_i)={\delta}_{(0),j}^{\sf neb}(\lambda),~y_i\ge 1
	$$ with $j\in\{1,\ldots,n\}$ such that $y_j=y_i-1$.
\end{definition}
\renewcommand*{\thedefinition}{3B}
\begin{definition}[ARE of ${\bm \delta}_{(0)}^{\sf neb}(\lambda)$ in the Binomial model]
	\label{def:are_bin_k0}	
	Suppose $Y_i~|~q_i\sim \texttt{Bin}(m,q_i)$ so that in equation \eqref{eq:DLE} $a_{y_i}={m \choose y_i}$ and $\theta_i=q_i/(1-q_i)$ in equation \eqref{eq:DLEmodel}. Under the loss $\ell^{(0)}(\theta_i,\cdot)$ an ARE of the true risk of ${\bm \delta}_{(0)}^{\sf neb}(\lambda)$ is
	$${\sf ARE}^{(0,\mc{B})}_n(\lambda)=\dfrac{1}{n}\Bigl\{ \sum_{i=1}^{n}\dfrac{y_i(y_i-1)}{(m-y_i+2)(m-y_i+1)}-2\sum_{i=1}^{n}y_i{\psi}_\lambda(y_i)+\sum_{i=1}^{n}\{{\delta}_{(0),i}^{\sf neb}(\lambda)\}^2\Bigr\}
	$$	where $$
	\psi_\lambda(y_i)={\delta}_{(0),j}^{\sf neb}(\lambda)/(m-y_i+1),~y_i\ge 1
	$$ with $j\in\{1,\ldots,n\}$ such that $y_j=y_i-1$.
\end{definition}
Note that if for some index $i$, $y_i-1$ is not available in the observed sample $\bm y$, ${\psi}_\lambda(y_i)$ can be calculated using cubic splines. 
We propose the following estimate of the tuning parameter $\lambda$ based on the ARE:
\begin{equation}
\label{eq:lam_k0}
\hat{\lambda}=\begin{cases}
\argmin_{\lambda\in{\Lambda}}{\sf ARE}^{(0,\mc{P})}_n(\lambda),~\text{if~}Y_i~|~\theta_i\stackrel{ind.}{\sim} \texttt{Pois}(\theta_i)\\
\argmin_{\lambda\in{\Lambda}}{\sf ARE}^{(0,\mc{B})}_n(\lambda),~\text{if~}Y_i~|~q_i\stackrel{ind.}{\sim} \texttt{Bin}(m,q_i)
\end{cases}
\end{equation}
where a choice of $\Lambda =[10,10^{2}]$  work well in the simulations and real data analyses of sections \ref{sec:numresults} and \ref{sec:realdata}. Lemmata \ref{lem:lam_pois_k1} and \ref{lem:lam_bin_k1} continue to provide the large-sample properties of the proposed ${\sf ARE}^{(0,\mc{P})}_n$, ${\sf ARE}^{(0,\mc{B})}_n$ criteria.

To analyze the quality of the estimates $\hat{\lambda}$ obtained from equation \eqref{eq:lam_k0}, we consider an oracle loss estimator $\bm \delta^{\sf or}_{(0)}\coloneqq{\bm \delta}_{(0)}^{\sf neb}(\lambda^{\sf orc}_0)$ where $$\lambda^{\sf orc}_0=\argmin_{\lambda\in\Lambda}\mc{L}_n^{(0)}(\bm \theta,{\bm \delta}_{(0)}^{\sf neb}(\lambda))$$ and Lemma \ref{lem:horc_k1} establishes the asymptotic optimality of $\hat{\lambda}$ obtained from equation \eqref{eq:lam_k0}. 
In theorem \ref{thm:bayesrisk_k0} below we provide decision theoretic guarantees on the \texttt{NEB} estimator and show that the largest coordinate-wise gap between  ${\bm \delta}^{\sf neb}_{(0)}({\hat{\lambda}})$ and $\bm \delta_{(0)}^{\pi}$ is asymptotically small.

\renewcommand*{\thetheorem}{3B} 
\begin{theorem}
	\label{thm:bayesrisk_k0}
	Under the conditions of Theorem \ref{thm:w_k0}, if {$\lim_{n\to\infty}c_nn^{-1/2}\log^3 n=0$} then, for the Poisson and the Binomial model,
	$$ c_n\Big\|{\bm\delta}_{(0)}^{\sf neb}({\hat{\lambda}})-\bm\delta_{(0)}^{\pi}\Big\|_\infty=o_p(1).$$
	Furthermore, under the same conditions, we have for the Poisson and the Binomial model,
$$\lim_{n\to\infty}\mathbb{P}\Big[\mc{L}_n^{(0)}\Bigl\{\bm \theta,{\bm \delta}_{(0)}^{\sf neb}(\hat{\lambda})\Bigr\}\ge \mc{L}_n^{(0)}(\bm \theta,\bm \delta_{(0)}^{\pi})+c_n^{-1}\epsilon \Big]=0\text{ for any }\epsilon>0.$$
\end{theorem}

\section{Technical Details and Proofs}
\label{sec:proofs}
We will begin this section with some notations and then state two lemmas that will be used in proving the statements discussed in Section \ref{sec:theory}. 

Let $c_0, c_1,\ldots$ denote some generic positive constants which may vary in different statements. Let $\mc{D}_n=\{0,1,2,\ldots,C\log n\}$ and given a random sample $(Y_1,\ldots,Y_n)$ from model \eqref{eq:DLEmodel} denote $B_n$ to be the event $\{\max_{1\le i\le n}Y_i\le C\log n\}$ where $C$ is the constant given by lemma \ref{lem:B} below under assumption $(A2)$. 
\renewcommand*{\thelemma}{A} 
\begin{lemma}
	\label{lem:B}
	Assumption $(A2)$ implies that with probability tending to $1$ as $n\to \infty$, $$\max(Y_1,\ldots,Y_n)\le C\log n$$ where $C>0$ is a constant depending on $\epsilon$.
\end{lemma}
Our next lemma below is a statement on the pointwise Lipschitz stability of the optimal solution $\hat{\bm h}_n^{(k)}(\lambda)$ under perturbations on the parameter $\lambda\in\Lambda$. See, for example, \cite{bonnans2013perturbation} for general results on the stability and sensitivity of parametrized optimization problems.
\renewcommand*{\thelemma}{B} 
\begin{lemma}
	\label{lem:C}
Let $\hat{\bm h}_n^{(k)}(\lambda_0)$ be the solution to problems \eqref{eq:quad_opt} and \eqref{eq:quad_opt_k0}, respectively, for $k\in\{0,1\}$ and for some $\lambda_0\in\Lambda$. Then, under Assumption $(A3)$, there exists a constant $L>0$ such that for any $\lambda \in \Lambda$ the solution $\hat{\bm h}_n^{(k)}(\lambda)$ to problems \eqref{eq:quad_opt} and \eqref{eq:quad_opt_k0} satisfies
$$\Big\|\hat{\bm h}_n^{(k)}(\lambda)-\hat{\bm h}_n^{(k)}(\lambda_0)\Big\|_2\le L|\lambda-\lambda_0|
$$

\end{lemma}

\subsection{Proof of Lemma \ref{lem:bayes}}
\label{sec:proof_lem1}
First note that for any coordinate $i$, the integrated Bayes risk of an estimator $\delta_{(k),i}$ of $\theta_i$ is $\sum_{y_i}\int p(y_i|\theta_i)\ell_n^{(k)}(\theta_i,\delta_{(k),i})dG(\theta_i)$ which is minimized with respect to $\delta_{(k),i}$ if for each $y_i$, $\delta_{(k),i}(y_i)$ is defined as $$\delta_{(k),i}^{\pi}(y_i)=\argmin_{\delta_{(k),i}}\int p(y_i|\theta_i)\ell_n^{(k)}(\theta_i,\delta_{(k),i})dG(\theta_i)$$ However, $\int p(y_i|\theta_i)\ell_n^{(k)}(\theta_i,\delta_{(k),i})dG(\theta_i)$ is a minimum with respect to $\delta_{(k),i}$ when $$\delta_{(k),i}^{\pi}(y_i)=\dfrac{\int p(y_i|\theta_i)\theta_i^{1-k}dG(\theta_i)}{\int p(y_i|\theta_i)\theta_i^{-k}dG(\theta_i)}$$ The result then follows by noting that $p(y_i-k)=\int a_{y_i-k}\theta_i^{y_i-k}/g(\theta_i)dG(\theta_i)$, and $p(y_i+1-k)=\int a_{y_i+1-k}\theta_i^{y_i+1-k}/g(\theta_i)dG(\theta_i)$ for $~y_i=k,k+1,\ldots$.

\subsection{Proof of Theorem \ref{thm:spq}}
\label{sec:proof_thm1}
Define $\tilde{\mathbb{M}}_\lambda(\bm h)=\sum_{i,j\in \mathcal{D}_n}\kappa_{\lambda}[h(i),h(j)](i,j)\mathbb{P}(Y=i)\mathbb{P}(Y=j)$ and re-write $\hat{\mathbb{M}}_{\lambda,n}(\bm h)$ as
$$\hat{\mathbb{M}}_{\lambda,n}(\bm h) = \dfrac{1}{n^2}\sum_{i,j\in \mathcal{D}_n}\kappa_{\lambda}[h(i),h(j)](i,j)\mc{C}_{ij}
,$$
where $\mc{C}_{ij}$ is the number of pairs $(Y_r,Y_s)$ in the sample that has $Y_r=i,Y_s=j$. Now, we have 
\begin{equation}
\label{eq:prf_thm1_1}
\sup_{\lambda\in\Lambda}\Big|\hat{\mathbb{M}}_{\lambda,n}(\bm h)-{\mathbb{M}}_{\lambda}(\bm h)\Big|\le \sup_{\lambda\in\Lambda}\Big|\hat{\mathbb{M}}_{\lambda,n}(\bm h)-\tilde{\mathbb{M}}_{\lambda}(\bm h)\Big|+\sup_{\lambda\in\Lambda}\Big|{\mathbb{M}}_\lambda(\bm h)-\tilde{\mathbb{M}}_{\lambda}(\bm h)\Big|
\end{equation}
Consider the first term on the right hand side of the inequality in equation \eqref{eq:prf_thm1_1} above and with $P_i\coloneqq\mathbb{P}(Y=i)$ note that assumption $A2$ and lemma \ref{lem:B} imply
\begin{eqnarray}
\label{eq:prf_thm1_2}
\mathbb{E}\sup_{\lambda\in\Lambda}\Big|\hat{\mathbb{M}}_{\lambda,n}(\bm h)-\tilde{\mathbb{M}}_{\lambda}(\bm h)\Big|\le \sum_{i,j\in \mathcal{D}_n}\mathbb{E}\Big[\sup_{\lambda\in\Lambda}\Big|\kappa_{\lambda}[h(i),h(j)](i,j)\Big|\Big|\dfrac{\mc{C}_{ij}}{n^2}-P_iP_j\Big|\Big]\Bigl\{1+o(1)\Bigr\}\nonumber\\
\le  \sum_{i,j\in \mathcal{D}_n}\Bigl\{\mathbb{E}\Big[\sup_{\lambda\in\Lambda}\Big|\kappa_{\lambda}[h(i),h(j)](i,j)\Big|\Big]^2\mathbb{E}\Big|\dfrac{\mc{C}_{ij}}{n^2}-P_iP_j\Big|^2\Bigr\}^{1/2}\Bigl\{1+o(1)\Bigr\}
\end{eqnarray}
In equation \eqref{eq:prf_thm1_2} above, $\mathbb{E}|n^{-2}{\mc{C}_{ij}}-P_iP_j|^2$ is $O(1/n)$ and assumption $A1$ together with the compactness of $\Lambda$ and the continuity of $\kappa_{\lambda}[h(i),h(j)](i,j)$ with respect to $\lambda$ imply that $\mathbb{E}[\sup_{\lambda\in\Lambda}|\kappa_{\lambda}[h(i),h(j)](i,j)|]^2<\infty$. Thus $\mathbb{E}\sup_{\lambda\in\Lambda}|\hat{\mathbb{M}}_{\lambda,n}(\bm h)-\tilde{\mathbb{M}}_{\lambda}(\bm h)|$ is $O(\log^2n/n^{1/2})$. Now consider the second term on the right hand side of the inequality in equation \eqref{eq:prf_thm1_1} and note that it is bounded above by the following tail sums 
$$
2\sum_{i\in\mc{D}_n,j\notin \mc{D}_n}\sup_{\lambda\in\Lambda}|\kappa_{\lambda}[h(i),h(j)](i,j)|P_iP_j+\sum_{i,j\notin \mc{D}_n}\sup_{\lambda\in\Lambda}|\kappa_{\lambda}[h(i),h(j)](i,j)|P_iP_j.
$$ 
But from Assumption (A1), $\mathbb{E}_p\sup_{\lambda\in\Lambda}|\kappa_{\lambda}[h(U),h(V)](U,V)|< \infty$ and together with assumption (A2) and proof of lemma \ref{lem:B}, it follows that the terms in the display above are $O(n^{-\nu})$ for some $\nu>1/2$. 

Now fix an $\epsilon>0$ and let $c_n=n^{1/2}/\log^2n$. Since $\mathbb{E}_p\sup_{\lambda\in\Lambda}|\hat{\mathbb{M}}_{\lambda,n}(\bm h)-\tilde{\mathbb{M}}_{\lambda}(\bm h)|$ is $O(\log^2n/n^{1/2})$ there exists a finite constant $M>0$ and an $N_1$ such that $c_n \mathbb{E}\sup_{\lambda\in\Lambda}|\hat{\mathbb{M}}_{\lambda,n}(\bm h)-\tilde{\mathbb{M}}_{\lambda}(\bm h)|\le M$ for all $n\ge N_1$. Moreover since $\sup_{\lambda\in\Lambda}|{\mathbb{M}}_\lambda(\bm h)-\tilde{\mathbb{M}}_{\lambda}(\bm h)|\to 0$ as $n\to\infty$, there exists an $N_2$ such that $\sup_{\lambda\in\Lambda}|{\mathbb{M}}_\lambda(\bm h)-\tilde{\mathbb{M}}_{\lambda}(\bm h)|\le M/c_n$ for all $n\ge N_2$. Thus with $t=4M/\epsilon$, we have $\mathbb{P}(c_n\sup_{\lambda\in\Lambda}|\hat{\mathbb{M}}_{\lambda,n}(\bm h)-{\mathbb{M}}_{\lambda}(\bm h)|>t)<\epsilon$ for all $n\ge \max(N_1,N_2)$ which suffices to prove the desired result.

\subsection{Proofs of Theorems \ref{thm:w} and \ref{thm:w_k0}}
\label{sec:proof_thm2}

We will first prove Theorem \ref{thm:w}. Note that from equation \eqref{eq:h(y)}, 
$$
\Big\|\hat{\bm w}_n^{(1)}(\lambda)-\bm w_p^{(1)}\Big\|_2=\Big\|\hat{\bm h}_n^{(1)}(\lambda)-\bm h_0^{(1)}\Big\|_2. 
$$
Now from assumption $A3$ and for any $\epsilon>0$, there exists a $\delta>0$ such that for any $\lambda\in\Lambda$, $$\mathbb{P}\Big[c_n\Big\|\hat{\bm h}_n^{(1)}(\lambda)-\bm h_0^{(1)}\Big\|_2\ge \epsilon\Big]\le\mathbb{P}\Big[c_n\Bigl\{\mathbb{M}_{\lambda}(\hat{\bm h}_n^{(1)})-\mathbb{M}_\lambda(\bm h_0^{(1)})\Bigr\}\ge \delta\Big].$$ But the right hand side is upper bounded by the sum of $\mathbb{P}\Big[c_n\Bigl\{\mathbb{M}_{\lambda}(\hat{\bm h}_n^{(1)})-\hat{\mathbb{M}}_{\lambda,n}(\hat{\bm h}_n^{(1)})\Bigr\}\ge \delta/3\Big]$, $\mathbb{P}\Big[c_n\Bigl\{\hat{\mathbb{M}}_{\lambda,n}(\hat{\bm h}_n^{(1)})-\hat{\mathbb{M}}_{\lambda,n}({\bm h}_0^{(1)})\Bigr\}\ge \delta/3\Big]$ and $\mathbb{P}\Big[c_n\Bigl\{\hat{\mathbb{M}}_{\lambda,n}({\bm h}_0^{(1)})-{\mathbb{M}}_{\lambda}({\bm h}_0^{(1)})\Bigr\}\ge \delta/3\Big]$. From theorem \ref{thm:spq}, the first and third terms go to zero as $n\to\infty$ while the second term is zero since $\hat{\mathbb{M}}_{\lambda,n}(\hat{\bm h}_n^{(1)})\le\hat{\mathbb{M}}_{\lambda,n}({\bm h}_0^{(1)})$. This proves the statement of theorem \ref{thm:w}. To prove theorem \ref{thm:w_k0} first note that from equation \eqref{eq:h(y)_k0}, $$\Big\|\hat{\bm w}_n^{(0)}(\lambda)-\bm w_p^{(0)}\Big\|_2^2\le \sum_{i=1}^{n}\Bigl\{\dfrac{Y_i+1}{\hat{h}_{n,i}^{(0)}(\lambda)h_{0,i}^{(0)}}\Bigr\}^2\Bigl\{\hat{h}_{n,i}^{(0)}(\lambda)-h_{0,i}^{(0)}\Bigr\}^2.$$ From assumption $A2$ and Lemma \ref{lem:B}, there exists a constant $c_0>0$ such that for large $n$, $\max_{1\le i\le n}(Y_i+1)\le c_0\log n$ with high probability. Moreover for $i=1,\cdots,n$, since $\hat{w}_{n,i}^{(0)}(\lambda)>0$ for every $\lambda\in\Lambda$ and $w_{p,i}^{(0)}>0$, lemma \ref{lem:B} and equation \eqref{eq:h(y)_k0} together imply $\hat{h}_{n,i}^{(0)}(\lambda)+c_0\log n> 0$ and ${h}_{0,i}^{(0)}+c_0\log n>0$. Thus, conditional on the event $\{\max_{1\le i\le n}(Y_i+1)\le c_0\log n\}$ and for any $\epsilon>0$, $$\mathbb{P}\Big[c_n\log n\Big\|\hat{\bm w}_n^{(0)}(\lambda)-\bm w_p^{(0)}\Big\|_2\ge \epsilon\Big]\le \mathbb{P}\Big[c_1c_n\Big\|\hat{\bm h}_n^{(0)}(\lambda)-\bm h_0^{(0)}\Big\|_2\ge \epsilon\Big]$$ for some constant $c_1>0$. The proof of the statement of theorem \ref{thm:w_k0} thus follows from the proof of theorem \ref{thm:w} above and lemma \ref{lem:B}.
\subsection{Proof of Lemma \ref{lem:lam_bin_k1}}
\label{sec:lem3_proof}
We will first prove the two statements of lemma \ref{lem:lam_bin_k1} under the scaled squared error loss. The proof for the squared error loss will follow from similar arguments and we will highlight only the important steps. Throughout the proof, we will denote $d_1\coloneqq \inf_{\lambda\in\Lambda}\inf_{1\le i\le n}(1-\hat{h}_{n,i}^{(1)}(\lambda))>0$ and $d_2\coloneqq\inf_{\lambda\in\Lambda}\inf_{1\le i\le n}\hat{w}_{n,i}^{(0)}(\lambda)>0$.
\\~\\
\textbf{Proof of statement 1 for the scaled squared error loss $(k=1)$}
\\~\\
\noindent Note that by triangle inequality $ \sup_{\lambda\in\Lambda}\Big|{\sf ARE}^{(1,\mc{B})}_n(\lambda,\bm Y)-\mc{R}_n^{(1)}(\bm \theta,{\bm \delta}_{(1)}^{\sf neb}(\lambda))\Big|$ is upper bounded by the following sum 
$$ 
\sup_{\lambda\in\Lambda}\Big|{\sf ARE}^{(1,\mc{B})}_n(\lambda,\bm Y)-\mathbb{E}{\sf ARE}^{(1,\mc{B})}_n(\lambda,\bm Y)\Big|+\sup_{\lambda\in\Lambda}\Big|\mathbb{E}{\sf ARE}^{(1,\mc{B})}_n(\lambda,\bm Y)-\mc{R}_n^{(1)}(\bm \theta,{\bm \delta}_{(1)}^{\sf neb}(\lambda))\Big|.
$$ 
Consider the first term. Using definition \ref{def:are_bin_k1}, this term is upper bounded by
\begin{eqnarray}
\label{eq:prf_lem2_1}
\Big|\dfrac{1}{n}\sum_{i=1}^{n}U_i\Big|&+&\sup_{\lambda\in\Lambda}\Big|\dfrac{2}{n}\sum_{i=1}^{n}\Bigl\{{\delta}_{(1),i}^{\sf neb}(\lambda)-\mathbb{E}{\delta}_{(1),i}^{\sf neb}(\lambda)\Bigr\}\Big|\nonumber\\
&+&\sup_{\lambda\in\Lambda}\Big|\dfrac{1}{n}\sum_{i=1}^{n}\Bigl\{(m-Y_i)\psi_\lambda(Y_i)-\mathbb{E}(m-Y_i)\psi_\lambda(Y_i)\Bigr\}\Big|
\end{eqnarray}
where $U_i=Y_i/(m-Y_i+1)-\theta_i$, $\mathbb{E}U_i=0$ and $\mathbb{E}U_i^2<\infty$ since $Y_i/(m-Y_i+1)\le m<\infty$ for all $i=1,\cdots,n$. So, $\Big|n^{-1}\sum_{i=1}^{n}U_i\Big|=O_p(n^{-1/2})$. 

Now consider the second term in equation \eqref{eq:prf_lem2_1} above and define $Z_n(\lambda)=n^{-1}\sum_{i=1}^{n}V_i(\lambda)$ where $V_i(\lambda)={\delta}_{(1),i}^{\sf neb}(\lambda)-\mathbb{E}{\delta}_{(1),i}^{\sf neb}(\lambda)$. Recall that for every $i$ such that $Y_i>0$, $1-\hat{h}_{n,i}^{(1)}(\lambda)>0$ for all $\lambda\in\Lambda$. Moreover, from definition \ref{def:neb_k1} for the Binomial model, 
$$
{\delta}_{(1),i}^{\sf neb}(\lambda)=[Y_i/(m-Y_i+1)]\left\{1-\hat{h}_{n,i}^{(1)}(\lambda)\right\}^{-1}.
$$ 
Thus, ${\delta}_{(1),i}^{\sf neb}(\lambda)\le m/d_1$ and $|V_i(\lambda)|\le 2m/d_1$ for all $i=1,\ldots,n$. Along with the fact that $\mathbb{E}V_i(\lambda)=0$ Hoeffding's inequality gives, for a fixed $\lambda$ and (for now) arbitrary $r_n>1$
\begin{equation}
\label{eq:hoeff1_lem2}
\mathbb{P}\Bigl\{|Z_n(\lambda)|>\dfrac{r_n}{\sqrt{n}}\Bigr\}\le 2\exp\Bigr\{-\dfrac{r_n^2d_1^2}{8m^2}\Bigl\}
\end{equation}
Next for a perturbation $\lambda'$ of $\lambda$ such that $(\lambda,\lambda')\in\Lambda\coloneqq[\lambda_l,\lambda_u]$, we will bound the increments $\Big|Z_n(\lambda)-Z_n(\lambda')\Big|$. To that effect, note that 
$$
n|Z_n(\lambda)-Z_n(\lambda')|\le \|{\bm \delta}_{(1)}^{\sf neb}(\lambda)-{\bm \delta}_{(1)}^{\sf neb}(\lambda')\|_1+\mathbb{E}\|{\bm \delta}_{(1)}^{\sf neb}(\lambda)-{\bm \delta}_{(1)}^{\sf neb}(\lambda')\|_1
$$
and $d_1^2\|{\bm \delta}_{(1)}^{\sf neb}(\lambda)-{\bm \delta}_{(1)}^{\sf neb}(\lambda')\|_1\le m\|\hat{\bm h}^{(1)}_n(\lambda)-\hat{\bm h}^{(1)}_n(\lambda')\|_1$. Now {from lemma \ref{lem:C} we know that $$\|\hat{\bm h}^{(1)}_n(\lambda)-\hat{\bm h}^{(1)}_n(\lambda')\|_1\le n^{1/2}c^{-1}|\lambda-\lambda'|\sup_{{\bm h}\in N_\delta(\hat{\bm h}_n^{(1)}(\lambda))}\|\nabla^2_{\bm h_n,\lambda}\hat{\mathbb{M}}_{\lambda,n}({\bm h})+o(1)\|_2.$$} Moreover, the Binomial model with $m<\infty$ and assumption $A4$ imply that the supremum in the display above is $O(\log n)$.
Thus so long as $|\lambda-\lambda'|\le\epsilon_n$, 
$$\Big|Z_n(\lambda)-Z_n(\lambda')\Big|\le c_0\epsilon_n\dfrac{\log n}{\sqrt{n}}.$$ 
Now choose $\lambda_j=j\epsilon_n\in\Lambda$ and note that $A_n=\{\sup_{\lambda\in\Lambda}|Z_n(\lambda)|>3r_n/\sqrt{n}\}\subset D_n\cup E_n$ where $$D_n=\{\sup_j |Z_n(\lambda_j)|>r_n/\sqrt{n}\} \text{ and } E_n=\{\sup_j\sup_{|\lambda-\lambda_j|\le\epsilon_n}|Z_n(\lambda)-Z_n(\lambda_j)|>2r_n/\sqrt{n}\}.$$ Choose $\epsilon_n$ so that $\epsilon_nn^{-1/2}\log n=o(n^{-1/2})$ and note that $\mathbb{P}(E_n)\le \exp\{-2r_n^2\}$ and $\mathbb{P}(D_n)\le 2(\lambda_u/\epsilon_n)\exp\{-r_n^2d_1^2/8m^2\}$ from equation \eqref{eq:hoeff1_lem2} and the cardinality of $\lambda_j$. Thus,
$$P(A_n)\le 2(\lambda_u/\epsilon_n)\exp\{-r_n^2d_1^2/8m^2\}+\exp\{-2r_n^2\}
$$
Set $r_n=(s\log n)^{1/2}(m\sqrt{8}/d_1)=O(\sqrt{\log n})$. Then $\mathbb{P}(A_n)\le n^{-s}3(\lambda_u/\epsilon_n)$ and thus the second term in equation \eqref{eq:prf_lem2_1} is $O_p(\sqrt{\log n/n})$.

We will now consider the third term in equation \eqref{eq:prf_lem2_1} and analyze it in a similar manner to the second term of equation \eqref{eq:prf_lem2_1}. Here we will assume that the set $\mathcal{I}_i=\{j:Y_j=Y_i+1\}$ is non-empty for every $i=1,\ldots,n$. Recall from definition \ref{def:are_bin_k1} that $(m-Y_i)\psi_\lambda(Y_i)\le m/d_1^2$ for all $Y_i=0,1,\ldots, m$. Define $Z_n(\lambda)=n^{-1}\sum_{i=1}^{n}V_i(\lambda)$ where $V_i(\lambda)=(m-Y_i)\psi_\lambda(Y_i)-\mathbb{E}(m-Y_i)\psi_\lambda(Y_i)$ with $\mathbb{E}V_i(\lambda)=0$ and $|V_i(\lambda)|\le 2m/d_1^2$. Hoeffding's inequality gives, for a fixed $\lambda$ and (for now) arbitrary $r_n>1$
\begin{equation*}
\label{eq:hoeff2_lem2}
\mathbb{P}\Bigl\{|Z_n(\lambda)|>\dfrac{r_n}{\sqrt{n}}\Bigr\}\le 2\exp\Bigr\{-\dfrac{r_n^2d_1^4}{8m^2}\Bigl\}
\end{equation*}
Moreover, for a perturbation $\lambda'$ of $\lambda$ such that $(\lambda,\lambda')\in\Lambda$, the increments $n|Z_n(\lambda)-Z_n(\lambda')|$ are bounded above by $m\|\bm \psi_\lambda-\bm\psi_{\lambda'}\|_1+m\mathbb{E}\|\bm\psi_{\lambda}-\bm\psi_{\lambda'}\|_1$ where $\bm \psi_{\lambda}=(\psi_{\lambda}(Y_1),\ldots,\psi_{\lambda}(Y_n))$ for any $\lambda\in\Lambda$. Now, note that from definition \ref{def:are_bin_k1} and for the Binomial model, 
$$
|\psi_\lambda(Y_i)-\psi_{\lambda'}(Y_i)|\le 2(m^2/d_1^3)|\hat{h}_{n,j}(\lambda)-\hat{h}_{n,j}(\lambda')|,
$$ 
where $j\in\mathcal{I}_i. $ Therefore
\begin{eqnarray}
|Z_n(\lambda)-Z_n(\lambda')|\le 2\dfrac{m^3}{nd_1^3}\Big[\|\hat{\bm h}^{(1)}_n(\lambda)-\hat{\bm h}^{(1)}_n(\lambda')\|_1+\mathbb{E}\|\hat{\bm h}^{(1)}_n(\lambda)-\hat{\bm h}^{(1)}_n(\lambda')\|_1\Big]\nonumber.
\end{eqnarray}
Now lemma \ref{lem:C} and assumption $A4$ imply that the right hand side of the inequality above is $|\lambda-\lambda'|O(n^{-1/2}\log n )$. So as long as $|\lambda-\lambda'|\le\epsilon_n$, choose $\lambda_j=j\epsilon_n\in\Lambda$. In a manner similar to the second term of equation \eqref{eq:prf_lem2_1} define the events $A_n,~D_n$ and $E_n$, and set $r_n=(s\log n)^{1/2}(m\sqrt{8}/d_1^2)$ to conclude that the third term of equation \eqref{eq:prf_lem2_1} continues to be is $O_p(\sqrt{\log n/n})$ which suffices to prove the statement of the result.
\\~\\
\textbf{Proof of statement 2 for the scaled squared error loss $(k=1)$}
\\~\\
\noindent Note that $\sup_{\lambda\in\Lambda}|{\sf ARE}^{(1,\mc{B})}_n(\lambda,\bm Y)-\mc{L}_n^{(1)}(\bm \theta,{\bm \delta}_{(1)}^{\sf neb}(\lambda))|$ is bounded above by the sum of:
\newline $\sup_{\lambda\in\Lambda}|{\sf ARE}^{(1,\mc{B})}_n(\lambda,\bm Y)-\mc{R}_n^{(1)}(\bm \theta,{\bm \delta}_{(1)}^{\sf neb}(\lambda))|$ and $\sup_{\lambda\in\Lambda}|\mc{L}_n^{(1)}(\bm \theta,{\bm \delta}_{(1)}^{\sf neb}(\lambda))-\mathbb{E}\mc{L}_n^{(1)}(\bm \theta,{\bm \delta}_{(1)}^{\sf neb}(\lambda))|$. The first term is $O_p(\sqrt{\log n/n})$ from statement 1. The second term is bounded above by
\begin{equation}
\label{eq:prf_lem2_2}
\dfrac{2}{n}\sup_{\lambda\in\Lambda}\Big|\sum_{i=1}^{n}\Bigl\{{\delta}_{(1),i}^{\sf neb}(\lambda)-\mathbb{E}{\delta}_{(1),i}^{\sf neb}(\lambda)\Bigr\}\Big|+\dfrac{1}{n}\sup_{\lambda\in\Lambda}\Big|\sum_{i=1}^{n}\theta_i^{-1}\Bigl\{\Big[{\delta}_{(1),i}^{\sf neb}(\lambda)\Big]^2-\mathbb{E}\Big[{\delta}_{(1),i}^{\sf neb}(\lambda)\Big]^2\Bigr\}\Big|
\end{equation}
where the first term in equation \eqref{eq:prf_lem2_2} is $O_p(\sqrt{\log n/n})$ from the proof of statement 1. Now consider the second term in equation \eqref{eq:prf_lem2_2} and define $Z_n(\lambda)=n^{-1}\sum_{i=1}^{n}V_i(\lambda)$ where $\theta_iV_i(\lambda)=[{\delta}_{(1),i}^{\sf neb}(\lambda)]^2-\mathbb{E}[{\delta}_{(1),i}^{\sf neb}(\lambda)]^2$. Recall that from definition \ref{def:neb_k1} for the Binomial model, 
$${\delta}_{(1),i}^{\sf neb}(\lambda)=[Y_i/(m-Y_i+1)]\left\{1-\hat{h}_{n,i}^{(1)}(\lambda)\right\}^{-1}.$$ Thus, $[{\delta}_{(1),i}^{\sf neb}(\lambda)]^2\le m^2/d_1^2$ and $|V_i(\lambda)|\le 2\theta_i^{-1}m^2/d_1^2$ for all $i=1,\ldots,n$. For an arbitrary $r_n>1$ and $\lambda$ fixed, we have from Hoeffding's inequality,
\begin{equation}
\label{eq:hoeff3_lem2}
\mathbb{P}\Bigl\{|Z_n(\lambda)|>\dfrac{r_n}{\sqrt{n}}\Bigr\}\le 2\exp\Bigr\{-\dfrac{r_n^2d_1^4 n}{8m^2\sum_{i=1}^{n}\theta_i^{-1}}\Bigl\}
\end{equation}
Moreover for a perturbation $\lambda'$ of $\lambda$ such that $(\lambda,\lambda')\in\Lambda$,
$$\sum_{i=1}^{n}\theta_i^{-1}\Big|[{\delta}_{(1),i}^{\sf neb}(\lambda)]^2-[{\delta}_{(1),i}^{\sf neb}(\lambda')]^2\Big|\le \dfrac{2m^2}{d_1^2}\Big\|\hat{\bm h}^{(1)}_n(\lambda)-\hat{\bm h}^{(1)}_n(\lambda')\Big\|_2\sum_{i=1}^{n}\theta_i^{-1}.
$$
Thus lemma \ref{lem:C}, assumption $A4$ and the display above together imply that $n|Z_n(\lambda)-Z_n(\lambda')|$ is bounded above by $c_1|\lambda-\lambda'|\log n\sum_{i=1}^{n}\theta_i^{-1}$. Now as long as $|\lambda-\lambda'|\le\epsilon_n$, choose $\epsilon_n$ so that $\epsilon_nn^{-1}\log n\sum_{i=1}^{n}\theta_i^{-1}=o(n^{-1/2})$ and along with equation \eqref{eq:hoeff3_lem2}, follow the steps outlined in the proof of the second term in equation \eqref{eq:prf_lem2_1} to conclude that the second term in equation \eqref{eq:prf_lem2_2} is $O_p(\sqrt{\log n/n})$ from which the desired result follows.
\\~\\
\textbf{Proof of statement 1 for the squared error loss $(k=0)$}
\\~\\
The proof of this statement is very similar to the proof of statement 1 under the scaled squared error loss and therefore we highlight the important steps here. To prove statement 1, we will only look at the term $ \sup_{\lambda\in\Lambda}\Big|{\sf ARE}^{(0,\mc{B})}_n(\lambda,\bm Y)-\mathbb{E}{\sf ARE}^{(0,\mc{B})}_n(\lambda,\bm Y)\Big|$ because under the Binomial model, it can be verified using definition \ref{def:are_bin_k0} that $\mathbb{E}{\sf ARE}^{(0,\mc{B})}_n(\lambda)=\mc{R}_n^{(0)}(\bm \theta,{\bm \delta}_{(0)}^{\sf neb}(\lambda))$. Now note that $ \sup_{\lambda\in\Lambda}\Big|{\sf ARE}^{(0,\mc{B})}_n(\lambda,\bm Y)-\mathbb{E}{\sf ARE}^{(0,\mc{B})}_n(\lambda,\bm Y)\Big|$ is bounded above by
\begin{eqnarray}
\label{eq:prf_lem2_3}
\Big|\dfrac{1}{n}\sum_{i=1}^{n}U_{i}\Big|&+&\sup_{\lambda\in\Lambda}\Big|\dfrac{2}{n}\sum_{i=1}^{n}\Bigl\{Y_i\psi_{\lambda}(Y_i)-\mathbb{E}Y_i\psi_{\lambda}(Y_i)\}\Big|\nonumber\\
&+&\dfrac{1}{n}\sup_{\lambda\in\Lambda}\Big|\sum_{i=1}^{n}\Bigl\{\Big[{\delta}_{(0),i}^{\sf neb}(\lambda)\Big]^2-\mathbb{E}\Big[{\delta}_{(0),i}^{\sf neb}(\lambda)\Big]^2\Bigr\}\Big|
\end{eqnarray}
where $U_i=Y_i(Y_i-1)/[(m-Y_i+2)(m-Y_i+1)]-\theta_i^2$, $\mathbb{E}U_i=0$ and $\mathbb{E}U_i^2<\infty$ since $|U_i|\le m^2<\infty$ for all $i=1,\cdots,n$. So, $\Big|n^{-1}\sum_{i=1}^{n}U_i\Big|=O_p(n^{-1/2})$. For the second term in equation \eqref{eq:prf_lem2_3} note that from definition \ref{def:are_bin_k0}, $Y_i\psi_{\lambda}(Y_i)\le m^2/d_2$ and for a perturbation $\lambda'$ of $\lambda$ such that $(\lambda,\lambda')\in\Lambda$,
$$\sum_{i=1}^{n}Y_i\Big|\psi_{\lambda}(Y_i)-\psi_{\lambda'}(Y_i)\Big|\le c_2|\lambda-\lambda'|\sqrt{n}\log n\{1+o(1)\}
.$$ The last inequality in the display above follows from definition \ref{def:are_bin_k0}, lemma \ref{lem:C} and assumption $A4$. Thus the upper bound on the second term of equation \eqref{eq:prf_lem2_3} and the corresponding upper bound on its increments over $(\lambda,\lambda')\in\Lambda$ suffice to show that this term is $O_p(\sqrt{\log n/n})$. Finally, the third term in equation \eqref{eq:prf_lem2_3} is bounded above by $4m^2/d_2^2$ and $\|{\bm \delta}_{(0)}^{\sf neb}(\lambda)-{\bm \delta}_{(0)}^{\sf neb}(\lambda')\|_1$ is $|\lambda-\lambda'|O(\sqrt{n}\log n)$ for $(\lambda,\lambda')\in\Lambda$ from which the desired follows that $ \sup_{\lambda\in\Lambda}\Big|{\sf ARE}^{(0,\mc{B})}_n(\lambda,\bm Y)-\mathbb{E}{\sf ARE}^{(0,\mc{B})}_n(\lambda,\bm Y)\Big|$ is $O_p(\sqrt{\log n/n})$.
\\~\\
\textbf{Proof of statement 2 for the squared error loss $(k=0)$}
\\~\\
We will only look at the term $\sup_{\lambda\in\Lambda}|\mc{L}_n^{(0)}(\bm \theta,{\bm \delta}_{(0)}^{\sf neb}(\lambda))-\mathbb{E}\mc{L}_n^{(0)}(\bm \theta,{\bm \delta}_{(0)}^{\sf neb}(\lambda))|$ and show that it is $O_p(\sqrt{\log n/n})$. Note that
\begin{equation}
\label{eq:prf_lem2_4}
\dfrac{2}{n}\sup_{\lambda\in\Lambda}\Big|\sum_{i=1}^{n}\theta_i\Bigl\{{\delta}_{(0),i}^{\sf neb}(\lambda)-\mathbb{E}{\delta}_{(0),i}^{\sf neb}(\lambda)\Bigr\}\Big|+\dfrac{1}{n}\sup_{\lambda\in\Lambda}\Big|\sum_{i=1}^{n}\Bigl\{\Big[{\delta}_{(0),i}^{\sf neb}(\lambda)\Big]^2-\mathbb{E}\Big[{\delta}_{(0),i}^{\sf neb}(\lambda)\Big]^2\Bigr\}\Big|
\end{equation}
is an upper bound on $\sup_{\lambda\in\Lambda}|\mc{L}_n^{(0)}(\bm \theta,{\bm \delta}_{(0)}^{\sf neb}(\lambda))-\mathbb{E}\mc{L}_n^{(0)}(\bm \theta,{\bm \delta}_{(0)}^{\sf neb}(\lambda))|$. Analogous to the preceding proofs under lemma \ref{lem:lam_bin_k1}, we will provide the upper bounds and bounds on the increments with respect to perturbations on $\lambda$ for the two terms in equation \eqref{eq:prf_lem2_4}. The rest of the proof will then follow from the proof of statement 1 in lemma \ref{lem:lam_bin_k1}. 

From definition \ref{def:neb_k0}, we know that ${\delta}_{(0),i}^{\sf neb}(\lambda)\le 2m/d_2$. Thus $$
\sum_{i=1}^{n}|\theta_i{\delta}_{(0),i}^{\sf neb}(\lambda)|\le 2(m/d_2)\sum_{i=1}^{n}\theta_i.
$$
Moreover, from Lemma \ref{lem:C} and Assumption (A4), 
$$\sum_{i=1}^{n}\theta_i\Big|{\delta}_{(0),i}^{\sf neb}(\lambda)-{\delta}_{(0),i}^{\sf neb}(\lambda')\Big|\le c_3|\lambda-\lambda'|\sum_{i=1}^{n}\theta_i\log n\{1+o(1)\}
$$ and 
$$\sum_{i=1}^{n}\Big|\Big[{\delta}_{(0),i}^{\sf neb}(\lambda)\Big]^2-[{\delta}_{(0),i}^{\sf neb}(\lambda')\Big]^2\Big|\le c_4|\lambda-\lambda'|\sqrt{n}\log n\{1+o(1)\}
$$
for $(\lambda,\lambda')\in\Lambda$. Hoeffding's inequality and the steps outlined in the proof of statement 1 ($k=1$) of lemma \ref{lem:lam_bin_k1} will then show that both the terms in equation \eqref{eq:prf_lem2_4} are $O_p(\sqrt{\log n/n})$ which proves the desired result.
\subsection{Proof of Lemma \ref{lem:lam_pois_k1}}
\label{sec:lem2_proof}
We will first prove the two statements of lemma \ref{lem:lam_pois_k1} under the scaled squared error loss. The proof for the squared error loss will follow from similar arguments and we will highlight only the important steps. Throughout the proof, we will denote $d_1\coloneqq \inf_{\lambda\in\Lambda}\inf_{1\le i\le n}(1-\hat{h}_{n,i}^{(1)}(\lambda))>0$ and $d_2\coloneqq\inf_{\lambda\in\Lambda}\inf_{1\le i\le n}\hat{w}_{n,i}^{(0)}(\lambda)>0$.
\\~\\
\textbf{Proof of statement 1 for the scaled squared error loss $(k=1)$}
\\~\\
\noindent Note that by triangle inequality $ \sup_{\lambda\in\Lambda}\Big|{\sf ARE}^{(1,\mc{P})}_n(\lambda,\bm Y)-\mc{R}_n^{(1)}(\bm \theta,{\bm \delta}_{(1)}^{\sf neb}(\lambda))\Big|$ is upper bounded by the following sum: 
$$\sup_{\lambda\in\Lambda}\Big|{\sf ARE}^{(1,\mc{P})}_n(\lambda,\bm Y)-\mathbb{E}{\sf ARE}^{(1,\mc{P})}_n(\lambda,\bm Y)\Big|+\sup_{\lambda\in\Lambda}\Big|\mathbb{E}{\sf ARE}^{(1,\mc{P})}_n(\lambda,\bm Y)-\mc{R}_n^{(1)}(\bm \theta,{\bm \delta}_{(1)}^{\sf neb}(\lambda))\Big|.
$$ 
Under the Poisson model and definition \ref{def:are_pois_k1} it can be shown that 
$$
\mathbb{E}\left\{{\sf ARE}^{(1,\mc{P})}_n(\lambda,\bm Y)\right\}= \mc{R}_n^{(1)}(\bm \theta,{\bm \delta}_{(1)}^{\sf neb}(\lambda)). 
$$
So the second term in the display above is zero. Now consider the first term and note that it is bounded above by
\begin{eqnarray}
\label{eq:prf_lem3_1}
\Big|\dfrac{1}{n}\sum_{i=1}^{n}U_i\Big|+\sup_{\lambda\in\Lambda}\Big|\dfrac{2}{n}\sum_{i=1}^{n}\Bigl\{{\delta}_{(1),i}^{\sf neb}(\lambda)-\mathbb{E}{\delta}_{(1),i}^{\sf neb}(\lambda)\Bigr\}\Big|+\sup_{\lambda\in\Lambda}\Big|\dfrac{1}{n}\sum_{i=1}^{n}\Bigl\{\psi_\lambda(Y_i)-\mathbb{E}\psi_\lambda(Y_i)\Bigr\}\Big|
\end{eqnarray}
where $U_i=Y_i-\theta_i$ and $\mathbb{E}U_i=0$. Then using Hoeffding's inequality on $\Big|n^{-1}\sum_{i=1}^{n}U_i\Big|$ along with assumption $A3$ and lemma \ref{lem:B} gives that the first term in equation \eqref{eq:prf_lem3_1} is $O_p(\log^{3/2}n/\sqrt{n})$. Now consider the second term in equation \eqref{eq:prf_lem3_1} and define $V_i(\lambda)={\delta}_{(1),i}^{\sf neb}(\lambda)-\mathbb{E}{\delta}_{(1),i}^{\sf neb}(\lambda)$ with $Z_n(\lambda)=n^{-1}\sum_{i=1}^{n}V_i(\lambda)$. Conditional on the event $B_n$, we have $|V_i(\lambda)|\le 2C\log n/d_1$. Moreover conditional on $B_n$, assumptions $A2-A4$ and lemma \ref{lem:C} give
$$\Big|Z_n(\lambda)-Z_n(\lambda')\Big|\le c_0\epsilon_n\dfrac{\log^4n}{\sqrt{n}}\{1+o(1)\}
$$
whenever $(\lambda,\lambda')\in\Lambda$ and $|\lambda-\lambda'|\le \epsilon_n$. The proof of statement 1 for lemma \ref{lem:lam_bin_k1} and the bounds in the display above along with those developed for $V_i(\lambda)$ establish that the second term in equation \eqref{eq:prf_lem3_1} is $O_p(\log ^{3/2}n /\sqrt{n})$. Now for the third term in equation \eqref{eq:prf_lem3_1}, define $V_i(\lambda)=\psi_{\lambda}(Y_i)-\mathbb{E}\psi_{\lambda}(Y_i)$. From assumption $A3$ and lemma \ref{lem:B}, there exists a constant $C'>0$ such that for large $n$, $\max_{1\le i\le n}(Y_i+1)\le C'\log n$ with high probability which gives, conditional on the event $B_n'=\{\max_{1\le i\le n}(Y_i+1)\le C'\log n\}$, $|V_i(\lambda)|\le c_1\log^2n$. Moreover with $Z_n(\lambda)=n^{-1}\sum_{i=1}^{n}V_i(\lambda)$ and conditional on the event $B_n'$, assumptions $A2-A4$ and lemma \ref{lem:C} give
$$\Big|Z_n(\lambda)-Z_n(\lambda')\Big|\le c_2\epsilon_n\dfrac{\log^5 n}{\sqrt{n}}\{1+o(1)\}
.$$ 
Now we mimic the proof of statement 1 for lemma \ref{lem:lam_bin_k1} to establish that the third term in equation \eqref{eq:prf_lem3_1} is $O_p(\log ^{5/2}n/\sqrt{n})$ which proves the statement of the lemma. 
\\~\\
\textbf{Proof of statement 2 for the scaled squared error loss $(k=1)$}
\\~\\
\noindent For the proof of this statement, we will show that $\sup_{\lambda\in\Lambda}|\mc{L}_n^{(1)}(\bm \theta,{\bm \delta}_{(1)}^{\sf neb}(\lambda))-\mathbb{E}\mc{L}_n^{(1)}(\bm \theta,{\bm \delta}_{(1)}^{\sf neb}(\lambda))|$ is $O_p(\log^{5/2}n/\sqrt{n})$. This term is bounded above by
\begin{equation}
\label{eq:prf_lem3_2}
\dfrac{2}{n}\sup_{\lambda\in\Lambda}\Big|\sum_{i=1}^{n}\Bigl\{{\delta}_{(1),i}^{\sf neb}(\lambda)-\mathbb{E}{\delta}_{(1),i}^{\sf neb}(\lambda)\Bigr\}\Big|+\dfrac{1}{n}\sup_{\lambda\in\Lambda}\Big|\sum_{i=1}^{n}\theta_i^{-1}\Bigl\{\Big[{\delta}_{(1),i}^{\sf neb}(\lambda)\Big]^2-\mathbb{E}\Big[{\delta}_{(1),i}^{\sf neb}(\lambda)\Big]^2\Bigr\}\Big|
\end{equation}
where the first term in equation \eqref{eq:prf_lem3_2} is $O_p(\log^{3/2} n/\sqrt{n})$ from the proof of statement 1. Now consider the second term in equation \eqref{eq:prf_lem3_2} and define $Z_n(\lambda)=n^{-1}\sum_{i=1}^{n}V_i(\lambda)$ where $\theta_iV_i(\lambda)=[{\delta}_{(1),i}^{\sf neb}(\lambda)]^2-\mathbb{E}[{\delta}_{(1),i}^{\sf neb}(\lambda)]^2$. Conditional on the event $B_n$, we have $|V_i(\lambda)|\le 2\theta_i^{-1}(C/d_1)^2\log^2n$. Moreover for a perturbation $\lambda'$ of $\lambda$ such that $(\lambda,\lambda')\in\Lambda$,
$$\sum_{i=1}^{n}\theta_i^{-1}\Big|[{\delta}_{(1),i}^{\sf neb}(\lambda)]^2-[{\delta}_{(1),i}^{\sf neb}(\lambda')]^2\Big|\le \dfrac{2C^2\log^2n}{d_1^3}\Big\|\hat{\bm h}^{(1)}_n(\lambda)-\hat{\bm h}^{(1)}_n(\lambda')\Big\|_2\sum_{i=1}^{n}\theta_i^{-1}
$$
conditional on $B_n$. Thus assumptions $A2-A4$, lemma \ref{lem:C} and the above display together imply that $n|Z_n(\lambda)-Z_n(\lambda')|$ is bounded above by $c_0\epsilon_n\log^5n\sum_{i=1}^{n}\theta_i^{-1}\{1+o(1)\}$ whenever $|\lambda-\lambda'|\le\epsilon_n$. Now we follow the steps outlined in the proof of statement 1 for lemma \ref{lem:lam_pois_k1} to conclude that the second term in equation \eqref{eq:prf_lem3_2} is $O_p(\log^{5/2} n/\sqrt{n})$ from which the desired result follows.
\\~\\
\textbf{Proof of statement 1 for the squared error loss $(k=0)$}
\\~\\
The proof of this statement is very similar to the proof of statement 1 under the scaled squared error loss and therefore we highlight the important steps here. To prove statement 1, we will only show that the term $ \sup_{\lambda\in\Lambda}\Big|{\sf ARE}^{(0,\mc{P})}_n(\lambda,\bm Y)-\mathbb{E}{\sf ARE}^{(0,\mc{P})}_n(\lambda,\bm Y)\Big|$ is $O_p(\log n^{5/2}/\sqrt{n})$ because under the Poisson model, it can be verified using definition \ref{def:are_pois_k0} that $\mathbb{E}{\sf ARE}^{(0,\mc{P})}_n(\lambda)=\mc{R}_n^{(0)}(\bm \theta,{\bm \delta}_{(0)}^{\sf neb}(\lambda))$. Now note that $ \sup_{\lambda\in\Lambda}\Big|{\sf ARE}^{(0,\mc{P})}_n(\lambda,\bm Y)-\mathbb{E}{\sf ARE}^{(0,\mc{P})}_n(\lambda,\bm Y)\Big|$ is bounded above by
\begin{eqnarray}
\label{eq:prf_lem3_3}
\Big|\dfrac{1}{n}\sum_{i=1}^{n}U_{i}\Big|&+&\sup_{\lambda\in\Lambda}\Big|\dfrac{2}{n}\sum_{i=1}^{n}\Bigl\{Y_i\psi_{\lambda}(Y_i)-\mathbb{E}Y_i\psi_{\lambda}(Y_i)\}\Big|\nonumber\\
&+&\dfrac{1}{n}\sup_{\lambda\in\Lambda}\Big|\sum_{i=1}^{n}\Bigl\{\Big[{\delta}_{(0),i}^{\sf neb}(\lambda)\Big]^2-\mathbb{E}\Big[{\delta}_{(0),i}^{\sf neb}(\lambda)\Big]^2\Bigr\}\Big|
\end{eqnarray}
where $U_i=Y_i(Y_i-1)-\theta_i^2$ and $\mathbb{E}U_i=0$. Then using Hoeffding's inequality on $\Big|n^{-1}\sum_{i=1}^{n}U_i\Big|$, along with assumption $A3$ and lemma \ref{lem:B}, gives that the first term in equation \eqref{eq:prf_lem3_3} is $O_p(\log^{5/2}n/\sqrt{n})$. Now for the second term in equation \eqref{eq:prf_lem3_3}, define $V_i(\lambda)=Y_i\psi_{\lambda}(Y_i)-\mathbb{E}Y_i\psi_{\lambda}(Y_i)$. Conditional on the event $B_n$, $|V_i(\lambda)|\le c_0\log^2n$. Moreover with $Z_n(\lambda)=n^{-1}\sum_{i=1}^{n}V_i(\lambda)$ and conditional on the event $B_n$, assumptions $A2-A4$ and lemma \ref{lem:C} give
$$|Z_n(\lambda)-Z_n(\lambda')|\le c_1\epsilon_n{n^{-1/2}\log^5 n}\{1+o(1)\}$$
whenever $|\lambda-\lambda'|\le \epsilon_n$ for $(\lambda,\lambda')\in\Lambda$. The proof of statement 1 ($k=1$ case) for lemma \ref{lem:lam_pois_k1} and the bounds in the display above along with those developed for $V_i(\lambda)$ establish that the second term in equation \eqref{eq:prf_lem3_3} is $O_p(\log ^{5/2}n /\sqrt{n})$. For the third term in equation \eqref{eq:prf_lem3_3}, we proceed in a similar manner and define $V_i(\lambda)=[{\delta}_{(0),i}^{\sf neb}(\lambda)]^2-\mathbb{E}[{\delta}_{(0),i}^{\sf neb}(\lambda)]^2$ and $Z_n(\lambda)=n^{-1}\sum_{i=1}^{n}V_i(\lambda)$. From Assumption (A3) and Lemma \ref{lem:B}, there exists a constant $C'>0$ such that for large $n$, $\max_{1\le i\le n}(Y_i+1)\le C'\log n$ with high probability which gives, conditional on the event $B_n'=\{\max_{1\le i\le n}(Y_i+1)\le C'\log n\}$, (i) $|V_i(\lambda)|\le c_2\log^2n$, and (ii) under Assumptions (A2)-(A4) and Lemma \ref{lem:C}, 
$$
|Z_n(\lambda)-Z_n(\lambda')|\le c_3\epsilon_nn^{-1/2}\log^4n\{1+o(1)\}
$$
whenever $|\lambda-\lambda'|\le \epsilon_n$ for $(\lambda,\lambda')\in\Lambda$. Thus the third term in equation \eqref{eq:prf_lem3_3} is $O_p(\log^{5/2}/\sqrt{n})$ which follows from the proof of statement 1 ($k=1$ case) for lemma \ref{lem:lam_pois_k1} together with the preceding bounds developed for $V_i(\lambda)$ and $|Z_n(\lambda)-Z_n(\lambda')|$, and suffices to prove the desired result.
\\~\\
\textbf{Proof of statement 2 for the squared error loss $(k=0)$}
\\~\\
We will only look at the term $\sup_{\lambda\in\Lambda}|\mc{L}_n^{(0)}(\bm \theta,{\bm \delta}_{(0)}^{\sf neb}(\lambda))-\mathbb{E}\mc{L}_n^{(0)}(\bm \theta,{\bm \delta}_{(0)}^{\sf neb}(\lambda))|$ and show that it is $O_p(\log^{5/2} n/\sqrt{n})$. Note that
\begin{equation}
\label{eq:prf_lem3_4}
\dfrac{2}{n}\sup_{\lambda\in\Lambda}\Big|\sum_{i=1}^{n}\theta_i\Bigl\{{\delta}_{(0),i}^{\sf neb}(\lambda)-\mathbb{E}{\delta}_{(0),i}^{\sf neb}(\lambda)\Bigr\}\Big|+\dfrac{1}{n}\sup_{\lambda\in\Lambda}\Big|\sum_{i=1}^{n}\Bigl\{\Big[{\delta}_{(0),i}^{\sf neb}(\lambda)\Big]^2-\mathbb{E}\Big[{\delta}_{(0),i}^{\sf neb}(\lambda)\Big]^2\Bigr\}\Big|
\end{equation}
is an upper bound on $\sup_{\lambda\in\Lambda}|\mc{L}_n^{(0)}(\bm \theta,{\bm \delta}_{(0)}^{\sf neb}(\lambda))-\mathbb{E}\mc{L}_n^{(0)}(\bm \theta,{\bm \delta}_{(0)}^{\sf neb}(\lambda))|$. Analogous to the preceding proofs under lemma \ref{lem:lam_pois_k1}, we will provide the upper bounds and bounds on the increments with respect to perturbations on $\lambda$ for the two terms in equation \eqref{eq:prf_lem3_4}. The rest of the proof will then follow from the proof of statement 1 in lemma \ref{lem:lam_pois_k1}. 

From definition \ref{def:neb_k0}, we know that under the Poisson model and conditional on the event $B_n'$, ${\delta}_{(0),i}^{\sf neb}(\lambda)\le C'\log n/d_2$. Moreover, from lemma \ref{lem:C}, 
$\sum_{i=1}^{n}\theta_i|{\delta}_{(0),i}^{\sf neb}(\lambda)-{\delta}_{(0),i}^{\sf neb}(\lambda')|$ is bounded above by $c_0|\lambda-\lambda'|\sum_{i=1}^{n}\theta_i\log^3n\{1+o(1)\}$ and 
$\sum_{i=1}^{n}|[{\delta}_{(0),i}^{\sf neb}(\lambda)]^2-[{\delta}_{(0),i}^{\sf neb}(\lambda')]^2|$ is bounded above by $c_1|\lambda-\lambda'|n^{1/2}\log^4n\{1+o(1)\}$
for $(\lambda,\lambda')\in\Lambda$. Hoeffding's inequality and the steps outlined in the proof of statement 1 ($k=1$) of lemma \ref{lem:lam_pois_k1} will then show that the first term in equation \eqref{eq:prf_lem3_4} is $O_p(\log^{3/2} n/\sqrt{n})$and the second term is $O_p(\log^{5/2} n/\sqrt{n})$ which proves the desired result.
\subsection{Proof of Lemma \ref{lem:horc_k1}}
\label{sec:lem4_proof}
The statement of this lemma follows from part (2) of Lemmata \ref{lem:lam_pois_k1} and \ref{lem:lam_bin_k1}. We will prove this lemma for the Poisson case first. Note that for any $\epsilon > 0$ and $k\in\{0,1\}$, the probability $\mathbb{P}\Big[\mc{L}_n^{(k)}(\bm \theta,\bm \delta_{(k)}^{\sf neb}(\hat{\lambda}))\ge \mc{L}_n^{(k)}(\bm \theta,\bm \delta_{(k)}^{\sf or})+c_n^{-1}\epsilon \Big]$ is bounded above by 
$$\mathbb{P}\Big[\mc{L}_n^{(k)}(\bm \theta,\bm \delta_{(k)}^{\sf neb}(\hat{\lambda}))-{\sf ARE}^{(1,\mc{P})}_n(\hat{\lambda},\bm Y)\ge \mc{L}_n^{(k)}(\bm \theta,\bm \delta_{(k)}^{\sf or})-{\sf ARE}^{(1,\mc{P})}_n(\lambda^{\sf orc},\bm Y)+c_n^{-1}\epsilon \Big],
$$ which converges to 0 by  part (2) of Lemma \ref{lem:lam_pois_k1}. For the Binomial case, similar arguments using part (2) of Lemma \ref{lem:lam_bin_k1} suffice.
 
\subsection{Proofs of Theorems \ref{thm:bayesrisk_k1}, \ref{thm:bayesrisk_k0}}
\label{sec:thm3_proof}
We will first prove Theorem \ref{thm:bayesrisk_k1}. Note that $||{\bm\delta}_{(1)}^{\sf neb}({\hat\lambda})-\bm\delta_{(1)}^{\pi}||_\infty^2\le ||{\bm\delta}_{(1)}^{\sf neb}({\hat\lambda})-\bm\delta_{(1)}^{\pi}||_2^2$ and
$$\Big\|{\bm\delta}_{(1)}^{\sf neb}({\hat\lambda})-\bm\delta_{(1)}^{\pi}\Big\|_2^2=\sum_{i=1}^{n}\Big[\dfrac{a_{Y_i-1}/a_{Y_i}}{\hat{w}^{(1)}_{n,i}(\hat{\lambda})w^{(1)}_{p,i}}\Big]^2\Big[\hat{w}^{(1)}_{n,i}(\hat{\lambda})-w^{(1)}_{p,i}\Big]^2
$$ 
Now, $\hat{w}_{n,i}^{(1)}(\lambda)>0$ for every $\lambda\in\Lambda$ and $w_{p,i}^{(1)}>0$. This fact along with assumption $A2$ and lemma \ref{lem:B} imply that there exists a constant $c_0>0$ such that $||{\bm\delta}_{(1)}^{\sf neb}({\hat\lambda})-\bm\delta_{(1)}^{\pi}||_2\le c_0\log n\|\hat{\bm w}_n^{(1)}(\hat\lambda)-\bm w_p^{(1)}\|_2$. The first result thus follows from the above inequality and Theorem \ref{thm:w}. To prove the second part of the theorem, note that $n^{-1}|\mc{L}_n^{(1)}(\bm \theta,\bm \delta_{(1)}^{\pi})-\mc{L}_n^{(1)}(\bm \theta,{\bm\delta}_{(1)}^{\sf neb}({\hat\lambda}))|$ is upper bounded by
\begin{equation}
\label{eq:prf_thm_3a}
\Big|\sum_{i=1}^{n}\theta_i^{-1}\Bigl\{\Big[{\delta}_{(1),i}^{\sf neb}(\hat\lambda)\Big]^2-\Big[{\delta}_{(1),i}^{\pi}\Big]^2\Bigr\}\Big|+2\Big|\sum_{i=1}^{n}\Big[{\delta}_{(1),i}^{\sf neb}(\hat\lambda)-{\delta}_{(1),i}^{\pi}\Big]\Big|
\end{equation}
Now use the fact that $\hat{w}_{n,i}^{(1)}(\lambda)>0$ for every $\lambda\in\Lambda$ and $w_{p,i}^{(1)}>0$ to deduce, from assumption $A2$ and Lemma \ref{lem:B}, that $|{\delta}_{(1),i}^{\sf neb}(\hat\lambda)+{\delta}_{(1),i}^{\pi}|\le c_1\log n$ for some constant $c_1>0$. Using this inequality we can now upper bound the display in equation \eqref{eq:prf_thm_3a} by:
$$\sum_{i=1}^{n}\Big|{\delta}_{(1),i}^{\sf neb}(\hat\lambda)-{\delta}_{(1),i}^{\pi}\Big|\Bigl\{2+\dfrac{c_1\log n}{\theta_i}\Bigr\}\le \Big\|{\bm\delta}_{(1)}^{\sf neb}({\hat\lambda})-\bm\delta_{(1)}^{\pi}\Big\|_\infty\Bigl\{2n+c_1\log n\sum_{i=1}^{n}\theta_i^{-1}\Bigr\}
$$
Thus, 
$$\Big|\mc{L}_n^{(1)}(\bm \theta,\bm \delta_{(1)}^{\pi})-\mc{L}_n^{(1)}(\bm \theta,{\bm\delta}_{(1)}^{\sf neb}({\hat\lambda}))\Big|\le  \Big\|{\bm\delta}_{(1)}^{\sf neb}({\hat\lambda})-\bm\delta_{(1)}^{\pi}\Big\|_\infty\Bigl\{2+c_1\dfrac{\log n}{n}\sum_{i=1}^{n}\theta_i^{-1}\Bigr\}
$$
Finally, the result follows from the above display and the first part of this theorem after noting that $\theta_i>0$ for all $i=1,2,\ldots$.

We will now prove Theorem \ref{thm:bayesrisk_k0}. Using theorem \ref{thm:w_k0}, the first part of theorem \ref{thm:bayesrisk_k0} follows along similar lines as the first part of theorem \ref{thm:bayesrisk_k1}. To prove the second part of the theorem, note that $n^{-1}\Big|\mc{L}_n^{(0)}(\bm \theta,\bm \delta_{(0)}^{\pi})-\mc{L}_n^{(0)}(\bm \theta,{\bm\delta}_{(0)}^{\sf neb}({\hat\lambda}))\Big|$ is upper bounded by
\begin{equation}
\label{eq:prf_thm_3b}
\Big|\sum_{i=1}^{n}\Bigl\{\Big[{\delta}_{(0),i}^{\sf neb}(\hat\lambda)\Big]^2-\Big[{\delta}_{(0),i}^{\pi}\Big]^2\Bigr\}\Big|+2\Big|\sum_{i=1}^{n}\theta_i\Big[{\delta}_{(0),i}^{\sf neb}(\hat\lambda)-{\delta}_{(0),i}^{\pi}\Big]\Big|
\end{equation}
and the display in equation \eqref{eq:prf_thm_3b} is less than or equal to 
$$
\|{\bm\delta}_{(0)}^{\sf neb}({\hat\lambda})-\bm\delta_{(0)}^{\pi}\|_\infty\{\sum_{i=1}^{n}2\theta_i+c_2n\log n\},
$$ 
where we have used the fact that  $\hat{w}_{n,i}^{(0)}(\lambda)>0$ for every $\lambda\in\Lambda$, $w_{p,i}^{(0)}>0$ and along with assumption $A2$ and Lemma \ref{lem:B}, $|{\delta}_{(0),i}^{\sf neb}(\hat\lambda)+{\delta}_{(0),i}^{\pi}|\le c_2\log n$ for some constant $c_2>0$. Thus for $n$ large,
$$\Big|\mc{L}_n^{(0)}(\bm \theta,\bm \delta_{(0)}^{\pi})-\mc{L}_n^{(0)}(\bm \theta,{\bm\delta}_{(0)}^{\sf neb}({\hat\lambda}))\Big|\le c_3\log n \Big\|{\bm\delta}_{(0)}^{\sf neb}({\hat\lambda})-\bm\delta_{(0)}^{\pi}\Big\|_\infty 
$$
from which the desired result follows.
\subsection{Proofs of Lemmata \ref{lem:B} and \ref{lem:C}}
\label{sec:lemmaAB_proof}
\noindent\textbf{{Proof of Lemma \ref{lem:B}}}
\\~\\
First note that from assumption $A4$ if $N(\delta,n)$ denotes the cardinality of the set $\{i:\theta_i\ge \epsilon^{-(1+\delta)}\log n\}$ for some $\delta>0$, then $N(\delta,n)\to 0$ as $n\to\infty$. We will now prove the statement of lemma \ref{lem:B} for the case when $Y_i|\theta_i\stackrel{ind.}{\sim} \texttt{Poi}(\theta_i)$. For distributions with bounded support, like the Binomial model, the lemma follows trivially. 

Under the Poisson model, we have $\mathbb{P}(Y_i\ge \theta_i+t)\le \exp\{-0.5t^2/(\theta_i+t)\}$ for any $t>0$. The above inequality follows from an application of Bennett inequality to the Poisson MGF (see \cite{pollard15}). Now consider $\mathbb{P}(\max_{i=1,\ldots,n}Y_i\le \theta_i+t)$ and note that since $Y_i$ are all independent, this probability is given by $\prod_{i=1}^{n}[1-\exp\{-0.5t^2/(\theta_i+t)\}]$. Take $t=s\log n$ where $s^2/\{s+\epsilon^{-(1+\delta)}\}>4$. Then with $\theta_i\le\epsilon^{-(1+\delta)}\log n$, the above probability is bounded below by $a_n=\{1-n^{-(1+\nu)}\}^n$ for some $\nu>0$. As $n\to \infty$, $a_n\to 1$ which proves the statement of the lemma.
\\~\\
\noindent\textbf{{Proof of Lemma \ref{lem:C}}}
\\~\\
We begin with some remarks on the optimization problems \eqref{eq:quad_opt} and \eqref{eq:quad_opt_k0}. Note that the feasible set $\bm H_n$ in equation \eqref{eq:quad_opt} (and \eqref{eq:quad_opt_k0})  is compact and independent of $\lambda$. Moreover, the optimization problem in definitions \ref{def:neb_k1} and \ref{def:neb_k0} is convex. Consequently, (i) for all $\lambda\in\Lambda$, the optimization takes place in a compact set, and (ii) the optimal solution set corresponding to any $\lambda\in\Lambda$ is a singleton, $\{\hat{\bm h}_n^{(k)}(\lambda)\}$. Now fix an $\epsilon>0$. Then for any $\lambda\in N_\epsilon(\lambda_0)\cap \Lambda$ there exists a $\delta>0$ such that the optimal solution ${\bm h}_n\coloneqq\hat{\bm h}_n^{(k)}(\lambda)\in N_\delta(\hat{\bm h}_n^{(k)}(\lambda_0))$ and $\hat{\mathbb{M}}_{\lambda,n}\{{\bm h}_n\}-\hat{\mathbb{M}}_{\lambda,n}\{\hat{\bm h}_n^{(k)}(\lambda_0)\}\le 0$. Moreover, we can re-write $\hat{\mathbb{M}}_{\lambda_0,n}\{{\bm h}_n\}-\hat{\mathbb{M}}_{\lambda_0,n}\{\hat{\bm h}_n^{(k)}(\lambda_0)\}$ as
\begin{eqnarray}
\hat{\mathbb{M}}_{\lambda_0,n}\{{\bm h}_n\}-\hat{\mathbb{M}}_{\lambda,n}\{{\bm h}_n\}-\hat{\mathbb{M}}_{\lambda_0,n}\{\hat{\bm h}_n^{(k)}(\lambda_0)\}+\hat{\mathbb{M}}_{\lambda,n}\{\hat{\bm h}_n^{(k)}(\lambda_0)\}+\hat{\mathbb{M}}_{\lambda,n}\{{\bm h}_n\}-\hat{\mathbb{M}}_{\lambda,n}\{\hat{\bm h}_n^{(k)}(\lambda_0)\}\nonumber
\end{eqnarray}
The last term in the display above is negative and thus we can upper bound $\hat{\mathbb{M}}_{\lambda_0,n}\{{\bm h}_n\}-\hat{\mathbb{M}}_{\lambda_0,n}\{\hat{\bm h}_n^{(k)}(\lambda_0)\}$ by
$$
\hat{\mathbb{M}}_{\lambda_0,n}\{{\bm h}_n\}-\hat{\mathbb{M}}_{\lambda,n}\{{\bm h}_n\}
-\hat{\mathbb{M}}_{\lambda_0,n}\{\hat{\bm h}_n^{(k)}(\lambda_0)\}+\hat{\mathbb{M}}_{\lambda,n}\{\hat{\bm h}_n^{(k)}(\lambda_0)\}
$$
Now apply the mean value theorem with respect to $\bm h_n$ to the function $\hat{\mathbb{M}}_{\lambda_0,n}\{{\bm h}_n\}-\hat{\mathbb{M}}_{\lambda,n}\{{\bm h}_n\}$ in the display above and notice that $\hat{\mathbb{M}}_{\lambda_0,n}\{{\bm h}_n^{(k)}\}-\hat{\mathbb{M}}_{\lambda_0,n}\{{\bm h}_n^{(k)}(\lambda_0)\}$ is bounded above by
$$\Big[\nabla_{\bm h_n}\Bigl\{\hat{\mathbb{M}}_{\lambda_0,n}(\bar{\bm h}_n)-\hat{\mathbb{M}}_{\lambda,n}(\bar{\bm h}_n)\Bigr\}\Big]^T\Big[\bm h_n-{\bm h}_n^{(k)}(\lambda_0)\Big]
$$ where $\bar{\bm h}_n=\hat{\bm h}_n^{(k)}(\lambda_0)+\tau\{\bm h_n-\hat{\bm h}_n^{(k)}(\lambda_0)\}$ for some $\tau\in(0,1)$ and $\nabla_{\bm h}\hat{\mathbb{M}}_{\lambda,n}({\bm h})$ is the partial derivative of $\hat{\mathbb{M}}_{\lambda,n}({\bm h})$ with respect to $\bm h$. Using $\nabla_{\bm h_n}[\hat{\mathbb{M}}_{\lambda_0,n}({\bm h}_n)-\hat{\mathbb{M}}_{\lambda,n}({\bm h}_n)]=\nabla^2_{\bm h_n,\lambda}\hat{\mathbb{M}}_{\lambda_0,n}({\bm h}_n)(\lambda-\lambda_0)+o(|\lambda-\lambda_0|)$ we get
$$\hat{\mathbb{M}}_{\lambda_0,n}\{{\bm h}_n\}-\hat{\mathbb{M}}_{\lambda_0,n}\{\hat{\bm h}_n^{(k)}(\lambda_0)\}\le \sup_{{\bm h}\in N_\delta(\hat{\bm h}_n^{(k)}(\lambda_0))}\Big[\Big\|\nabla^2_{\bm h_n,\lambda}\hat{\mathbb{M}}_{\lambda_0,n}({\bm h})+o(1)\Big\|_2\Big]\Big|\lambda-\lambda_0\Big|\Big\|{\bm h}_n-\hat{\bm h}_n^{(k)}(\lambda_0)\Big\|_2
$$ Moreover assumption $A3$ implies that
$$\hat{\mathbb{M}}_{\lambda_0,n}\{{\bm h}_n^{(k)}\}-\hat{\mathbb{M}}_{\lambda_0,n}\{\hat{\bm h}_n^{(k)}(\lambda_0)\}\ge c\Big\|{\bm h}_n^{(k)}-\hat{\bm h}_n^{(k)}(\lambda_0)\Big\|_2^2
$$
The desired result thus follows from the above two displays with  $$L=\sup_{{\bm h}\in N_\delta(\hat{\bm h}_n^{(k)}(\lambda_0))}\|\nabla^2_{\bm h_n,\lambda}\hat{\mathbb{M}}_{\lambda_0,n}({\bm h})+o(1)\|_2/c.$$

%
\end{appendices}
\newpage
%
%
%
%
%
%
	\def\spacingset#1{\renewcommand{\baselinestretch}%
	{#1}\small\normalsize} \spacingset{1}
\spacingset{0.5}
\footnotesize
\bibliography{paper-ref}

\begin{thebibliography}{35}
\providecommand{\natexlab}[1]{#1}
\providecommand{\url}[1]{\texttt{#1}}
\expandafter\ifx\csname urlstyle\endcsname\relax
  \providecommand{\doi}[1]{doi: #1}\else
  \providecommand{\doi}{doi: \begingroup \urlstyle{rm}\Url}\fi

\bibitem[Aizer and Doyle~Jr(2015)]{aizer2015juvenile}
Anna Aizer and Joseph~J Doyle~Jr.
\newblock Juvenile incarceration, human capital, and future crime: Evidence
  from randomly assigned judges.
\newblock \emph{The Quarterly Journal of Economics}, 130\penalty0 (2):\penalty0
  759--803, 2015.

\bibitem[Bandari et~al.(2012)Bandari, Asur, and Huberman]{bandari2012pulse}
Roja Bandari, Sitaram Asur, and Bernardo~A Huberman.
\newblock The pulse of news in social media: Forecasting popularity.
\newblock \emph{ICWSM}, 12:\penalty0 26--33, 2012.

\bibitem[Bonnans and Shapiro(2013)]{bonnans2013perturbation}
J~Fr{\'e}d{\'e}ric Bonnans and Alexander Shapiro.
\newblock \emph{Perturbation analysis of optimization problems}.
\newblock Springer Science \& Business Media, 2013.

\bibitem[Brown(2008)]{brown2008season}
Lawrence~D Brown.
\newblock In-season prediction of batting averages: A field test of empirical
  bayes and bayes methodologies.
\newblock \emph{The Annals of Applied Statistics}, pages 113--152, 2008.

\bibitem[Brown and Greenshtein(2009)]{brown2009nonparametric}
Lawrence~D Brown and Eitan Greenshtein.
\newblock Nonparametric empirical bayes and compound decision approaches to
  estimation of a high-dimensional vector of normal means.
\newblock \emph{The Annals of Statistics}, pages 1685--1704, 2009.

\bibitem[Brown et~al.(2013)Brown, Greenshtein, and Ritov]{brown2013poisson}
Lawrence~D Brown, Eitan Greenshtein, and Ya'acov Ritov.
\newblock The poisson compound decision problem revisited.
\newblock \emph{Journal of the American Statistical Association}, 108\penalty0
  (502):\penalty0 741--749, 2013.

\bibitem[Chwialkowski et~al.(2016)Chwialkowski, Strathmann, and
  Gretton]{chwialkowski2016kernel}
Kacper Chwialkowski, Heiko Strathmann, and Arthur Gretton.
\newblock A kernel test of goodness of fit.
\newblock JMLR: Workshop and Conference Proceedings, 2016.

\bibitem[Clevenson and Zidek(1975)]{clevenson1975simultaneous}
M~Lawrence Clevenson and James~V Zidek.
\newblock Simultaneous estimation of the means of independent poisson laws.
\newblock \emph{Journal of the American Statistical Association}, 70\penalty0
  (351a):\penalty0 698--705, 1975.

\bibitem[Damm and Dustmann(2014)]{damm2014does}
Anna~Piil Damm and Christian Dustmann.
\newblock Does growing up in a high crime neighborhood affect youth criminal
  behavior?
\newblock \emph{American Economic Review}, 104\penalty0 (6):\penalty0 1806--32,
  2014.

\bibitem[Efron(2011)]{efron2011tweedie}
Bradley Efron.
\newblock Tweedie's formula and selection bias.
\newblock \emph{Journal of the American Statistical Association}, 106\penalty0
  (496):\penalty0 1602--1614, 2011.

\bibitem[Efron(2012)]{efron2012large}
Bradley Efron.
\newblock \emph{Large-scale inference: empirical Bayes methods for estimation,
  testing, and prediction}, volume~1.
\newblock Cambridge University Press, 2012.

\bibitem[Fourdrinier and Robert(1995)]{fourdrinier1995intrinsic}
Dominique Fourdrinier and Christian~P Robert.
\newblock Intrinsic losses for empirical bayes estimation: A note on normal and
  poisson cases.
\newblock \emph{Statistics \& probability letters}, 23\penalty0 (1):\penalty0
  35--44, 1995.

\bibitem[Fourdrinier et~al.(2018)Fourdrinier, Strawderman, and
  Wells]{fourdrinier2018shrinkage}
Dominique Fourdrinier, William~E Strawderman, and Martin~T Wells.
\newblock \emph{Shrinkage Estimation}.
\newblock Springer, 2018.

\bibitem[Fu et~al.(2017)Fu, Narasimhan, and Boyd]{fu2017cvxr}
Anqi Fu, Balasubramanian Narasimhan, and Stephen Boyd.
\newblock Cvxr: An r package for disciplined convex optimization.
\newblock \emph{arXiv preprint arXiv:1711.07582}, 2017.

\bibitem[Fu et~al.(2018)Fu, James, and Sun]{funonparametric2018}
Luella Fu, Gareth James, and Wenguang Sun.
\newblock Nonparametric empirical bayes estimation on heterogeneous data.
\newblock 2018.

\bibitem[James and Stein(1961)]{james1961estimation}
William James and Charles Stein.
\newblock Estimation with quadratic loss.
\newblock In \emph{Proceedings of the fourth Berkeley symposium on mathematical
  statistics and probability}, volume~1, pages 361--379, 1961.

\bibitem[Jiang et~al.(2009)Jiang, Zhang, et~al.]{jiang2009general}
Wenhua Jiang, Cun-Hui Zhang, et~al.
\newblock General maximum likelihood empirical bayes estimation of normal
  means.
\newblock \emph{The Annals of Statistics}, 37\penalty0 (4):\penalty0
  1647--1684, 2009.

\bibitem[Klenke(2014)]{Klenke2014}
Achim Klenke.
\newblock \emph{Probability Theory: A Comprehensive Course}, pages 331--349.
\newblock Springer London, 2014.

\bibitem[Koenker and Gu(2017)]{koenker2017rebayes}
Roger Koenker and Jiaying Gu.
\newblock Rebayes: An r package for empirical bayes mixture methods.
\newblock \emph{Journal of Statistical Software}, 82\penalty0 (1):\penalty0
  1--26, 2017.

\bibitem[Koenker and Mizera(2014)]{koenker2014convex}
Roger Koenker and Ivan Mizera.
\newblock Convex optimization, shape constraints, compound decisions, and
  empirical bayes rules.
\newblock \emph{Journal of the American Statistical Association}, 109\penalty0
  (506):\penalty0 674--685, 2014.

\bibitem[Koski et~al.(2018)Koski, Bowers, and Costanza]{koski2018state}
Susan~V Koski, David Bowers, and SE~Costanza.
\newblock State and institutional correlates of reported victimization and
  consensual sexual activity in juvenile correctional facilities.
\newblock \emph{Child and Adolescent Social Work Journal}, 35\penalty0
  (3):\penalty0 243--255, 2018.

\bibitem[Liu and Wang(2016)]{liu2016stein}
Qiang Liu and Dilin Wang.
\newblock Stein variational gradient descent: A general purpose bayesian
  inference algorithm.
\newblock In \emph{Advances In Neural Information Processing Systems}, pages
  2378--2386, 2016.

\bibitem[Liu et~al.(2016)Liu, Lee, and Jordan]{liu2016kernelized}
Qiang Liu, Jason~D Lee, and Michael~I Jordan.
\newblock A kernelized stein discrepancy for goodness-of-fit tests.
\newblock In \emph{Proceedings of the International Conference on Machine
  Learning (ICML)}, 2016.

\bibitem[Moniz and Torgo(2018)]{moniz2018multi}
Nuno Moniz and Lu{\'\i}s Torgo.
\newblock Multi-source social feedback of online news feeds.
\newblock \emph{arXiv preprint arXiv:1801.07055}, 2018.

\bibitem[Noack(1950)]{noack1950class}
Albert Noack.
\newblock A class of random variables with discrete distributions.
\newblock \emph{The Annals of Mathematical Statistics}, 21\penalty0
  (1):\penalty0 127--132, 1950.

\bibitem[Oates et~al.(2017)Oates, Girolami, and Chopin]{oates2017control}
Chris~J Oates, Mark Girolami, and Nicolas Chopin.
\newblock Control functionals for monte carlo integration.
\newblock \emph{Journal of the Royal Statistical Society: Series B (Statistical
  Methodology)}, 79\penalty0 (3):\penalty0 695--718, 2017.

\bibitem[Pollard(2015)]{pollard15}
David Pollard.
\newblock A few good inequalities.
\newblock \emph{available at:
  \url{http://www.stat.yale.edu/~pollard/Books/Mini/Basic.pdf}}, 2015.

\bibitem[Robbins(1956)]{robbins1956empirical}
Herbert Robbins.
\newblock An empirical bayes approach to statistics.
\newblock Technical report, COLUMBIA UNIVERSITY New York City United States,
  1956.

\bibitem[Sato and Ken-Iti(1999)]{sato1999levy}
Ken-iti Sato and Sato Ken-Iti.
\newblock \emph{L{\'e}vy processes and infinitely divisible distributions}.
\newblock Cambridge university press, 1999.

\bibitem[Serfling(2009)]{serfling2009approximation}
Robert~J Serfling.
\newblock \emph{Approximation theorems of mathematical statistics}, volume 162.
\newblock John Wiley \& Sons, 2009.

\bibitem[Shmueli et~al.(2005)Shmueli, Minka, Kadane, Borle, and
  Boatwright]{shmueli2005useful}
Galit Shmueli, Thomas~P Minka, Joseph~B Kadane, Sharad Borle, and Peter
  Boatwright.
\newblock A useful distribution for fitting discrete data: revival of the
  conway--maxwell--poisson distribution.
\newblock \emph{Journal of the Royal Statistical Society: Series C (Applied
  Statistics)}, 54\penalty0 (1):\penalty0 127--142, 2005.

\bibitem[{US Department of Justice and Federal Bureau of
  Investigation}(2014)]{us2012uniform}
{US Department of Justice and Federal Bureau of Investigation}.
\newblock {Uniform Crime Reporting Program Data: County-Level Detailed Arrest
  and Offense Data, United States, 2012}.
\newblock \emph{Inter-University Consortium for Political and Social Research
  Ann Arbor, MI}, 2014.
\newblock \doi{https://doi.org/10.3886/ICPSR35019.v1}.

\bibitem[{US Department of Justice and Federal Bureau of
  Investigation}(2017)]{us2014uniform}
{US Department of Justice and Federal Bureau of Investigation}.
\newblock {Uniform Crime Reporting Program Data: County-Level Detailed Arrest
  and Offense Data, United States, 2014}.
\newblock \emph{Inter-University Consortium for Political and Social Research
  Ann Arbor, MI}, 2017.
\newblock \doi{https://doi.org/10.3886/ICPSR36399.v2}.

\bibitem[Weinstein et~al.(2018)Weinstein, Ma, Brown, and
  Zhang]{weinstein2018group}
Asaf Weinstein, Zhuang Ma, Lawrence~D Brown, and Cun-Hui Zhang.
\newblock Group-linear empirical bayes estimates for a heteroscedastic normal
  mean.
\newblock \emph{Journal of the American Statistical Association}, pages 1--13,
  2018.

\bibitem[Xie et~al.(2012)Xie, Kou, and Brown]{xie2012sure}
Xianchao Xie, SC~Kou, and Lawrence~D Brown.
\newblock Sure estimates for a heteroscedastic hierarchical model.
\newblock \emph{Journal of the American Statistical Association}, 107\penalty0
  (500):\penalty0 1465--1479, 2012.

\end{thebibliography}

\end{document}